\documentclass{article}
\usepackage{psfrag}
\usepackage[parfill]{parskip}
\usepackage{float, graphicx}
\usepackage[]{epsfig}
\usepackage{amsmath, amsthm, amssymb,}
\usepackage{epsfig}
\usepackage{verbatim}
\usepackage{multicol}
\usepackage{url}
\usepackage{latexsym}
\usepackage{mathrsfs}
\usepackage[colorlinks, bookmarks=true]{hyperref}
\usepackage{graphicx}
\usepackage{amsmath}
\usepackage{enumerate}
\usepackage[normalem]{ulem}
\usepackage{bm}
\usepackage{dsfont}
\usepackage{stmaryrd}
\usepackage{enumitem}
\usepackage{color}
\usepackage{xcolor}
\usepackage{hyperref}% http://ctan.org/pkg/hyperref
\usepackage[margin=1.1in]{geometry}

\usepackage[noadjust]{cite}

\makeatletter

\def\@settitle{\begin{center}%
    \bfseries
 \normalfont\LARGE\@title
  \end{center}%
}
\def\@setauthors{\begin{center}%
 \normalsize\@author
  \end{center}%
}

\makeatother

\numberwithin{equation}{section}

\def\e{\mathrm{e}}

\def\rr{\mathbb{R}}
\def\ee{\mathbb{E}}

\renewcommand{\cal}{\mathcal}
\newcommand\cA{{\mathcal A}}

\newcommand{\cC}{{\cal C}}
\newcommand{\cD}{{\cal D}}

\newcommand{\cF}{{\cal F}}

\newcommand{\cN}{{\cal N}}

\newcommand{\cX}{{\mathcal X}}

\newcommand{\fb}{{\mathfrak b}}
\newcommand{\fc}{{\mathfrak c}}
\newcommand{\fd}{{\mathfrak d}}

\newcommand{\fo}{{\mathfrak o}}

\newcommand{\bmk}{{\bm{k}}}

\newcommand{\bmv}{{\bm{v}}}

\newcommand{\bmy}{{\bm{y}}}

\newcommand{\rd}{{\rm d}}

\newcommand{\ri}{\mathrm{i}}

\newcommand{\bC}{{\mathbb C}}
\newcommand{\bE}{\mathbb{E}}

\newcommand{\bN}{\mathbb{N}}
\newcommand{\bP}{\mathbb{P}}

\newcommand{\bR}{{\mathbb R}}

\newcommand{\al}{\alpha}

\newcommand{\la}{\lambda}

\DeclareMathOperator{\Tr}{Tr}

\DeclareMathOperator{\supp}{supp}
\DeclareMathOperator{\dist}{dist}

\DeclareMathOperator{\OO}{O}
\DeclareMathOperator{\oo}{o}

\renewcommand{\Re}{\mathop{\mathrm{Re}}}
\renewcommand{\Im}{\mathop{\mathrm{Im}}}

\newcommand{\deq}{\mathrel{\mathop:}=} %define :=
 
\renewcommand{\leq}{\leqslant}
\renewcommand{\geq}{\geqslant}

\newcommand{\td}{\tilde}
\newcommand{\del}{\partial}

\newcommand{\qq}[1]{[\![{#1}]\!]}

\newcommand{\beq}{\begin{equation}}
\newcommand{\eeq}{\end{equation}}
\theoremstyle{plain} %plain, definition, remark
\newtheorem{theorem}{Theorem}[section]
\newtheorem*{theorem*}{Theorem}

\newtheorem*{lemma*}{Lemma}

\newtheorem*{corollary*}{Corollary}
\newtheorem{proposition}[theorem]{Proposition}
\newtheorem*{proposition*}{Proposition}

\newtheorem*{assumption*}{Assumption}
\newtheorem{claim}[theorem]{Claim}

\newtheorem{definition}[theorem]{Definition}
\newtheorem*{definition*}{Definition}
\newtheorem{example}[theorem]{Example}
\newtheorem*{example*}{Example}
\newtheorem{remark}[theorem]{Remark}

\newtheorem*{remark*}{Remark}
\newtheorem*{remarks*}{Remarks}

\def\author#1{\par
    {\centering{\authorfont#1}\par\vspace*{0.05in}}
}

\def\titlefont{\fontsize{13}{15}\bfseries\boldmath\selectfont\centering{}}
\def\authorfont{\fontsize{13}{15}}

\let\affiliationfont\rhfont

\def\address#1{\par
    {\centering{\affiliationfont#1\par}}\par\vspace*{11pt}
}

\def\body{
\setcounter{footnote}{0}
\def\thefootnote{\alph{footnote}}
\def\@makefnmark{{$^{\rm \@thefnmark}$}}
}

\def\title#1{
    \thispagestyle{plain}
    \vspace*{-14pt}
    \vskip 79pt
    {\centering{\titlefont #1\par}}%
    \vskip 1em
}

\def\ev{\boldsymbol{e}}
\def\rhosc{\rho_{\text{sc}}}

\newcommand{\bfi}{{\bf i}}

\begin{document}

\title{Edge Universality of Sparse Random Matrices}
%\title{Extreme Eigenvalues of Sparse Erd{\H o}s-R{\'e}nyi graphs: A Transition from Gaussian to Tracy-Widom}

\vspace{1.2cm}

\noindent \begin{minipage}[c]{0.5\textwidth}
 \author{Jiaoyang Huang}
\address{University of Pennsylvania\\
   E-mail: huangjy@wharton.upenn.edu}
 \end{minipage}
\begin{minipage}[c]{0.5\textwidth}
 \author{Horng-Tzer Yau}
\address{Harvard University \\
   E-mail: htyau@math.harvard.edu}

 \end{minipage}

\begin{abstract}
We consider the statistics of the extreme eigenvalues of sparse random matrices, a class of random matrices that includes the normalized adjacency matrices of the Erd{\H o}s-R{\'e}nyi graph $G(N,p)$.
Recently, it was shown in \cite{lee2021higher}, up to an explicit random shift, the optimal rigidity of extreme eigenvalues holds, provided the averaged degree grows with the size of the graph, $pN>N^\varepsilon$.  We prove in the same regime, 
(i) Optimal rigidity holds for all eigenvalues with respect to an explicit random measure. 
(ii) Up to an explicit random shift, the  fluctuations of the extreme eigenvalues are given the Tracy-Widom distribution.

\end{abstract}

\section{Introduction}\label{s:intro}

In this work we study the statistics of eigenvalues at the edge of  the spectrum of sparse random matrices. A natural example is the  adjacency matrix of the  Erd\H{o}s-R\'enyi graph $G(N,p)$, which is the random undirected graph on $N$ vertices in which each edge appears independently with probability $p$. Introduced in \cite{MR0108839,MR0120167}, the Erd\H{o}s-R\'enyi graph $G(N,p)$ has numerous applications in  graph theory, network theory, mathematical physics and  combinatorics. For further information, we refer the reader to the monographs \cite{MR1782847,MR1864966}.
Many interesting properties of graphs are revealed by the eigenvalues and eigenvectors of their adjacency matrices. 

The adjacency matrices of Erd\H{o}s-R\'enyi graphs have typically $pN$ nonzero entries in each column and are sparse if $p \ll 1$. %As the matrix entries take values $0$ or $1$, the mean of the entries is not  $0$.  
%Up to the fact that the mean of the matrix entries is non-zero, 
When $p$ is of constant order, the Erd\H{o}s-R\'enyi matrix is essentially a Wigner matrix (up to a non-zero mean of the matrix entries).  When $p \to 0$ as $N \to \infty$, the law of the matrix entries is highly concentrated at $0$, and the Erd\H{o}s-R\'enyi matrix can be viewed as a singular Wigner matrix.
%The adjacency matrices of Erd\H{o}s-R\'enyi graphs can be thought of as singular Wigner matrices.  The singularity is due to the fact that the distribution of the elements if highly concentrated around $0$.  
The singular nature of this ensemble can be expressed by the fact that the $k$-th moment of a matrix entry (in the scaling that the bulk of the eigenvalues lie in an interval of order $1$) decays like ($k \geq 2$)
\begin{align} \label{eqn:decaymom}
N^{-1} (p N)^{- (k-2)/2}.
\end{align}
When $p \ll 1$, this decay in $k$ is much slower than the $N^{-k/2}$ case of Wigner matrices, and is the main source of difficulties in studying sparse ensembles with random matrix methods.

The class of random matrices whose moments decay like \eqref{eqn:decaymom}  were introduced in the works \cite{MR3098073,MR2964770} as a natural generalization of the sparse Erd\H{o}s-R\'enyi graph and encompass many other sparse ensembles.  This is the class we study in this work.

The global statistics of the eigenvalues of the Erd\H{o}s-R\'enyi graph are well understood. The empirical eigenvalue distribution converges to the semi-circle distribution provided $p\gg 1/N$, which follows from Wigner's original proof. It was proven in \cite{MR2384414,BFK1, benaych2019largest,latala2018dimension}, the spectral norm of normalized adjacency matrices of Erd\H{o}s-R\'enyi graphs converges to $2$ if $p\gg \log N/N$. Some of their results also extend to nonhomogeneous Erd\H{o}s-R\'enyi graphs. Later, sharp transition was identified in \cite{alt2021extremal, tikhomirov2021outliers}. It was proved that there exists a critical $b_*\approx 2.59$, if $p>b_* \log N/N$ then the spectral norm converges to $2$,  and if $p<b_* \log N/N$, then the extreme eigenvalues are determined by the largest degrees. For the global eigenvalue fluctuations, it was proven in \cite{MR2953145} that the linear statistics, after normalizing by $p^{1/2}$, converge to a Gaussian random variable.

Bulk universality of sparse random matrices (the statement that local statistics inside the bulk are asymptotically the same as those of Gaussian matrices)   was proven in \cite{MR2964770, MR3429490} for $p\geq N^{\varepsilon-1}$ for any $\varepsilon>0$. Universality for the edge statistics of sparse random matrices (the statement that the distribution of the extreme eigenvalues converge to the Tracy-Widon law)  was more intricate. Edge universality for sparse random matrices was proven first in the regime $p\gg N^{-2/3}$ in  \cite{MR2964770, Lee2016}.  Later, it was observed in \cite{huang2020transition} there is a transition in the behavior at $p=N^{-2/3}$. More precisely, it was proved for $N^{-2/9}\ll p\ll N^{-2/3}$, the extreme eigenvalues behave like 
\[
\cX + N^{-2/3}  \xi,  \quad \cX \sim  \frac{\text{Gaussian} }{N\sqrt p},  \quad  \xi \sim \mbox {Tracy-Widon law},
\]
where $\cX$ is related to  the total number of edges. 

%When $N^{-2/9}\ll p\ll N^{-2/3}$, the extreme eigenvalues have Gaussian fluctuation of order $(p^{1/2} N)^{-1}$ instead of Tracy-Widom fluctuation.  
In the regime $N^{-2/9}\ll p\ll N^{-2/3}$, the Gaussian term dominates and the leading behavior changes from the Tracy-Widon law to Gaussian at 
 $p=N^{-2/3}$.  This phenomenon that the leading behavior is Gaussian was  extended down to the optimal scale $p\geq N^{\varepsilon-1}$ in \cite{he2021fluctuations}. For the sparser regime when $(\log \log N)^4/N\ll p< b_*\log N/N$, it was proven in \cite{alt2023poisson},  the eigenvalues near the spectral edges form asymptotically a Poisson point process.
There is however a nature question that  \emph{what is the next order fluctuation term, and can we recover edge universality  by subtracting all these higher order fluctuations}?

This question was partially settled  in \cite{lee2021higher} %by Jaehun Lee, 
where it was proved that, for  $N^{\varepsilon-1}\leq p\ll N^{-2/3}$,  there exists a sequence of explicit random correction terms, which capture higher (sub-leading) order fluctuations of extreme eigenvalues  %The leading order is given by the fluctuation of the total number of edges. 
and after subtracting these explicit correction terms, the optimal rigidity of extreme eigenvalues holds. It was also explicitly 
conjectured in \cite{lee2021higher} that up to this explicit random shift, the  fluctuations of the extreme eigenvalues are given by the Tracy-Widom distribution. One main result of this paper proves this edge universality conjecture. As a consequence, the gaps between extreme eigenvalues of sparse random matrices with $p\geq N^{\varepsilon-1}$ are given asymptotically by the gaps of Airy point process.

%
%Recently, it was shown in \cite{lee2021higher} %by Jaehun Lee, 
%there exists a sequence of explicit random correction terms, which capture higher (sub-leading) order fluctuations of extreme eigenvalues when $N^{\varepsilon-1}\leq p\ll N^{-2/3}$. The leading order is given by the fluctuation of the total number of edges. It was also proven in  \cite{lee2021higher} after subtracting these explicit correction terms, the optimal rigidity of extreme eigenvalues holds, and conjectured \cite[Conjecture 2.12]{lee2021higher} that up to this explicit random shift, the  fluctuations of the extreme eigenvalues are given the Tracy-Widom distribution. One main result of this paper proves this edge universality conjecture. As a consequence, the gaps between extreme eigenvalues of sparse random matrices with $p\geq N^{\varepsilon-1}$ are given asymptotically by the gaps of Airy point process.

Universality for the edge statistics of Wigner matrices was first established by the moment method \cite{MR1727234} under certain symmetry assumptions on the distribution of the matrix elements. The moment method was further developed in \cite{MR2475670,MR2647136} and \cite{MR2726110}. A different approach to edge universality for Wigner matrices  based on the direct comparison with corresponding Gaussian ensembles was developed in \cite{MR2669449,MR2871147, tao2010random}.  

For sparse random matrices, because of those random correction terms to the extreme eigenvalues, naive moment methods or direct comparison with Gaussian ensembles fail. 
Our proof of edge universality utilizes a three-step dynamical approach, which was originally developed  to prove the bulk universality of Wigner matrices in a series of papers \cite{MR2537522,MR2481753,MR3068390,MR2981427,MR2919197,MR2871147,MR2810797,kevin3,landon2019fixed,MR3687212,MR3541852,ajanki2}.  This strategy is as follows: i)
Establish a local semicircle law controlling the number of eigenvalues in windows of size $N^{\delta-1}$, where $\delta>0$ is arbitrarily small.
ii) Analyze the local ergodicity of Dyson Brownian motion to obtain universality after adding a small Gaussian noise to the ensemble.
iii) A density argument comparing a general matrix to one with a small Gaussian component.

For edge universality, in the first step, a local semicircle law is not enough. For edge universality to be true, it is necessary that the extreme eigenvalues, up to an explicit (random) shift, fluctuate on scale $\OO(N^{-2/3})$. Such optimal rigidity estimate for extreme eigenvalues was obtained recently in \cite{lee2021higher}. The proof is based on first constructing a higher order self-consistent equation for the Stieltjes transform of the empirical eigenvalue distributions, then computing the moments of the self-consistent equation by a recursive moment estimate. The rigidity estimates follow from a careful analysis of the recursive moment estimate. 

This approach was first introduced in \cite{Lee2016}, where a deterministic higher order self-consistent equation was constructed and used to prove the optimal edge rigidity (with respect to a deterministically shifted edge) provided $p\gg N^{-2/3}$. Later, a higher order self-consistent equation with one random correction term was constructed in \cite{huang2020transition}, and used to prove the optimal edge rigidity (with respect to a randomly shifted edge)  provided $p\gg N^{-2/9}$.
By including sufficiently many random correction terms, \cite{lee2021higher} proves an optimal edge rigidity down to the optimal scale $p\gg N^{\varepsilon-1}$, and gives a full description of the randomly shifted edge.
We revisit the recursive moment estimate for the random self-consistent equation introduced in \cite{lee2021higher}. By exploring a splitting phenomenon in the expansion, see Proposition \ref{e:one-off}, we improve the error in the recursive moment estimate, and obtain an optimal rigidity for all eigenvalues with respect to an explicit random measure. This is also  crucial for the third step when we do the comparison.

For the second step, edge universality for ensembles with a small Gaussian noise was established in \cite{landon2017edge}.  This work proves for wide classes of initial data, the edge statistics of Dyson Brownian motion coincides with  Gaussian matrices. Moreover \cite{landon2017edge} finds the optimal time to equilibrium $t \sim N^{-1/3}$ for sufficiently regular initial data.

In order to complete the three-step strategy, we need to compare sparse ensembles to Gaussian divisible ensembles, which is a sparse ensemble with a small Gaussian component. Similarly to \cite{Lee2016, huang2020transition}, we can interpolate them by considering the Dyson matrix flow. The main challenge is to keep track of the change of the randomly shifted edge along the interpolation, and show the change of the Stieltjes transform over time is offset by the shift of the edge. The error term in the change of the Stieltjes transform around the randomly shifted edge has a $(pN)^{-1/2}$ expansion. In \cite{Lee2016, huang2020transition}, it was directly checked that the expansion vanishes up to the third order, with an error $\OO((pN)^{-3/2})$. A general principle of the cancellations up to arbitrary order is needed in order to solve the general case $p\geq N^{\varepsilon-1}$. In this paper, we prove a version of such general cancellation principle. First we show the Stieljes transform of sparse ensembles with a small Gaussian component satisfies a modified self-consistent equation, which can be used to precisely characterize the change of the randomly shifted edge along the interpolation.
Next, we prove the change of the Stieltjes transform around the randomly shifted edge (up to negligible error) is given by the derivative of the modified self-consistent equation. Then a similar argument as in the first step, can be used to show its expectation is small up to arbitrary order.

It is worth to comparing the results with random $d$-regular graphs. Bulk universality of random $d$-regular graphs was proven in \cite{bauerschmidt2017bulk} for $d\geq N^{\varepsilon}$ for any $\varepsilon>0$. 
For edge statistics, those Gaussian fluctuations from degree fluctuations in Erd{\H o}s--R{\'e}nyi  are absent in regular graphs. The eigenvalues of random regular graphs are more rigid than those of Erd{\H o}s--R{\'e}nyi graphs of the same average degree. We do not expect any shift of the spectral edge. It was proven that the law of the second largest eigenvalue (after shifting by $2$ and proper normalization) converges to the Tracy-Widom distribution, for $N^{-2/9}\ll d\ll N^{1/3}$ in \cite{bauerschmidt2020edge}; for $d\gg N^{2/3}$ in \cite{he2022spectral}; and
for $1\ll d\ll N^{1/3}$ in \cite{huang2023edge}. 
The edge universality was conjectured in \cite{miller2008distribution} to be true down to $d=3$. For fixed degree $d$, even to show the concentration of extreme eigenvalues around the spectral edges $\pm 2$ requires significant work.  It was first conjectured  by Alon \cite{alon1986eigenvalues}  and proven later in  \cite{friedman2008proof, bordenave2015new, huang2021spectrum}.

\noindent\textbf{Organization.} We define the model and present the main results in the rest of Section \ref{s:intro}. 
In Section \ref{s:rig}, we prove the optimal rigidity estimates for all eigenvalues with respect to an explicit random measure in the regime $pN\geq N^{\varepsilon}$. In Section \ref{sec:ht}, we recall the results from \cite{landon2017edge, adhikari2020dyson} for  edge universality of Gaussian divisible ensembles. 
In Section \ref{s:comparison} we analyze the Stieltjes transform to compare a sparse ensemble to a Gaussian divisible ensemble and establish our results about Tracy-Widom fluctuations.

\noindent\textbf{Notations.} We use $C$ to represent large universal constants, and $c$ small universal constants, which may be different from line by line. Let $Y\geq 0$. %We write that $X=\OO_{\prec}(Y)$ or $X\prec Y$, if there exists an \emph{arbitrarily} small exponent $\fc>0$ such that $X\leq  N^{\fc}Y$ for $N\geq N(\fc)$ large enough. 
%We write $X\ll Y$ or $Y\gg X$, if there exists a small exponent $\fc>0$ such that $X\leq  N^{-\fc}Y$ for $N\geq N(\fc)$ large enough. 
We write $X\lesssim Y$, $Y\gtrsim X$ or $X=\OO(Y)$ if there exists a constant $C>0$, such that $X\leq CY$. 
We write $X\asymp Y$ or $X=\Omega(Y)$ if there exists a constant $C>0$ such that $Y/C\leq X\leq Y/C$.
We use $\bC_+$ to represent the upper half plane. We denote $\qq{a,b}=\{a, a+1,\cdots, b\}$. For any index set $\bm m=\{m_1, m_2,\cdots, m_r\}$,
we write
\begin{align*}
\sum^*_{\bm m}=\sum_{m_1,m_2,\cdots, m_r\in \qq{1,N}\atop \text{distinct}}, \quad \sum_{\bm m}=\sum_{m_1,m_2,\cdots, m_r\in \qq{1,N}}.
\end{align*}
For two index set $\bm m$ and $\bm v$, We denote $\bmv\cup \bm m$ as $\bmv\bm m$. If $\bmv=\{i\}$ has only one element, we  simply write $\{i\}\cup \bm m$ as $i\bm m$. 

\noindent\textbf{Acknowledgements }
The research of J.H. is supported by the Simons Foundation as a Junior Fellow at
the Simons Society of Fellows, and NSF grant DMS-2054835. The research of H-T.Y. is supported by NSF grants DMS-1855509, DMS-2153335 and a Simons Investigator award.

\subsection{Sparse random matrices}

In this section we introduce the class of sparse random matrices that we consider.  This class was introduced in \cite{MR3098073,MR2964770} and we repeat the discussion appearing there.

The Erd\H{o}s-R\'enyi graph is the undirected random graph in which each edge appears with probability $p$. % is the symmetric matrix with independent elements that are $1$ with probability $p$ and $0$ with probability $1-p$. 
 It is notationally convenient to replace the parameter $p$ with $q$ defined through
\begin{align*}
pN = q^2,\quad q=\sqrt{pN}.%\frac{q^2}{N}.
\end{align*}
We allow $q$ to depend on $N$.   We denote by $A$ the adjacency matrix of the Erd\H{o}s-R\'enyi graph.  $A$ is an $N\times N$ symmetric matrix whose entries $A_{ij}$ above the main diagonal are independent and distributed according to
\begin{align*}
A_{ij} =  \begin{cases} 1 & \mbox{ with probability } q^2/N \\ 0 & \mbox{ with probability } 1 - q^2/N \end{cases}.
\end{align*}

We extract the mean of each entry and rescale the matrix so that the limiting eigenvalue distribution is roughly supported on $[-2,2]$.  We introduce the matrix $H$ by %We define the \emph{normalized} adjacency matrix $H$ as
\begin{align*}
H\deq \frac{A-q^2\vert \ev \rangle \langle \ev \vert
}{q\sqrt{1-q^2/N}},
\end{align*}
where $\bm e$ is the unit vector
\begin{align*}
\ev=  ( 1, ..., 1 )^T /\sqrt{N}.
\end{align*}
%
%For the Erd\H{o}s-R\'enyi matrix we define $A$ to be the $N\times N$ symmetric matrix whose entries $a_{ij}$ are independent up $a_{ij} = a_{ji}$ and each element is distributed according to
%\begin{align}
%a_{ij} = \frac{\gamma}{q} \begin{cases} 1 & \mbox{ with probability } \frac{q^2}{N} \\ 0 & \mbox{ with probability } 1 - \frac{q^2}{N} \end{cases}.
%\end{align}
%Here we have defined
%\begin{align}
%\gamma := \left( 1- \frac{q^2}{N} \right)^{-1/2}.
%\end{align}
%
%We further extract the mean of each entry and write
%\begin{align}
%A = H + \gamma q \vert \ev \rangle \langle \ev \vert
%\end{align}
%where $\ev$ is the unit vector
%\begin{align}
%\ev := \frac{1}{ \sqrt{N}} ( 1, ..., 1 )^T.
%\end{align}
%Note that the matrix elements of $H$ are centered. 
It is easy to check that the matrix elements of $H$ (in the upper half triangle) have mean zero $\bE[h_{ij}]=0$, variance $\bE[h_{ij}^2]=1/N$, and satisfy the moment bounds 
\begin{align*}
\bE [  h_{ij}^k ] = \frac{1}{ Nq^{k-2} }\left[\left(1-\frac{q^2}{N}\right)^{-k/2+1}\left(\left(1-\frac{q^2}{N}\right)^{k-1}+(-1)^k\left(\frac{q^2}{N}\right)^{k-1}\right)\right]=\frac{\Omega(1)}{Nq^{k-2}},
\end{align*}
for $k\geq 2$.  This motivates the following definition. %We are prompted to make the following definitions of sparse random matrices.

\begin{definition}[Sparse random matrices] {\label{asup}}
We assume that   $H=(h_{ij})$ is an $N\times N$ random matrix whose entries are real and independent up to the symmetry constraint $h_{ij}=h_{ji}$. We further assume that $(h_{ij})$ satisfies $\bE[h_{ij}]=0$, $\bE[h_{ij}^2]=1/N$ and  that for any $k\geq 2$, the $k$-th cumulant of $h_{ij}$ is given by 
\begin{align}\label{e:asup}
\frac{\cal (k-1)!\cC_k}{Nq^{k-2}},
\end{align}
where $q=q(N)$ is the sparsity parameter, such that $0< q\lesssim \sqrt{N}$.  For $\cC_k$ (which may depend on $q, N$) we make the following assumptions:
\begin{enumerate}[label={\normalfont(\arabic*)}]
\item $| \cC_k | \leq C_k$ for some constant $C_k>0$. 
\item $ \cC_4 \geq c$ \label{item:4clower}
\end{enumerate}% and $\cC_k=\cC_k(q,N)\asymp1$.
\end{definition}
\begin{remark} By the defintion, $\cC_2=1$. The lower bound, $ \cC_4 \geq c$, ensures that the scaling by $q$ for the ensemble $H$ is ``correct.''
\end{remark}

%{\cor
%\begin{definition}[Sparse noncentered random matrices] \label{def:snrm} $A$ is a sparse noncentered random matrix with mean $f$ if it is of the form
%\begin{align}
%A = H + f \vert \ev \rangle \langle \ev \vert
%\end{align}
%where $H$ is a sparse centered random matrix as in Assumption \ref{asup} and $f$ is a deterministic number satisfying
%\begin{align}
%1 <  f \leq N^{1/2}.
%\end{align}
%\end{definition}
%}

We denote the eigenvalues of $H$ by
$
\lambda_1\geq \lambda_2\geq\cdots\geq \lambda_N,
$
corresponding eigenvectors $u_1, u_2,\cdots, u_N$, and the Green's function of $H$ by
\begin{align*}
G(z):=(H-z)^{-1}=\sum_{\al=1}^N\frac{u_\al u_\al^*}{\la_\al-z}.
\end{align*}
The Stieltjes transform of the empirical eigenvalue distribution is denoted by
\begin{align*}
m_N(z):=\frac{1}{N}\Tr G(z)=\frac{1}{N}\sum_{\al=1}^N\frac{1}{\lambda_\al-z}.
\end{align*}
%
%We fix a large constant $\fb>0$ and define the spectral domain,
%\begin{align}\begin{split}\label{e:domain}
%%\cal D_*&=\{z=E+\ri\eta: \;-\fb\leq E\leq \fb, \; 0< \eta\leq \fb\}, \\
%\cal D&=\{z=E+\ri\eta: \;-\fb\leq E\leq \fb, \; 1/N\ll \eta\leq \fb\}. 
%\end{split}\end{align}
%
\begin{definition}[overwhelming probabiltiy]
We say an event $\Omega$ holds with overwhelming probability, if for any $D>0$, $\bP(\Omega)\geq 1-N^{-D}$ for $N\geq N(D)$ large enough. 
\end{definition}

\begin{definition}[Stochastic dominant]
For $N$-dependent random (or deterministic) variables $A$ and $B$, we say $B$ stochastically dominate $A$, if for any $\varepsilon>0$ and $D>0$, then
\begin{equation}
 \bP(A \geq N^\varepsilon B) \leq N^{-D},
\end{equation}
for $N\geq N(\varepsilon, D)$ large enough, and we write   $A \prec B$ or $A=\OO_\prec(B)$.
\end{definition}

%Moreover, for deterministic $N$-dependent quantities $A$, $B$, we write  $A \ll B$, 
%\begin{equation}
% \quad \text{if} \quad A =  O_\varepsilon(N^{-\varepsilon} B) \text{ for some $\varepsilon>0$.}
%\end{equation}

\subsection{Main Results}

We first recall the edge rigidity estimates for the sparse random matrices from \cite[Theorem 2.10]{lee2021higher}.
\begin{theorem}{\normalfont (\cite[Theorem 2.10]{lee2021higher})}\label{t:edgerigidity}
Let $H$ be as in Definition \ref{asup} with  $N^{\varepsilon}\leq q\lesssim N^{1/2}$ and eigenvalues given by $\la_1\geq \la_2\geq \cdots\geq \lambda_N$.  There exists an explicit random measure $\tilde \rho$ supported on $[-{\tilde L}, {\tilde L}]$ (depending on certain averaged quantities of $h_{ij}$ as defined in Proposition \ref{p:minfty}). We have for any fixed index $k\geq 1$,
\begin{align*}
|\la_k-{\tilde L}|\prec \frac{1}{N^{2/3}}.
\end{align*}
Analogous results hold for the smallest eigenvalues.
\end{theorem}

In the following Theorem, we improve the optimal edge rigidity results from \cite[Theorem 2.10]{lee2021higher} to also include the bulk eigenvalues. The proof follows that of  \cite[Theorem 2.10]{lee2021higher} with some modifications. We  give the proof in Section \ref{s:rigidity}.

\begin{theorem}\label{t:rigidity}
Let $H$ be as in Definition \ref{asup} with  $N^{\varepsilon}\leq q\lesssim N^{1/2}$ and eigenvalues given by $\la_1\geq \la_2\geq \cdots\geq \lambda_N$.  There exists an explicit random measure $\tilde \rho$ supported on $[-{\tilde L}, {\tilde L}]$ (depending on certain averaged quantities of $h_{ij}$ as defined in Proposition \ref{p:minfty}).
We denote the classical eigenvalue locations of $\tilde \rho$ as $\gamma_1>\gamma_2>\cdots>\gamma_N$,
\begin{align*}
    \frac{k-1/2}{N}=\int_{\gamma_k}^{\tilde L} \tilde \rho(x)\rd x,\quad 1\leq k\leq N. 
\end{align*}
Then  we have the following optimal rigidity estimates
 \begin{align}\label{e:largeeig}
 |\lambda_k-\gamma_k|\prec \frac{1}{N^{2/3}\min\{k, N-k+1\}^{1/3}},\quad 1\leq k\leq N.
 \end{align}
\end{theorem}

The following theorem concerns the edge fluctuations for sparse random matrices as defined in Definition \ref{asup}. We prove that the extreme eigenvalues have asymptotically Tracy-Widom fluctuation after subtracting the random edge location $\tilde L$ of $\tilde \rho$. As a consequence, the gaps between extreme eigenvalues are asymptotically the same as the gaps between the Airy point process.

\begin{theorem} \label{thm:Tracy-Widom}
Let $H$ be as in Definition \ref{asup} with  $N^{\varepsilon}\leq q\lesssim N^{1/2}$ and eigenvalues given by $\la_1\geq \la_2\geq \cdots\geq \lambda_N$.  There exists an explicit random measure $\tilde \rho$ supported on $[-{\tilde L}, {\tilde L}]$ (depending on certain averaged quantities of $h_{ij}$ as defined in Proposition \ref{p:minfty}).  Fix an integer $k\geq 1$, let $F : \rr^k \to \rr$ be a bounded test function with bounded derivatives.  There is a universal constant $\fc>0$ so that,
\begin{align}\begin{split} \label{eqn:htedgeb1}
&\phantom{{}={}}\ee_{H}[ F (N^{2/3} ( \lambda_1 - {\tilde L}), \cdots , N^{2/3} ( \lambda_k- {\tilde L} ) ]\\
& = \ee_{GOE}[ F (N^{2/3} ( \mu_1 - 2 ), \cdots, N^{2/3} ( \mu_k - 2) ) ]+\OO\left(N^{-\fc}\right).
\end{split}\end{align}
The second expectation  is with respect to a GOE matrix with eigenvalues $\mu_1\geq\mu_2\geq\cdots\geq \mu_N$. 
Analogous results hold for the smallest eigenvalues.
\end{theorem}

\section{Optimal Rigidity Estimates}\label{s:rig}
In this section, we recall the main ingredients for the optimal edge rigidity estimates from  \cite{lee2021higher}, and  improve it to the bulk eigenvalues.

\subsection{Local law for sparse random graphs}

In this section, we recall the following entrywise local semicircle law for sparse random matrices from \cite{MR3098073}. We also collect some estimates following from the local law, which will be used in the rest of this paper.

\begin{theorem}{\normalfont (\cite[Theorem 2.8]{MR3098073})}\label{locallaw}
Let $H$ be as in Definition \ref{asup}.  Let $\fb >0$ be any large constant. % and fix a large constant $\fb>0$. %, and the spectral domain $\cD$ is as defined in \ref{e:domain}. 
Then %for any $D>0$, and $N\geq N(D)$, with probability at least $1-N^{-D}$, 
uniformly for any $z=E+\ri \eta$ such that $-\fb\leq E\leq \fb$ and $0< \eta\leq \fb$,
%\begin{align}
%z\in \{E+\ri\eta: \;-\fb\leq E\leq \fb, \; 1/N\ll \eta\leq \fb\},
%\end{align}
%$z\in \cD$ 
we have
\begin{align}\label{indLocLaw}
\max_{i,j}|G_{ij}(z)-\delta_{ij}m_{ sc}(z)|\prec\left(\frac{1}{q}+\sqrt{\frac{\Im m_{sc}(z)}{N\eta} }+\frac{1}{N\eta}\right),
\end{align}
where $m_{sc}(z)$ is the Stieltjes transform of the semi-circle distribution.
\end{theorem}
As an easy consequence of Theorem \ref{locallaw}, we have that the eigenvectors of $H$ are completely delocalized: with overwhelming probability, uniformly for $1\leq \al\leq N$, 
\begin{align}\label{e:delocal}
\|u_\al\|_\infty\prec \frac{1}{\sqrt N}.
\end{align}
The following estimates utilizing Ward identity and delocalization of eigenvectors \eqref{e:delocal} will be used repeated in the rest of the paper
\begin{align}\begin{split}\label{e:WdI}
    \sum_{j=1}^N|G_{ij}(z)|^2=\frac{\Im[G_{ii}(z)]}{\Im[z]}
    &=\frac{1}{\Im[z]} \Im\left[\sum_{\al=1}^N \frac{u^2_\al(i)}{\la_\al-z}\right]
    =\sum_{\al=1}^N \frac{u^2_\al(i)}{|\la_\al-z|^2}\\
    &\prec \frac{1}{N}\sum_{\al=1}^N \frac{1}{|\la_\al-z|^2}
    =\frac{\Im[m_N(z)]}{\Im[z]}.
\end{split}\end{align}

\subsection{Higher order self-consistent equation}
\label{s:construct}
In this section, we recall the higher order self-consistent equations for sparse random matrices, and some useful estimates on the equilibrium measure and its Stieltjes transform from \cite{lee2021higher,huang2020transition}.
We recall from \cite{lee2021higher}, the random equilibrium measure $\tilde \rho$ in Theorem \ref{t:edgerigidity} and its Stieltjes transform $\tilde m$ 
\begin{align*}
    \tilde m(z)=\int \frac{\tilde \rho(x)\rd x}{x-z},\quad z\in \bC_+.
\end{align*}
Both $\tilde \rho$ and $\tilde m$ are random and depend on $N$ and certain averaged quantities of $h_{ij}$.
They are characterized by a random polynomial $P(z,m)$,
\begin{align}\label{e:deftmz0}
    P(z,\tilde m(z))=0,\quad z\in \bC_+.
\end{align}

Explicitly, $P(z,m)$ is given by
\begin{align}\label{e:defP}
P(z,m)=1+zm+m^2+Q(m),
\end{align}
where 
\begin{align}\label{e:defQ}
Q(m)=a_2m^2+ a_4 m^4 +a_6 m^6+\cdots+a_{2L}m^{2L},
\end{align}
is an even polynomial of $m$, with degree $2L$ (where $L$ is a sufficiently large integer, which will be chosen later). The coefficients $a_{2\ell}$ for $1\leq \ell\leq L$ are explicit random polynomials in the variables $h_{ij}$. To construct them, we need to introduce some notations.

\begin{definition}(Weighted forest).\label{d:forest}
By a weighted forest we mean a finite simple graph which is a union of trees: $\cF=(V(\cF), E(\cF), W(\cF))=(V,E, W)$. Here $V$ is a finite set of vertices, $E$ is a finite set of edges, and 
each edge $e\in E$ connects $\{\al_e, \beta_e\}\in V$. $W$ is a set of edge weights: each edge $e\in E$ is associated with a positive  \emph{odd}  integer $s_e\in \bN$. We denote the number of connected components of $\cF$ as $\theta(\cF)$. 
\end{definition}

We remark that by definition a weighted forest can have arbitrary degrees and weights. But given the total sum of weights $\sum_e(s_e+1)=2\ell$, there are only finite number of weighted forests. We treat the vertices as indeterminate. We will later assign them to be numbers in $\qq{1,N}$, and sum over them. 

For any $1\leq \ell \leq L$, the coefficient $a_{2\ell}$ is a linear combination of terms (with bounded coefficients) in the following form
\begin{align}\label{e:tt}
    w(\cF)= 
    \sum^*_{x_1,x_2,\cdots, x_{|V(\cF)|}}
    \frac{1}{N^{\theta(\cF)}}  \prod_{e\in E(\cF)} w(h_{\al_e \beta_e}; s_e), \quad w(h; s):= h^{s+1}-\frac{\bm1(s=1)}{N},
\end{align}
where  $\cF$ is a weighted forest as in Defintion \ref{d:forest}, $x_1, x_2,\cdots, x_{|V(\cF)|}$ enumerate the vertices of $\cF$ and  weights $s_e$ satifisfy $\sum_e (s_e+1)=2\ell$. In \eqref{e:tt}, we slightly misuse the notations. The summation means we assign each vertex a distinct value in $\qq{1,N}$, $x_1\neq x_2\neq \cdots\neq x_{|V(\cF)|}\in \qq{1,N}$,  and sum over all possible assignments. The summation  can be viewed as a sum over all the possible embeddings of $\cF$ to the complete graph on $N$ vertices $\{1,2,\cdots, N\}$. 

\begin{example}\label{ex:wF}
    If the graph $\cF$ consists of a single edge, with weight $s=1$, then \eqref{e:tt} simplifies to
    \begin{align*}
    w(\cF)= 
    \sum_{1\leq i\neq j\leq N}
    \frac{1}{N}  \left(h_{ij}^2-\frac{1}{N}\right).
\end{align*}
Later one can see from the proof, we will have $a_2=\sum^*_{ij}(h_{ij}^2-1/N)/N$. 
\end{example}

It has been proven in \cite{lee2021higher}, with high probability the random coefficients $a_{2\ell}$ are small and satisfy $|a_{2\ell}|\prec 1/q$. We remark that this bound is not sharp. From the expression of $a_2$ as in Example \ref{ex:wF}, it is not hard to check that  $|a_2|\prec 1/(q\sqrt{N})$. As long as $a_{2\ell}$ all go to zero as $N$ goes to infinite, we can view \eqref{e:deftmz0} as a small perturbation of the equation $1+zm_{sc}(z)+m_{sc}(z)^2=0$, where $m_{sc}(z)$ is the Stieltjes transform of the semi-circle distribution. By a perturbation argument, the solution $\tilde m(z)$ of  $P(z, \tilde m(z))=0$ defines a holomorphic function from the upper half plane $\bC_+$ to itself. It turns out that it is the Stieltjes transform of a probability measure $\tilde \rho$. The following proposition from \cite[Proposition 2.5]{huang2020transition} collects some properties of $\tilde m(z)$ and the measure $\tilde \rho$.

\begin{proposition}\label{p:minfty}
There exists an algebraic function $\tilde m: \bC_+\rightarrow \bC_+$, which depends on the coefficients $a_2, a_4, \cdots, a_{2L}$ of $Q$, such that the following holds:
\begin{enumerate}
\item $\tilde m$ is the solution of the polynomial equation, $P(z, \tilde m(z))=0$.
\item $\tilde m$ is the Stieltjes transform of a symmetric probability measure $\tilde \rho$, with $\supp \tilde \rho=[-{\tilde L}, {\tilde L}]$, where ${\tilde L}$ depends smoothly on the coefficients of $Q$, and its derivatives with respect to the coefficients of $Q$ are uniformly bounded. Moreover,  $\tilde \rho$ is strictly positive on $(-{\tilde L},{\tilde L})$ and has square root behavior at the edge. 
\item We have the following estimate on the imaginary part of $\tilde m$, 
\begin{align} \label{eqn:imtilminf}
\Im[\tilde m(E+\ri \eta)]\asymp \left\{
\begin{array}{cc}
\sqrt{\kappa+\eta}, & \text{if $E\in [-{\tilde L}, {\tilde L}]$,}\\
\eta/\sqrt{\kappa+\eta}, & \text{if $E\not\in [- {\tilde L}, {\tilde L}]$,}
\end{array}
\right.
\end{align}
and 
\begin{align*}
|\del_2 P(z,\tilde m(z))|\asymp \sqrt{\kappa+\eta},\quad 
\del_2^2 P(z,\tilde m(z))=1+\OO(1/q),
\end{align*}
where $\kappa =\dist(\Re[z], \{- {\tilde L},{\tilde L}\})$.
\end{enumerate}
\end{proposition}

\begin{remark}
Since the coefficients of the polynomial $Q$ are random, the edge location $\tilde L$ is also random, depending on certain averaged quantities of $h_{ij}$. We can write $\tilde L$ as its mean plus its fluctuation $\tilde L=L+\Delta L$, where $L=\bE[\tilde L]$. Then $L$ is close to $2$ and has a $1/q$ expansion: $L=2+6\cC_4/q^2+(120\cC_6-81\cC_4^2)/q^4+\cdots$. The leading flucutation of $\tilde L$ is given by $\sum_{ij}(h_{ij}^2-1/N)/N$ which is of size $1/(\sqrt N q)$.
\end{remark}

\subsection{Optimal Rigidity Estimates}
\label{s:rigidity}

In this section we compute higher order moments of the self-consistent equation $P(z, m_N(z))$. This
gives us a recursive moment estimate for the Stieltjes transform $m_N (z)$. The rigidity estimates Theorem \ref{t:rigidity} follow from a careful analysis of the recursive moment estimate and an iteration argument.

We only analyze the behavior of the Stieltjes transform $m_N(z)$ with $z$ close to the right edge $  {\tilde L}$ (or bounded away from $-\tilde L$). The case that $z$ is close to the left edge can be analyzed in the same way. Fix a small constant $\fc>0$, we define the shifted spectral domain $\cD$ as
\begin{align}\label{e:recalldefD}
    \cD=\{\kappa+\ri \eta: |\kappa|\leq 3/2, 0<\eta\leq 1, N\eta \sqrt{|\kappa|+\eta}\geq N^\fc\}.
\end{align}

\begin{proposition}\label{p:DSE}
Let $H$ be as in Definition \ref{asup} with $N^{\varepsilon}\leq q\lesssim N^{1/2}$. There exist a finite number of random correction terms $a_2, a_4, \cdots, a_{2L}$ and polynomial $P$ as in \eqref{e:defQ}. For $z={\tilde L}+w$ where ${\tilde L}$ is from Proposition \ref{p:minfty} and $w=\kappa+\ri \eta\in \cD$, let $m_N(z)$ be the Stieltjes transform of the empirical eigenvalue density of $H$, then we have
\begin{align}\label{e:Pmoment}
\bE[|P( z, m_{N}( z)|^{2r}]\prec \bE[\Phi_r(w)],
\end{align}
where $\Phi_r(w)$ is defined as
\begin{align}\label{e:defPhir}
\begin{split}
\Phi_r(w):&= \sum_{s=1}^{2r}\left(\frac{\Im[m_N(z)]}{N\eta}\right)^s|P( z, m_{N}( z)|^{2r-s}\\
&+ \sum_{s=1}^{2r}\frac{1}{N\eta}\left( \frac{\Im[m_N(z)]|\del_2 P(z,m_N(z))|}{N\eta}+\frac{1}{N}\right)^s |P( z, m_N( z))|^{2r-s-1}.
\end{split}\end{align}
\end{proposition}

Proposition \ref{p:DSE} improves \cite[Proposition 3.1]{lee2021higher}. In the control parameter $\Phi_r$ \eqref{e:defPhir}, comparing with  \cite[Proposition 3.1]{lee2021higher} we no longer have the term
\begin{align}\label{e:badterm}
     \sum_{s=1}^{2r}\sqrt{\frac{\Im[m_N(z)]}{N\eta}}\left( \frac{\Im[m_N(z)]|\del_2 P(z,m_N(z))|}{N\eta}+\frac{1}{N}\right)^s |P( z, m_N( z))|^{2r-s-1}.
\end{align}
Inside the bulk of the spectral, namely $\kappa=\Omega(1)$, we have $\Im[m_N(z)]=\Omega(1)$. The error \eqref{e:badterm} is $\sqrt{N\eta}$ bigger than the second term on the righthand side of \eqref{e:defPhir}. Getting rid of the error \eqref{e:badterm} is crucial to the optimal rigidity estimates for eigenvalues inside the bulk of the spectral.

In this proof we write, for simplicity of notation, $z=\tilde L+w$ with $w=\kappa+\ri \eta\in \cD$, $P=P(  z, m_N(  z))$, $G=G(  z)$, $m_N=m_N(  z)$, $P'=\del_{2}P(  z, m_N(  z))$, $\del_{ij}=\del_{h_{ij}}$, $\del_{ij}G=\del_{ij}(G(  z))$, $\del_{ij}^pm_N=\del_{ij}(m_N(  z))$. We also denote $D_{ij}G=(\del_{ij}G)(  z)$ and $D_{ij}m_N=(\del_{ij}  m_N)(  z)$, where the derivatives do not hit $z$.

A central object in our proof is the following notion of a polynomial in the entries of the Green's function. 
\begin{definition} \label{d:evaluation}
  Let $R = R(\{x_{st}\}_{s,t = 1}^r, y)$ be a monomial in the $r^2+1$ abstract variables $\{x_{st}\}_{s,t = 1}^r, y$. We denote its degree by $\deg(R)$. For ${\bf i} \in \qq{N}^r$, we define its \emph{evaluation} on the Green's function and the Stieltjes transform by
  \begin{equation}
    R_{\bf i} = R(\{G_{i_s i_t}\}_{s,t = 1}^r, m),
  \end{equation}
  and say that $R_{\bf i}$ is a \emph{monomial in the Green's function entries} $\{G_{i_s i_t}\}_{s,t = 1}^r$. 
  We denote the number of off-diagonal entries of $R_\bfi$ as  $\chi(R_\bfi)$.
    %  \item
%  Let $F = F(\{x_{st}\}_{s,t = 1}^r)$ be a monomial in $r^2$ variables.
%  Then the number of off-diagonal entries of $F$ is the total degree of variables $x_{st}$ with $s \neq t$. If the number of off-diagonal entries of $F$ is zero then we define $\chi_F = 1$, otherwise we define $\chi_F = 0$.
\end{definition}

In \cite{lee2021higher}, the terms with one off-diagonal Green's function are bounded using the Ward identity \eqref{e:WdI}, which leads to an error in the form \eqref{e:badterm}. Instead of using the Ward identity, we estimate such terms using the $L_2$ norm of the Green's function, i.e. $\|G\|_2\leq 1/\eta$, which leads to the following estimate:
\begin{proposition}\label{p:one-off}
Adopt the assumptions of Proposition \ref{p:DSE}. Given a weighted forest $\cF$ with vertex set $V(\cF)=ij\bm m$, where $\bm m=\{m_1, m_2,\cdots, m_r\}$. Then for any monomial $R_{ij\bm m}$ as in Definition \ref{d:evaluation} and nonnegative integers  $\bm p=\{p_{e}\}_{e\in E(\cF)}$ such that $p_e\geq s_e$, we have
\begin{align}\label{e:one-off}
    \frac{1}{N^{r+2}}\sum_{ij\bm m}^*\bE\left[ G_{ij}R_{ij\bm m}\left(\prod_{e\in E(\cF)}\del_{\al_e\beta_e}^{p_{e}-s_e}\right)(P^{r-1}\bar P^{r}) \right]\prec \bE[\Phi_r],
\end{align}
and for any choice of $\bm m\in \qq{1,N}$ (possibly not distinct), it holds
\begin{align}\label{e:derbb2}
    \left(\frac{1}{N}+\frac{\Im[m_N]}{N\eta}\right)\left|\left(\prod_{e\in E(\cF)}\del_{\al_e\beta_e}^{p_{e}-s_e}\right)(P^{r-1}\bar P^{r}) \right|\prec \Phi_r,
\end{align}
where $\Phi_r$ is as defined in \eqref{e:defPhir}.
\end{proposition}

Using Proposition \ref{p:one-off} as input, in the following, we give a proof of Proposition \ref{p:DSE} following \cite[Proposition 3.1]{lee2021higher}. 
  We postpone the proof of Proposition \ref{p:one-off} and some estimates to next section. 
The proof uses the cumulant expansion to compute the expectations. It turns out, all the terms we will get  in the expansion are in the following form
\begin{align}\label{e:ordero}
    \frac{1}{q^\fo}\times \frac{1}{N^r}\sum_{\bm m}^*\bE\left[R_{\bm m} \left(\prod_{e\in E(\cF)}\del_{\al_e\beta_e}^{p_{e}-s_e}\right)(P^{r-1}\bar P^r)\right],
\end{align}
where $\cF$ is a weighted forest  with vertex set $V(\cF)=\bm m=\{m_1, m_2,\cdots, m_r\}$, $R$ is a monomial as in Definition \ref{d:evaluation},   $\bm p=\{p_{e}\}_{e\in E(\cF)}$ are nonnegative integers, and  $\fo\geq 0$ is the order parameter. Since the second factor in \eqref{e:ordero} can be trivially bounded by $\OO_\prec(1)$, the whole expression can be bounded by $\OO_\prec(1/q^\fo)$.  For terms with order at least $M$, we will trivially bound them by $\OO_\prec(1/q^M)$.

\begin{proof}[Proof of Proposition \ref{p:DSE}]
We divide the proof of Proposition \ref{p:DSE} into three steps.

\noindent\textit{Step $1$ (eliminate off-diagonal Green's function terms) }
The starting point is the following identity,
\begin{align*}
1+  z m_N(  z)=\sum_{ij}h_{ij}G_{ij}(  z).
\end{align*}
Using the cumulant expansion, we can write the moment of $P(  z, m_N(  z))$ as,
\begin{align}
\begin{split}\label{e:Pmoment}
&\phantom{{}={}}\bE[|P(  z,m_N(  z))|^{2r}]
=\bE\left[(1+zm_N(z))P^{r-1}\bar P^{r}\right]+\bE[(m_N^2+Q) P^{r-1}\bar P^{r}]\\
&=\frac{1}{N}\bE\left[\sum_{ij}h_{ij}G_{ij}P^{r-1}\bar P^{r}\right]+\bE[(m_N^2+Q) P^{r-1}\bar P^{r}]\\
&=\frac{1}{N}\sum_{p=1}^{ M}\sum_{ij}\frac{\cC_{p+1}}{Nq^{p-1}}\bE[\del_{ij}^{p}(G_{ij}P^{r-1}\bar P^r)]+\OO_\prec\left(\frac{1}{q^{ M}}\right)+\bE[(m_N^2+Q) P^{r-1}\bar P^{r}]\\
&=\frac{1}{N}\sum_{ij}\sum_{p=1}^{ M}\sum_{s= 0}^p\frac{\cC_{p+1}}{Nq^{p-1}}{p\choose s}\bE[\del_{ij}^{s}G_{ij}\del_{ij}^{p-s}(P^{r-1}\bar P^r)]+\OO_\prec\left(\bE[\Phi_r]\right)+\bE[(m_N^2+Q) P^{r-1}\bar P^{r}],
\end{split}
\end{align}
where $ M$ is large enough such that $1/q^ M\prec \Phi_r$.
For the last line in \eqref{e:Pmoment}, we will prove later in Proposition \ref{p:dermN}, $\del_{ij}^s G_{ij}=D_{ij}^s G_{ij}+\OO_\prec(\Im[m_N]/N\eta)\prec 1$. We can replace $\del_{ij}$ by $D_{ij}$ and error can be bounded by $\bE[\Phi_r]$. Thus combining with \eqref{e:derbb2}, we can bound the terms with $i=j$ in \eqref{e:Pmoment} by $\bE[\Phi_r]$
\begin{align}\label{e:i=jterm}
\frac{1}{N}\sum_{i=j}\sum_{p=1}^{ M}\sum_{s= 0}^p\frac{\cC_{p+1}}{Nq^{p-1}}{p\choose s}\bE[\del_{ij}^{s}G_{ij}\del_{ij}^{p-s}(P^{r-1}\bar P^r)]
\prec \sum_{p=1}^{ M}\sum_{s= 0}^p\sum_{i}\frac{1}{N^2}\bE[|\del_{ii}^{p-s}(P^{r-1}\bar P^r)|]\prec \bE[\Phi_r].
\end{align}
With \eqref{e:i=jterm}, we can rewrite the summation on the righthand side of \eqref{e:Pmoment} as 
 \begin{align}\label{e:Pmomentstar}
\frac{1}{N}\sum^*_{ij}\sum_{p=1}^{ M}\sum_{s= 0}^p\frac{\cC_{p+1}}{Nq^{p-1}}{p\choose s}\bE[\del_{ij}^{s}G_{ij}\del_{ij}^{p-s}(P^{r-1}\bar P^r)]+\bE[(m_N^2+Q) P^{r-1}\bar P^{r}]+\OO_\prec\left(\bE[\Phi_r]\right).
\end{align}
The derivatives $D_{ij}^s G_{ij}$ is a sum of terms in the form $G^a_{ii}G^b_{jj}G^c_{ij}$. Thanks to Proposition \ref{p:one-off}, terms with $c\geq 1$ is bounded by $\OO_\prec(\bE[\Phi_r])$, and 
\begin{align*}
D_{ij}^s G_{ij}=-\bm1(s \text{ is odd})s! G_{ii}^{\frac{s+1}{2}}G_{jj}^{\frac{s+1}{2}}+\{\text{terms with off-diagonal entries}\}.
\end{align*}
Therefore, the leading terms in \eqref{e:Pmoment} are those which do not contain any off-diagonal Green's function terms,
\begin{align}\label{e:nooff}\begin{split}
    &\phantom{{}={}}\frac{1}{N}\sum^*_{ij}\sum_{p=1}^{ M}\sum_{s= 0}^p\frac{\cC_{p+1}}{Nq^{p-1}}{p\choose s}\bE[\del_{ij}^{s}G_{ij}\del_{ij}^{p-s}(P^{r-1}\bar P^r)]+\OO_\prec\left(\bE[\Phi_r]\right)+\bE[(m_N^2+Q) P^{r-1}\bar P^{r}]\\
    &=-\frac{1}{N}\sum_{p=1}^{ M}\sum_{s\text{ odd}}\sum^*_{ij} \frac{\cC_{p+1}s!}{Nq^{p-1}}{p\choose s}\bE\left[G_{ii}^{\frac{s+1}{2}}G_{jj}^{\frac{s+1}{2}}\del_{ij}^{p-s}(P^{r-1}\bar P^r)\right]+\bE[(m_N^2+Q) P^{r-1}\bar P^r]+\OO_\prec(\bE[\Phi_r])\\
    &=-\frac{1}{N}\sum_{p=2}^{ M}\sum_{s\text{ odd}}\sum^*_{ij} \frac{\cC_{p+1}s!}{Nq^{p-1}}{p\choose s}\bE\left[G_{ii}^{\frac{s+1}{2}}G_{jj}^{\frac{s+1}{2}}\del_{ij}^{p-s}(P^{r-1}\bar P^r)\right]+\bE[Q P^{r-1}\bar P^r]+\OO_\prec(\bE[\Phi_r]).
\end{split}
\end{align}
Here in the last line, we used that the term corresponding to $p=1,s=1$
\begin{align*}
-\frac{1}{N}\sum^*_{ij} \frac{\cC_{2}}{N}\bE\left[G_{ii}G_{jj}(P^{r-1}\bar P^r)\right]
&=-\bE[m^2P^{r-1}\bar P^r]+\OO\left(\frac{1}{N}\bE[|P|^{2r-1}] \right)\\
&=-\bE[m^2P^{r-1}\bar P^r]+\OO(\bE[\Phi_r]),
\end{align*}
which cancels with $\bE[m^2P^{r-1}\bar P^r]$ in the second term.

We remark that for any given $p$, the term on the righthand side of \eqref{e:nooff} is in the form (up to some bounded multiplicative factor)
\begin{align}\label{e:final-1}
  \frac{1}{N^{|V(\cF)|}}\sum_{ 2\leq p\leq M}
  \sum^*_{\bm m}\bE\left[R_{\bm m}
  \prod_{e\in E(\cF)}\frac{\cC_{p_{e}+1}s_e!}{q^{p_e-1}}{p_e\choose s_e}  \del_{\al_e \beta_e}^{p_e-s_e}(P^{r-1}\bar P^r)\right],
\end{align}
where  the weighted forest $\cF$ as in Definition \ref{d:forest} is a single edge $e=\{m_1, m_2\}$ with vertex set $\bm m=\{m_1, m_2\}$. The edge $e$ has weight $s_e$. The monomial $R_{\bm m}$ as in \eqref{d:evaluation} has no off-diagonal entries, i.e. $\chi(R_{\bm m})=0$, and total degree $\deg (R_{\bm m})=s_e+1$.

\noindent\textit{Step $2$ (replace diagonal Green's function entries by $m_N$) }

For the diagonal terms $G_{ii}^{\frac{s+1}{2}}G_{jj}^{\frac{s+1}{2}}$ in \eqref{e:nooff} we can replace them by diagonal terms with different indices using the following proposition. 

\begin{proposition}\label{p:diag}
Adopt the assumptions of Proposition \ref{p:DSE}. Given a weighted forest $\cF$ with vertex set $V(\cF)=i\bm m$, where $\bm m=\{m_1, m_2,\cdots, m_r\}$. Then for any monomial $R_{\bm m}$ as in Definition \ref{d:evaluation} with no off-diagonal Green's function entries, i.e. $\chi(R_{\bm m})=0$, and  integers  $\bm p=\{p_{e}\}_{e\in E(\cF)}$ with $p_e-s_e\geq 0$, we have
\begin{align}\begin{split}\label{e:diag}
    &\frac{1}{N^{r+1}}\sum^*_{i\bm m}\bE\left[ G^\al_{ii}R_{\bm m}V\right]
    =\frac{1}{N^{r+1}}\sum^*_{i\bm m}\bE\left[ G^{\al-1}_{ii} m_N R_{\bm m}V\right]
     + \Omega_1+\Omega_2+\OO_{\prec}(\bE[\Phi_r]),\\
    &V=\left(\prod_{e\in E(\cF)}\del_{\al_e \beta_e}^{p_e-s_e}\right)(P^{r-1}\bar P^{r}), 
\end{split}\end{align}
where the $\Omega_1, \Omega_2$ are given by
\begin{align}\begin{split}\label{e:hod}
&\Omega_1=-\frac{1}{N^{r+2}}\sum_{p=2}^M\sum_{s \text{ odd}} \frac{\cC_{p+1}s!}{Nq^{p-1}}{p\choose s}
   \sum^*_{ij k \bm m} \bE[G^\al_{ii}G_{j j}^{\frac{s+1}{2}}G_{kk}^{\frac{s+1}{2}} R_{\bm m}\del^{p-s}_{j k} V)],\\
& \Omega_2=  \frac{1}{N^{r+1}}\sum_{p=2}^M\sum_{s \text{odd}} \frac{\cC_{p+1}}{Nq^{p-1}}{p\choose s}s! {\frac{s-1}{2}+\al-1\choose \al-1}
    \sum^*_{i k \bm m}\bE[m_NG^{\al+\frac{s-1}{2}}_{ii}G^{\frac{s+1}{2}}_{kk}R_{\bm m}\del^{p-s}_{ik} V)].\end{split}\end{align}
%\begin{align}\begin{split}\label{e:hod}
%-\frac{1}{N^{r+2}}\sum_{ij k \bm m}\sum_{p=2}^M\sum_{s \text{odd}} \frac{\cC_{p+1}s!}{Nq^{p-1}}{p\choose s}
%    \left(\bE\left[G^\al_{ii}G_{j j}^{\frac{s+1}{2}}G_{kk}^{\frac{s+1}{2}} \left(\prod_{x\in i\bm m}G^{\al_x}_{xx}\right)\left(\del^{p-s}_{j k} \prod_{x\neq y\in i\bm m}\del_{xy}^{\beta_{xy}}\right)(P^{r-1}\bar P^{r})\right]\right.\\
%-   \left.{\frac{s-1}{2}+\al-1\choose \al-1}
%    \bE\left[G^{\al+\frac{s-1}{2}}_{ii}G_{jj}G^{\frac{s+1}{2}}_{kk}\left(\prod_{x\in i\bm m}G^{\al_x}_{xx}\right) \left(\del^{p-s}_{ik}\prod_{x\neq y\in i\bm m}\del_{xy}^{\beta_{xy}}\right)(P^{r-1}\bar P^{r})\right]\right).
%\end{split}\end{align}
Comparing these terms $\Omega_1, \Omega_2$ with \eqref{e:diag}, since $p\geq 2$ they are of order at least $1$ (recall from \eqref{e:ordero}).
\end{proposition}

\begin{remark}
These terms $\Omega_1, \Omega_2$ in \eqref{e:hod} are in the same form as in \eqref{e:diag}. With given $p$, the term in $\Omega_1$  is associated with a weighted forest $\cF_1$, which is from $\cF$ by adding vertices $j,k$ and an edge $\{j,k\}$ with weight $s$. In total $\Omega_1$ has $r+2$ vertices.  $\Omega_2$ is associated with a weighted forest $\cF_2$, which is from $\cF$ by adding one vertex $k$ and an edge $\{i,k\}$ with weight $s$. In total $\Omega_2$ has $r+1$ vertices.
For given $s$, both $\Omega_1, \Omega_2$ have an extra derivative $\del^{p-s}$ and the total number of diagonal Green's function entries increases by $s+1$. 
\end{remark}

Thanks to Proposition \ref{p:diag}, we can replace one copy of $G_{ii}$ in the last line of \eqref{e:nooff} by $m_N$: 
\begin{align}\begin{split}\label{e:cp1}
&\phantom{{}={}}-\frac{1}{N}\sum_{p=2}^{ M}\sum_{s\text{ odd}}\sum^*_{ij} \frac{\cC_{p+1}s!}{Nq^{p-1}}{p\choose s}\bE\left[G_{ii}^{\frac{s+1}{2}}G_{jj}^{\frac{s+1}{2}}\del_{ij}^{p-s}(P^{r-1}\bar P^r)\right]\\
&=-\frac{1}{N}\sum_{p=2}^{ M}\sum_{s\text{ odd}}\sum^*_{ij} \frac{\cC_{p+1}s!}{Nq^{p-1}}{p\choose s}\bE\left[G_{ii}^{\frac{s-1}{2}}m_NG_{jj}^{\frac{s+1}{2}}\del_{ij}^{p-s}(P^{r-1}\bar P^r)\right]+\OO_\prec(\bE[\Phi_r])\\
&+\frac{1}{N^2}\sum_{{s, s'} \text{ odd}}\sum_{{p,p'}=2}^M\frac{\cC_{p+1}s!}{Nq^{p-1}}{p\choose s}\frac{\cC_{{p'}+1}s'!}{Nq^{{p'}-1}}{p'\choose s'}
   \sum^*_{ijk\ell} \bE[(G^{\frac{s+1}{2}}_{ii}G^{\frac{s+1}{2}}_{jj}G^{\frac{{s'}+1}{2}}_{kk}G_{\ell\ell}^{\frac{{s'}+1}{2}}\del^{{p'}-{s'}}_{\ell k}\del_{ij}^{p-s}(P^{r-1}\bar P^r))]\\
&    - \frac{1}{N}\sum_{{s, s'} \text{ odd}}\sum_{{p,p'}=2}^M\frac{\cC_{p+1}s!}{Nq^{p-1}}{p\choose s}\frac{\cC_{{p'}+1}s'!}{Nq^{{p'}-1}}{p'\choose s'}
   \sum^*_{ijk}{\frac{s+s'}{2}-1\choose \frac{s-1}{2}}
    \bE[(m_NG^{\frac{s+s'}{2}}_{ii}G^{\frac{s+1}{2}}_{jj}G_{kk}^{\frac{s'+1}{2}}\del^{{p'}-{s'}}_{i k}\del_{ij}^{p-s}(P^{r-1}\bar P^r))].
\end{split}\end{align}
We remark that since $p'\geq 2$, the last two terms on the righthand side of \eqref{e:cp1} are of higher order than the first term on the righthand side of \eqref{e:cp1}. By repeating this procedure, we can replace $G_{ii}^{(s+1)/2}G_{jj}^{(s+1)/2}$ in \eqref{e:nooff} by $m_N^{s+1}$, to get the leading terms
\begin{align}\label{e:nextt}\begin{split}
%    &\phantom{{}={}}-\frac{1}{N}\sum_{p=2}^{ M}\sum_{s\text{ odd}}\sum_{ij} \frac{\cC_{p+1}s!}{Nq^{p-1}}{p\choose s}\frac{1}{N^{s+1}}\sum_{1\leq k_1,k_2,\cdots,k_{s+1}\leq N}\bE\left[G_{k_1k_1}G_{k_2k_2}\cdots G_{k_{s+1}k_{s+1}}\del_{ij}^{p-s}(P^{r-1}\bar P^r)\right]\\
     -\frac{1}{N}\sum_{p=2}^{ M}\sum_{s\text{ odd}}\sum^*_{ij} \frac{\cC_{p+1}s!}{Nq^{p-1}}{p\choose s}\bE\left[m_N^{s+1}\del_{ij}^{p-s}(P^{r-1}\bar P^r)\right],
\end{split}\end{align}
with higher order terms (which have at least one more copies of $1/q$) as linear combination of terms (with bounded coefficients) in the form 
\begin{align}\label{e:final0}
  \frac{1}{N^{r}}\sum_{ 2\leq p_{e_1},\cdots, p_{e_{|E(\cF)|}}\leq M}
  \sum^*_{\bm m}\bE\left[R_{\bm m}
  \prod_{e\in E(\cF)}\frac{\cC_{p_{e+1}}s_e!}{q^{p_e-1}}{p_e\choose s_e}  \del_{\al_e \beta_e}^{p_e-s_e}(P^{r-1}\bar P^r)\right],
\end{align}
where  $\cF$ is a weighted forest as in Definition \ref{d:forest} with vertex set $\bm m=\{m_1,m_2,\cdots, m_r\}$; the monomial $R_{\bm m}$ has no off-diagonal entries, i.e. $\chi(R_{\bm m})=0$ and $\deg(R_{\bm m})=\sum_{e\in E(\cF)}(s_e+1)$.  

We can repeat Step 2 for these higher order terms \eqref{e:final0}. If a term has order bigger than $M$, we can trivially bound it by $1/q^M\prec \Phi_r$ as in \eqref{e:ordero}.  
The final expression is a linear combinations (with bounded coefficients) of terms in the form: 
\begin{align}\begin{split}\label{e:final}
   &\bE[m_N^{2\ell} L_\cF(P^{r-1}\bar P^{r})],\\
   &L_\cF= \sum_{2\leq p_{e_1},\cdots, p_{e_{|E(\cF)|}}\leq M}\frac{1}{N^{|V(\cF)|}}\sum_{ x_1,\cdots, x_{|V(\cF)|}}^*\prod_{e\in E(\cF)}\frac{\cC_{p_{e+1}} s_e!}{Nq^{p_e-1}}{p_e\choose s_e}  \del_{\al_e\beta_e}^{p_e-s_e},
\end{split}\end{align}
where  $\cF$ is a forest as in Defintion \ref{d:forest}, $x_1, x_2,\cdots, x_{|V(\cF)|}$ enumerate the vertices of $\cF$. Moreover, all the weights $s_e$ are odd positive integers, and the total weights satisfies $\sum_e (s_e+1)=2\ell$. The above discussion leads to the following claim.

\begin{claim}\label{c:step2}
Under the assumptions of Proposition \ref{p:DSE}, with an error $\OO_\prec(\bE[\Phi_r])$, the first term on the righthand side of \eqref{e:nooff} is a linear combinations of terms in the form,
\begin{align}\label{e:finalk} 
    \bE[m_N^{2\ell} L_\cF(P^{r-1}\bar P^{r})],
\end{align}
where $\cF$ is a forest as in Defintion \ref{d:forest}, and $2\ell=\sum_{e\in E(\cF)}(s_e+1)$ and $L_\cF$ is as defined in \eqref{e:final}.
\end{claim}
\noindent\textit{Step $3$ (rewrite differential operators as an expectation) }
Finally by the cumulant expansion we have, 
\begin{align*}
    &\phantom{{}={}}\bE[h_{ij}^{s+1}m_N^{2\ell}(P^{r-1}\bar P^{r})]
    =\sum_{p=1}^M\frac{\cC_{p+1}}{Nq^{p-1}} \bE[\del_{ij}^{p}(h_{ij}^{s}m_N^{2\ell}P^{r-1}\bar P^{r})]+\OO_\prec(\bE[\Phi_r])\\
    &=\sum_{p=1}^M\frac{\cC_{p+1}}{Nq^{p-1}}{p\choose s',s''} \bE[\del_{ij}^{s'}(h_{ij}^{s})\del_{ij}^{s''}(m_N^{2\ell})\del_{ij}^{p-s'-s''}(P^{r-1}\bar P^{r})]+\OO_\prec(\bE[\Phi_r]).
\end{align*}
If $s''\geq 1$, then we have $|\del_{ij}^{s''}(m_N^{2\ell})|\prec \Im[m_N/N\eta]$ from \eqref{e:dmN} in Proposition \ref{p:dermN}. Then it follows
\begin{align*}
\del_{ij}^{s'}(h_{ij}^{s})\del_{ij}^{s''}(m_N^{2\ell})\del_{ij}^{p-s'-s''}(P^{r-1}\bar P^{r})\prec \frac{\Im[m_N]}{N\eta}\del_{ij}^{p-s'-s''}(P^{r-1}\bar P^{r})\prec \Phi_r,
\end{align*}
where we used \eqref{e:derbb2} in the last inequality. If $s''=0$ and $s'<s$, then $\del_{ij}^{s'}(h_{ij}^{s})=s(s-1)\cdots (s-s'+1) h_{ij}^{s-s'}$, we can do another cumulant expansion
\begin{align*}
&\phantom{{}={}}\bE[\del_{ij}^{s'}(h_{ij}^{s})\del_{ij}^{p-s'}(P^{r-1}\bar P^{r})]\\
&=\sum_{p'=1}^M\frac{\cC_{p'+1}s\cdots (s-s'+1) }{Nq^{p'-1}} \bE[
\del_{ij}^{p'}(h_{ij}^{s-s'-1}\del_{ij}^{p-s'}(P^{r-1}\bar P^{r}) )]+\OO_\prec\left(\frac{1}{q^M}\right)
\prec \bE[\Phi_r],
\end{align*}
where we used \eqref{e:derbb2} in the last inequality.  The remaining terms correspond to $s'=s$ and $s''=0$. we conclude 
\begin{align*}
    \bE\left[h_{ij}^{s+1}m_N^{2\ell}(P^{r-1}\bar P^{r})\right]
    =\sum_{p=1}^M\frac{\cC_{p+1}s!}{Nq^{p-1}}{p\choose s} \bE[m_N^{2\ell}\del_{ij}^{p-s}(P^{r-1}\bar P^{r})]+\OO_\prec(\bE[\Phi_r]).
\end{align*}
And by moving the term corresponding to $p=1$ to the left,  we can rewrite the above equation as
\begin{align}\label{e:ediff}
    \bE\left[\left(h_{ij}^{s+1}-\frac{\bm1(s=1)}{N}\right)m_N^{2\ell}(P^{r-1}\bar P^{r})\right]
    =\sum_{p=2}^M\frac{\cC_{p+1}s!}{Nq^{p-1}}{p\choose s} \bE[m_N^{2\ell}\del_{ij}^{p-s}(P^{r-1}\bar P^{r})]+\OO_\prec(\bE[\Phi_r]).
\end{align}
By repeatedly using the relation \eqref{e:ediff}, we can rewrite terms as in \eqref{e:final}, and have proved the following claim.
\begin{claim}\label{c:step3}
Under the assumptions of Proposition \ref{p:DSE}, we can rewrite \eqref{e:finalk} as:
\begin{align}\label{e:final2}\begin{split}
    &\phantom{{}={}}\bE[m_N^{2\ell}L_\cF(P^{r-1}\bar P^{r})]\\
    &=\prod_{e\in E(\cF)} \bE\left[\sum_{  x_1, \cdots, x_{|V(\cF)|}}^*\left(\frac{1}{N^{\theta(\cF)}}  \prod_{e\in E(\cF)} \left(h_{\al_e\beta_e}^{s_e+1}-\frac{\bm1_{s_e=1}}{N}\right)\right) m_N^{2\ell}(P^{r-1}\bar P^{r})\right]+\OO_\prec(\bE[\Phi_r])\\
    &=\prod_{e\in E(\cF)} \bE\left[w(\cF) m_N^{2\ell}(P^{r-1}\bar P^{r})\right]+\OO_\prec(\bE[\Phi_r]),
\end{split}\end{align}
where  $\cF$ is a forest as in Defintion \ref{d:forest}, and $x_1, x_2,\cdots, x_{|V(\cF)|}$ enumerate the vertices of $\cF$. $L_\cF$ and $w(\cF)$ are as defined in \eqref{e:final} and \eqref{e:tt}. Moreover, all the weights $s_e$ are odd positive integers with $\sum_e (s_e+1)=2\ell$.  

\end{claim}
Thanks to Claims \eqref{c:step2} and \eqref{c:step3}, up to an error $\OO_\prec(\bE[\Phi_r])$,
the first term on the righthand side of \eqref{e:nooff} is in the form 
\begin{align}\label{e:ama}
    -\bE[(a_2m_N^2+a_4 m_N^4+\cdots+a_{2L}m_N^{2L})(P^{r-1}\bar P^r)]+\OO_\prec(\bE[\Phi_r]),
\end{align}
where
$a_{2\ell}$ is a sum of terms in the form $w(\cF)$ as in \eqref{e:tt}, where  $\cF$ is a forest as in Defintion \ref{d:forest}. Moreover, all the weights $s_e$ are odd positive integers with $\sum_e (s_e+1)=2\ell$.   
We can use the expression \eqref{e:ama} as the definition of the polynomial $Q$. Thus 
the term \eqref{e:ama} cancels with $\bE[Q P^{r-1}\bar P^r]$ in \eqref{e:Pmoment}, and we conclude Proposition \ref{p:DSE}.

\end{proof}
\subsection{Proof of Proposition \ref{p:one-off} and \ref{p:diag}}

Before proving Propositions \ref{p:one-off} and \ref{p:diag}, we collect some useful estimates of the derivatives of $m_N(z)$ and $a_{2\ell}$ in Propositions \ref{p:a2kstr} and \ref{p:dermN}.

\begin{proposition}\label{p:a2kstr}
Adopt the assumptions in Proposition \ref{p:DSE}, and take $z=\tilde L+w$. Fix distinct indices $i,j,\bm m=\{m_1, m_2,\cdots, m_r\}$. We consider the differential operator $\del^{\bm \beta}$ with $\bm\beta=\{\beta_{uv}\}_{u,v\in ij\bm m}$
\begin{align*}
    \del^{\bm \beta}=\prod_{u,v\in ij\bm m}\del_{uv}^{\beta_{uv}},\quad D^{\bm \beta}=\prod_{u,v \in ij\bm m} D_{uv}^{\beta_{uv}}, \quad |\bm\beta|=\sum_{u,v\in \bm m} \beta_{uv}\geq 1.
\end{align*}
\begin{enumerate}

\item The derivative $D^{\bm\beta} m_N(z)$ of $m_N(z)$ is a linear combination of terms in the following form (with bounded coefficients)
\begin{align}\label{e:dermN}
      \frac{\Im[m_N(z)]}{N\eta}\sum_{\bmk}G^a_{ij} X_{\bmk\bm m i}Y_{\bmk\bm m j},
\end{align}
where $a\geq 0$, $\bmk$ is an index set,  and $\sum_{\bmk}|X_{ \bmk \bm m i}|^2,\sum_{\bmk} |Y_{\bmk \bm m j}|^2=\OO_\prec(1)$. %More precisely, 
%\begin{align}
%D^{\bm\beta} m_N(z)=\sum_{\gamma} C_\gamma \frac{\Im[m_N]}{N\eta}\sum_{\bmk_\gamma}G_{ij}^{a_\gamma} X^{(\gamma)}_{\bmk_\gamma\bm m i}Y^{(\gamma)}_{\bmk_\gamma\bm m j}.
%\end{align}

\item  The derivative $\del^{\bm\beta} a_{2\ell}$ of $a_{2\ell}$ (as defined in  \eqref{e:tt} ) 
is a linear combination of terms in the following form (with bounded coefficients)
\begin{align}\label{e:derhij}
   \frac{1}{N}\sum_{\bmk:ij\bm m \bmk \atop \text{distinct}}h_{ij}^a  X_{ \bmk \bm m i}  Y_{\bmk\bm m j},
\end{align}
where $a\geq 0$, $\bmk$ is an index set,  $\sum_{\bmk}|X_{\bmk \bm m i}|^2,\sum_{\bmk} |Y_{\bmk\bm m j}|^2=\OO_\prec(1)$. 
\end{enumerate}
\end{proposition}
\begin{proof}
The derivative $D^{\bm\beta}m_N$ is a linear combination of terms in the following form
\begin{align}\label{e:dermN2}
    \frac{1}{N}\sum_{k=1}^nG_{kv_1}G_{v_2 v_3}\cdots G_{v_{2\ell-2} v_{2\ell-1}} G_{v_{2\ell}k},
\end{align}
where $v_1,v_2,\cdots, v_{2\ell}\in ij\bm m$. For the Green's function entries in \eqref{e:dermN2} we can regroup them depending if they contain indices $i,j$
\begin{align*}
   \frac{1}{N} \sum_{k=1}^N G_{ij}^a \tilde X_{k\bm m i} \tilde Y_{k\bm m j}, 
\end{align*}
where for each $G_{xy}$ in \eqref{e:dermN2}, if the index set $\{x,y\}$ only contains $i$ we put it in $\tilde X_{k\bm m i}$; if it only contains $j$ we put it in $\tilde Y_{k\bm m j}$; if it does not contain $i,j$ we simply put it in $\tilde X_{k\bm m j}$. 
%$G_{k v_1}$ if $v_1=i$ we put it in $\tilde X_{k\bm m i}$, otherwise we put it in $\tilde Y_{k\bm m i}$; for $G_{v_{2\ell}k }$ if $v_{2\ell}=i$ we put it in $\tilde X_{k\bm m i}$, otherwise we put it in $\tilde Y_{k\bm m i}$; for $G_{v_{s}v_{s+1}}$ if $\{v_s, v_{s+1}\}=\{i,j\}$ we seperate it out, otherwise if $i\in \{v_s, v_{s+1}\}$ we put it in $\tilde X_{k\bm m i}$, if $i\notin \{v_s, v_{s+1}\}$ we put it in $\tilde X_{k\bm m j}$. 

There are two cases. In the first case, each of $\tilde X_{k\bm m i}, \tilde Y_{k\bm m j}$ contains  one of $G_{k v_1}, G_{v_{2\ell}k}$; in the second case, both $G_{k v_1}, G_{v_{2\ell}k}$ are in $\tilde X_{k\bm m i}$ or $\tilde Y_{k\bm m j}$.

In the first case, say $G_{k v_1}$ is in $\tilde X_{k\bm m i}$, using  $|G_{v_s v_{s+1}}|\prec 1$ from Theorem \ref{locallaw} and Ward identity \eqref{e:WdI} we have
\begin{align*}
    \sum_{k=1}^N |\tilde X_{k\bm m i}|^2 \prec \sum_{k=1}^N |G_{k v_1}|^2\prec \frac{\Im[m_N]}{\eta}, \quad 
    \sum_{k=1}^N |\tilde Y_{k\bm m j}|^2 \prec \sum_{k=1}^N |G_{k v_{2\ell}}|^2\prec \frac{\Im[m_N]}{\eta}.
\end{align*}
The claim \eqref{e:dermN} follows by taking $\bmk=\{k\}$ and $\sqrt{\Im[m_N]/\eta}X_{\bmk \bm m i}=\tilde X_{k\bm m i}, \sqrt{\Im[m_N]/\eta}Y_{\bmk \bm m j}=\tilde Y_{k\bm m j}$.

In the second case, say both $G_{k v_1}, G_{v_{2\ell}k}$ are in $\tilde X_{k\bm m i}$. Then $\tilde Y_{k\bm m j}=\tilde Y_{\bm m j}$ does not depend on the index $k$. Moreover, using  $|G_{v_s v_{s+1}}|\prec 1$ from Theorem \ref{locallaw} and Ward identity \eqref{e:WdI}
\begin{align*}
    \left|\sum_k \tilde X_{k\bm m i}\right|\prec \sum_{k=1}^N |G_{k v_1}G_{v_{2\ell}k}|\leq \frac{1}{2}\sum_{k=1}^N |G_{k v_1}|^2+|G_{v_{2\ell}k}|^2\prec \frac{\Im[m_N]}{\eta}.
\end{align*}
The claim \eqref{e:dermN} follows by taking $\bmk=\emptyset$ and $(\Im[m_N]/\eta) X_{\bmk \bm m i}=\sum_k \tilde X_{k\bm m i}, Y_{\bmk \bm m j}=\tilde Y_{\bm m j}$.
This finishes the proof of \eqref{e:dermN}.

Next we prove \eqref{e:derhij}. We recall the following estimates from \cite[Appendix A]{lee2021higher}, which follows from computing high moments
\begin{align}\label{e:sumbb}
\sum_{j=1}^N \left|h_{ij}^{2}-\frac{1}{N}\right|\prec 1,
\quad 
    \sum_{j=1}^N h_{ij}^{\al}\prec \frac{1}{q^{\al-2}},\quad \al\geq 2.
\end{align}
We recall from \eqref{e:tt} that $a_{2\ell}$ is a linear combination of terms in the form 
\begin{align}\label{e:ttcopy}
    w(\cF)= 
    \sum_{x_1,x_2,\cdots, x_{|V(\cF)|}}^*
    \frac{1}{N^{\theta(\cF)}}  \prod_{e\in E(\cF)} w(h_{\al_e \beta_e}; s_e), \quad w(h; s):= h^{s+1}-\frac{\bm1(s=1)}{N},
\end{align}
where  $\cF$ is a weighted forest as in Defintion \ref{d:forest}, $x_1, x_2,\cdots, x_{|V(\cF)|}$ enumerate the vertices of $\cF$ and  weights $s_e$ satifisfy $\sum_e (s_e+1)=2\ell$.

When we compute $\del^{\bm\beta} w(\cF)$,  the derivative $\del_{uv}$ may hit $w(h_{x y};s)$, which is nonzero if and only if $e=\{x,y\}=\{u,v\}$, and $\del_{uv}w(h_{uv};s)=(s+1) h_{uv}^s$. Therefore up to certain constant, the derivative $\del^{\bm\beta} w(\cF)$ is also in the form \eqref{e:tt}, i.e. a product over $e\in E(\cF)$. The difference is that we need to fix some edges to be $\{u,v\}$ with $u\neq v\in ij\bm m$ and no longer need to sum over these indices. 

We denote the vertex set which are not fixed to be $ij\bm m$ as $\bmv:=V(\cF)\setminus ij\bm m$.
After fixing some vertices to be $ij\bm m$ in $\cF$, the edge set of $E(\cF)$ decomposes into two sets
$\{i,j\}\cup E_1(\cF)=\{e=\{x,y\}: x, y\in ij\bm m\}$ and $E_2(\cF)=E(\cF)\setminus (E_1(\cF)\cup\{i,j\})$. Then for each edge $\{x,y\}\in E_2(\cF)$, at least one of $x,y$ is in $\bmv$. Then  $\del^{\bm\beta} w(\cF)$ is a linear combination of terms 
\begin{align}\begin{split}\label{e:wcF}
&\frac{h_{ij}^a}{N^{\theta(\cF)}}\prod_{\{x,y\}\in E_1(\cF)}h^{a_{xy}}_{xy}\left(h_{xy}^2-\frac{1}{N}\right)^{b_{xy}} \sum_{\bmv: ij\bm m\bmv\atop\text{distinct}}
\prod_{e\in E_2(\cF)}w(h_{e},s_{e}).
\end{split}\end{align}
We can further rewrite the first product in \eqref{e:wcF} as
\begin{align}\label{e:wcF2}
\prod_{\{i,x\}\in E_1(\cF)}h^{a_{ix}}_{ix}\left(h_{ix}^2-\frac{1}{N}\right)^{b_{ix}}
\prod_{\{j,x\}\in E_1(\cF)}h^{a_{jx}}_{jx}\left(h_{jx}^2-\frac{1}{N}\right)^{b_{jx}}
\prod_{\{x,y\}\in E_1(\cF)\atop x,y\not\in\{i,j\}}h^{a_{xy}}_{xy}\left(h_{xy}^2-\frac{1}{N}\right)^{b_{xy}}.
\end{align}
For each edge $e\in E_2(\cF)$ there are several possibilities: i) $e=\{i,x\}$ where $x\in \bmv$; iii) $e=\{j,x\}$ where $x\in \bmv$; iv) $e=\{x,y\}$ where $x,y\not\in \{i,j\}$. We can rewrite the last sum in \eqref{e:wcF} as
\begin{align}\begin{split}\label{e:wcF3}
& \sum_{\bmv: ij\bm m\bmv\atop\text{distinct}}\prod_{\{i,x\}\in E_2(\cF)}w(h_{ix},s_{\{ix\}})
\prod_{\{j,x\}\in E_2(\cF)}w(h_{jx},s_{\{jx\}})
\prod_{\{x,y\}\in E_2(\cF),\atop x,y\not\in \{i,j\}}w(h_{xy},s_{\{xy\}}).
\end{split}\end{align}

Using \eqref{e:wcF2}, \eqref{e:wcF3}, we can rewrite \eqref{e:wcF} in the following form 
\begin{align}\begin{split}\label{e:wcF4}
\frac{1}{N} h_{ij}^{a}\sum_{\bmv: ij\bm m\bmv\atop\text{distinct}}\tilde X_{i\bm m \bmv_0 \bmv_1}\tilde Y_{i\bm m \bmv_0 \bmv_2},
\end{split}\end{align}
where
\begin{align}\begin{split}\label{e:Xexp}
\tilde X_{i\bm m \bmv_0 \bmv_1}
= N^{-\frac{\theta(\cF)-1}{2}}
&\prod_{\{i,x\}\in E_1(\cF)}h^{a_{ix}}_{ix}\left(h_{ix}^2-\frac{1}{N}\right)^{b_{ix}}
\prod_{\{x,y\}\in E_1(\cF)\atop x,y\not\in\{i,j\}}h^{a_{xy}/2}_{xy}\left(h_{xy}^2-\frac{1}{N}\right)^{b_{xy}/2}\\
\times &\prod_{\{i,x\}\in E_2(\cF)}w(h_{ix},s_{\{ix\}})
\prod_{\{x,y\}\in E(\cF),\atop x,y\not\in \{i,j\}}\sqrt{w(h_{xy},s_{\{xy\}})}\\
\tilde Y_{j\bm m \bmv_0 \bmv_2}
= N^{-\frac{\theta(\cF)-1}{2}}
&\prod_{\{j,x\}\in E_1(\cF)}h^{a_{jx}}_{jx}\left(h_{jx}^2-\frac{1}{N}\right)^{b_{jx}}
\prod_{\{x,y\}\in E_1(\cF)\atop x,y\not\in\{i,j\}}h^{a_{xy}}_{xy}\left(h_{xy}^2-\frac{1}{N}\right)^{b_{xy}}\\&\prod_{\{j,x\}\in E_2(\cF)}w(h_{jx},s_{\{jx\}})
\prod_{\{x,y\}\in E(\cF),\atop x,y\not\in \{i,j\}}\sqrt{w(h_{xy},s_{\{xy\}})}.
\end{split}\end{align}
In \eqref{e:Xexp}, we have divided the vertex set $\bmv= \bmv_0\bmv_1 \bmv_2$, where $\bmv_1$ is the set of leaf vertices adjacent to $i$: $\bmv_1=\{x: x \text{ is a leaf vertex}, \{i,x\}\in E(\cF)\}$; $\bmv_2$ is the set of leaf vertices adjacent to $j$: $\bmv_2=\{x: x \text{ is a leaf vertex}, \{j,x\}\in E(\cF)\}$; and $\bmv_0$ is the set of remaining vertices. In this way, $\tilde X_{i\bm m \bmv_0 \bmv_1}$ does not depend on the vertex set $\bmv_2$ and $\tilde Y_{j\bm m \bmv_0 \bmv_2}$ does not depend on the vertex set $\bmv_1$.

Next we show that 
\begin{align}\label{e:smallb}
\sum_{\bmv_0}\left(\sum_{\bmv_1}|\tilde X_{i\bm m \bmv_0 \bmv_1}|\right)^2\prec 1,\quad
\sum_{\bmv_0}\left(\sum_{\bmv_2}|\tilde Y_{j\bm m \bmv_0 \bmv_2}|\right)^2\prec 1.
\end{align}

We prove \eqref{e:smallb} for the case that $\cF$ is a tree. The case that $\cF$ is a union of trees is the same. In this case, \eqref{e:Xexp} simplifies
\begin{align}\label{e:newXexp}
|\tilde X_{i\bm m \bmv_0 \bmv_1}|
\prec \prod_{\{i,x\}\in E_2(\cF),\atop x\in \bmv_1}|w(h_{ix},s_{\{ix\}})|\prod_{\{i,x\}\in E_2(\cF),\atop x\not\in \bmv_1 }|w(h_{ix},s_{\{ix\}})|
\prod_{\{x,y\}\in E_2(\cF),\atop x,y\not\in \{i,j\}}\sqrt{|w(h_{xy},s_{\{xy\}})|}.
\end{align}
The last two factors in \eqref{e:newXexp} do not depend on the indices $\bmv_1$. 

For any vertex $x$ in $\cF$, and a subset of its neighborhoods $\bmy=\{y_1, y_2,\cdots, y_n\}$. 
%the product in \eqref{e:tt} containing vertex $x$ is given by
%\begin{align}
%    \prod_{y: e=\{x,y\}\in E(\cF)}w(h_{xy},s_e),
%\end{align}
%where the product is over vertices $y$ adjacent to $x$.
Using \eqref{e:sumbb} and $|w(h_{xy},s_{xy})|\prec 1$ for any $\{x,y\}\in V(\cF)$, we have
\begin{align}\label{e:sumbb2}
    \sum_{x=1}^N \prod_{i=1}^n|w(h_{xy_i},s_{xy_i})|
    \leq \sum_{x=1}^N |w(h_{xy_1},s_{xy_1})|\prec 1,
\end{align}
provided $n\geq 1$.
Using the bound \eqref{e:sumbb2}, we can sum over the indices $\bmv_1$ in \eqref{e:newXexp}
\begin{align*}
\left(\sum_{\bmv_1}|\tilde X_{i\bm m \bmv_0 \bmv_1}|\right)^2\prec\prod_{\{i,x\}\in E_2(\cF),\atop x\not\in \bmv_1 }|w(h_{ix},s_{\{ix\}})|^2
\prod_{\{x,y\}\in E_2(\cF),\atop x,y\not\in \{i,j\}}|w(h_{xy},s_{\{xy\}})|.
\end{align*}
Then we can further sum over vertices $\bmv_0$, by repeatedly using \eqref{e:sumbb2}, 
\begin{align*}
\sum_{\bmv_0}\left(\sum_{\bmv_1}|\tilde X_{i\bm m \bmv_0 \bmv_1}|\right)^2\prec 1.
\end{align*}
This finishes the proof of the claim \eqref{e:Xexp}. 

The summation in \eqref{e:wcF} is over indices $\bmv=\bmv_0\bmv_1\bmv_2\in\qq{1,N}$ such that $ij\bm m \bmv$ are all distinct. We can first sum over $\bmv_0$ then $\bmv_1, \bmv_2$
\begin{align}\label{e:rewritesum}
\sum_{\bmv: ij\bm m\bmv\atop\text{distinct}}=\sum_{\bmv_0: ij\bm m\bmv_0\atop\text{distinct}}\sum_{\bmv_1\bmv_2: ij\bm m\bmv_0\bmv_1\bmv_2\atop\text{distinct}},
\end{align}
where in the second summation $\bmv_1, \bmv_2$ are distinct. By the inclusion-exclusion principle, we can rewrite the second summation in \eqref{e:rewritesum} as a linear combination of terms with $\bmv_1=\{\bmv_0', \bmv_1'\}$, $\bmv_2=\{\bmv_0', \bmv_2'\}$ where $\bmv_1', \bmv_2'$ may not be distinct:
\begin{align*}
\sum_{\bmv_0': ij\bm m\bmv_0\bmv_0'\atop\text{distinct}}\sum_{\bmv_1': ij\bm m\bmv_0\bmv_0'\bmv_1'\atop\text{distinct}}\sum_{\bmv_2': ij\bm m\bmv_0\bmv_0'\bmv_2'\atop\text{distinct}}.
\end{align*}
In this way, we can combine $\bmk=\{\bmv_0, \bmv_0'\}$, and conclude that the summation in \eqref{e:wcF} is a linear combination of terms in the form
\begin{align*}
\frac{1}{N} h_{ij}^{a}\sum_{\bmk: ij\bm m\bmk\atop\text{distinct}}X_{i\bm m\bmk}Y_{j\bm m\bmk},\quad X_{i\bm m\bmk}:=\sum_{\bmv_1': ij\bm m\bmk\bmv_1'\atop\text{distinct}}\tilde X_{i\bm m \bmk \bmv'_1},\quad Y_{i\bm m\bmk}:=\sum_{\bmv_2': ij\bm m\bmv''_0\bmv_2'\atop\text{distinct}}\tilde Y_{i\bm m \bmk \bmv'_2}. 
\end{align*}
Thanks to \eqref{e:smallb}
\begin{align*}
\sum_{\bmk}|X_{i\bm m\bmk}|^2\leq \sum_{\bmk}\left(\sum_{\bmv_1'}|\tilde X_{i\bm m \bmk \bmv'_1}|\right)^2
&=\sum_{\bmv_0\bmv_0'}\left(\sum_{\bmv_1'}|\tilde X_{i\bm m \bmv_0 \bmv_0' \bmv'_1}|\right)^2\\
&\leq \sum_{\bmv_0}\left(\sum_{\bmv_0'\bmv_1'}|\tilde X_{i\bm m \bmv_0 \bmv_0' \bmv'_1}|\right)^2
=\sum_{\bmv_0}\left(\sum_{\bmv_1}|\tilde X_{i\bm m \bmv_0\bmv_1}|\right)^2\prec 1.
\end{align*} 
And  we have the same estimate for $\sum_{\bmk}|Y_{j\bm m\bmk}|^2\prec 1$. This finishes the claim \eqref{e:derhij}.

\end{proof}

\begin{proposition}\label{p:dermN}
Adopt the assumptions in Proposition \ref{p:DSE}, and take $z=\tilde L+w$. Fix distinct indices $\bm m=\{m_1, m_2,\cdots, m_r\}$. We consider the differential operator $\del^{\bm \beta}$ with $\bm\beta=\{\beta_{uv}\}_{u,v\in ij\bm m}$
\begin{align*}
    \del^{\bm \beta}=\prod_{u,v\in \bm m}\del^{\beta_{uv}}_{uv},\quad D^{\bm \beta}=\prod_{u,v\in \bm m}D^{\beta_{uv}}_{uv}, \quad |\bm\beta|=\sum_{u,v\in \bm m} \beta_{uv}\geq 1.
\end{align*}
We have the following estimates
\begin{align}\label{e:impp}
    \del^{\bm \beta}z=\del^{\bm \beta} {\tilde L}\prec \frac{1}{N}, \quad \del^{\bm \beta} G_{ab}(z)=D^{\bm\beta}G_{ab}(z)+\OO_\prec \left(\frac{\Im[m_N]}{N\eta}\right),\quad \del^{\bm \beta}m_N(z)\prec \frac{\Im[m_N(z)]}{N\eta},
\end{align}
and
\begin{align}\label{e:dmN}
    &\del^{\bm \beta} m_N(z)=D^{\bm \beta} m_N(z)+\del_z m_N(z) \del^{\bm \beta} z+\OO_\prec\left(\frac{\Im[m_N]}{(N\eta)^2} \right),\\
\label{e:dP}\begin{split}
    &\del^{\bm \beta} P(z, m_N(z))=(D^{\bm \beta} m_N(z)+\del_z m_N(z) \del^{\bm \beta} z)P'\\
    &{\phantom{\del^{\bm \beta} P={}}}+((\del^{\bm \beta} z ) m_N(z)+ \sum_{k=1}^L (\del^{\bm \beta} a_{2\ell}) m_N^{2\ell}(z))+\OO_\prec\left(\frac{\Im[m_N(z)](|P'|+\Im[m_N(z)])}{(N\eta)^2} \right).\end{split}
\end{align}

\end{proposition}

\begin{proof}
From \eqref{e:derhij}, we have $\del^{\bm \beta}a_{2\ell}$, is a sum of terms in the form
\begin{align*}
    \left|\frac{1}{N}\sum_{\bmk:ij\bm m \bmk \atop \text{distinct}}h_{ij}^a   X_{ \bmk \bm m i}  Y_{\bmk\bm m j}\right|\prec \frac{1}{N}\sum_{\bm k}|X_{\bm k\bm m i}||Y_{\bm k \bm m j}|\prec \frac{1}{N}.
\end{align*}
Thus $|\del^{\bm \beta}a_{2\ell}|\prec 1/N$. For $\del^{\bm \beta} z$, by the chain rule we have it is a linear combination of terms in the form
\begin{align*}
   \left|\left(\left(\prod_{i=1}^{m}\del_{a_{2\ell_i}}\right) z\right)  \prod_{i=1}^{m} (\del^{\bm \beta_i} a_{2\ell_i})\right|\prec \frac{1}{N},\quad\quad \sum_{i=1}^m \bm\beta_i=\bm\beta,\quad m\geq 1,
\end{align*}
where we used Proposition \ref{p:minfty} that $\left(\prod_{i=1}^{m}\del_{a_{2\ell_i}}\right) z$ is bounded.

For the derivative $\del^{\bm \beta} G_{ab}$, each $\del_{uv}$ either it hits $G_{ab}$ or it hits $z$. The derivative $\del^{\bm \beta} G_{ab}$ is a sum of terms in the form
\begin{align*}
    (D^{\bm\beta'}\del_z^m G_{ab}) \prod_{i=1}^m (\del^{\bm\beta_i}z)
    \prec \frac{1}{N^m} |D^{\bm\beta'}\del_z^m G_{ab}|,\quad \bm\beta'+\sum_{i=1}^m\bm\beta_i=\bm\beta,\quad m\geq 1,
\end{align*}
where we used $|\del^{\beta_i}z|\prec 1/N$. 
We notice that for any Green's function term $G_{ij}$, its derivative satisfies
\begin{align}\label{e:dzGii}
|\del^m_z G_{ij}|=m!|(G^{m+1})_{ij}|
    =m!\left|\sum_\al \frac{u_\al(i)u_\al(j)}{(z-\la_\al)^{m+1}}\right|
    \prec \frac{1}{N\eta^{m-1}}\sum_\al \frac{1}{|z-\la_\al|^2}= \frac{\Im[m_N(z)]}{\eta^m},\quad m\geq 1,
\end{align}
where we used the delocalization of eigenvectors $\|u_\al\|_\infty\prec 1/\sqrt N$ from \eqref{e:delocal}. Since $D^{\bm \beta'}G_{ab}$ is a polynomial of Green's function entries, and $|G_{ij}|, |m_N|\prec 1$ from Theorem \ref{locallaw}, \eqref{e:dzGii} also implies for $m\geq 1$
\begin{align}\label{e:Gabdd}
    (D^{\bm\beta'}\del_z^m G_{ab}) \prod_{i=1}^m \del^{\bm\beta_i}z
    \prec \frac{1}{N^m} |\del_z^m D^{\bm\beta'} G_{ab}|\prec \frac{\Im[m_N]}{(N\eta)^{m}}.
\end{align}
The second statement in \eqref{e:impp} follows. 
For the bound of $\del^{\bm\beta} m_N$, we have 
\begin{align*}
    \del^{\bm\beta} m_N=\frac{1}{N}\sum_{i=1}^N\del^{\bm\beta} G_{ii}\prec \frac{\Im[m_N]}{N\eta}.
\end{align*}

Next for \eqref{e:dmN},  the derivative $\del^{\bm \beta} m_N$ is a sum of terms in the form
\begin{align}\label{e:mNbdd}
    (D^{\bm\beta'}\del_z^m m_N) \prod_{i=1}^m (\del^{\bm\beta_i}z)
    \prec \frac{1}{N^m} |D^{\bm\beta'}\del_z^m m_N|,\quad \bm\beta'+\sum_{i=1}^m\bm\beta_i=\bm\beta,\quad m\geq 1.
\end{align}
The first two terms in \eqref{e:dmN} correspond to the case when $m=0$ and $m=1, \bm \beta'=\emptyset$. When $m\geq 2$, using \eqref{e:Gabdd}, we have that \eqref{e:mNbdd} is bounded by $\OO_\prec(\Im[m_N]/(N\eta)^2)$. Next we estimate \eqref{e:mNbdd} for $m=1$ and $|\bm \beta'|\geq 1$,
\begin{align*}
    \frac{1}{N} |D^{\bm\beta'}\del_z m_N|= \frac{1}{N^2} |D^{\bm\beta'}\Tr G^2|
    =\frac{1}{N^2}\sum_{\bm\beta_1'+\bm\beta_2'=\bm\beta'} \left|\sum_{ij}D^{\bm\beta_1'}G_{ij} D^{\bm\beta_2'}G_{ij}\right|.
\end{align*}
There are two cases, 1) $|\bm\beta_1'|, |\bm\beta_2'|\geq 1$; 2) one of $|\bm\beta_1'|, |\bm\beta_2'|$ is zero. In the first case, we have that both $D^{\bm\beta_1'}G_{ij}, D^{\bm\beta_2'}G_{ij}$ are sums of terms in the form $G_{ix_1}G_{x_2 x_3}\cdots G_{x_{2\ell} j}$ where $\ell\geq 1$ and $x_1, x_2, \cdots, x_{2\ell}\in \bm m$. Then by the Ward indentity \eqref{e:WdI}, we have
\begin{align*}
\frac{1}{N^2}\left|\sum_{ij}D^{\bm\beta_1'}G_{ij} D^{\bm\beta_2'}G_{ij}\right|
\leq  \frac{1}{N^2}\sqrt{\sum_{ij}|D^{\bm\beta_1'}G_{ij}|^2\sum_{ij}|D^{\bm\beta_2'}G_{ij}|^2}\prec \frac{\Im[m_N]^2}{(N\eta)^2}.
\end{align*}
In the second case, say $|\bm\beta_1'|=0$. Since the $L_2$ norm of $G$ is bounded by $1/\eta$, we have
\begin{align*}
   \frac{1}{N^2} \left|\sum_{ij}G_{ij} G_{ix_1}G_{x_2 x_3}\cdots G_{x_{2\ell} j}\right|
   \leq \frac{1}{\eta N^2} \sqrt{\sum_{i}|G_{ix_1}|^2\sum_j|G_{x_2 x_3}\cdots G_{x_{2\ell} j}|^2 }
   \prec \frac{\Im[m_N]}{(N\eta)^2}.
\end{align*}
This finishes the proof of \eqref{e:dmN}. The claim \eqref{e:dP} follows from \eqref{e:dmN}.

\end{proof}

\begin{proof}[Proof of Proposition \ref{p:one-off}]
We first prove \eqref{e:derbb2}. The left hand side of \eqref{e:derbb2} is a linear combination of terms (with bounded coefficients) in the form
\begin{align*}
\left(\frac{1}{N}+\frac{\Im[m_N]}{N\eta}\right)I_1I_2\cdots I_t  P^{r-1-t_1}\bar P^{r-t_2}), \quad t_1+t_2=t,
\end{align*}
where for $1\leq s\leq t$,  $I_s=\del^{\bm\beta_s}P$ or $I_s=\del^{\bm\beta_s}\bar P$ with $|\bm\beta_s|\geq 1$.
Thanks to Proposition \ref{p:dermN}, $|\del^{\bm\beta}P|\prec \Im[m_N]/N\eta$. Therefore we can conclude that
\begin{align*}
    \left(\frac{1}{N}+\frac{\Im[m_N]}{N\eta}\right)I_1I_2\cdots I_t  P^{r-1-t_1}\bar P^{r-t_2})\lesssim \sum_{s\geq 1} \bE\left[\left(\frac{\Im[m_N]}{N\eta}\right)^s |P^{2r-s}|\right]\prec \Phi_r.
\end{align*}
This gives \eqref{e:derbb2}.

For \eqref{e:one-off}, there are two cases: 1) $\sum_{e\in E(\cF)}(p_e-s_e)\geq 1$, 2) $\sum_{e\in E(\cF)}(p_e-s_e)=0$. We first study the first case. As we will show later, the second case can be reduced to the first case.
We recall the derivatives of $P$ from \eqref{e:dP},
\begin{align}\begin{split}\label{e:DPP}
    \del^{\bm \beta} P&=(D^{\bm \beta} m_N+(\del_z m_N) \del^{\bm \beta} z)P'+((\del^{\bm \beta} z) m_N+ \sum_{\ell=1}^L(\del^{\bm \beta} a_{2\ell}) m_N^{2\ell})+\OO_\prec\left(\frac{\Im[m_N](|P'|+\Im[m_N])}{(N\eta)^2} \right)\\
    &=(D^{\bm \beta} m_N)P'+ \sum_{\ell=1}^L(\del^{\bm \beta} a_{2\ell}) ((\del_z m_N P'+m_N)\del_{a_{2\ell}}z+m_N^{2\ell})+\OO_\prec\left(\frac{\Im[m_N](|P'|+\Im[m_N])}{(N\eta)^2} \right).
\end{split} \end{align}

We  can rewrite \eqref{e:one-off} as a sum of terms in the form
\begin{align}\label{e:PII}
    \frac{1}{N^{r+2}}\sum^*_{ij\bm m} \bE[G_{ij}R_{ij\bm m} I_1 I_2\cdots I_t P^{r-1-t_1}\bar P^{r-t_2}],\quad t_1+t_2=t,
\end{align}
where $I_s$ for $1\leq s\leq t$ corresponds to terms in \eqref{e:DPP} defined in the following: 
\begin{enumerate}
   \item\label{e:1} From Proposition \ref{p:a2kstr}, the first term in \eqref{e:DPP} $(D^{\bm \beta} m_N) P'$ is a sum of terms in the form:
   \begin{align}\label{e:x11}
    \frac{\Im[m_N]P'}{N\eta}\sum_{\bmk_s}G_{ij}^{a_s} X_{\bmk_s \bm m i} Y_{\bmk_s \bm m j}.
\end{align}
   where $a_s\geq 1$ and $\bmk_s$ is some index set, and $\sum_{\bmk_s}|X_{\bmk_s \bm m i}^2|,\sum_{\bmk_s} |Y_{\bmk_s\bm m j}^2|=\OO_\prec(1)$.  We take $I_s$ to be \eqref{e:x11} or its complex conjugate.
   \item \label{e:2}From Proposition \ref{p:a2kstr}, the second term in \eqref{e:DPP} $\del^{\bm \beta} a_{2\ell} ((\del_z m_N P'+m_N)\del_{a_{2\ell}}z+m_N^{2\ell})$ is a sum of terms in the form:
   \begin{align}\label{e:x12}\begin{split}
   &\phantom{{}={}}\frac{((\del_z m_N P'+m_N)\del_{a_{2\ell}}z+m_N^{2\ell})}{N}\sum_{\bmk_s: ij\bm m\bmk_s \atop\text{distinct}}h_{ij}^{a_s}  X_{ \bmk_s \bm m i}  Y_{\bmk_s\bm m j}\\
   &=\OO\left(\frac{\Im[m_N]|P'|}{N\eta}+\frac{1}{N}\right)\sum_{\bmk_s: ij\bm m\bmk_s \atop\text{distinct}}h_{ij}^{a_s}  X_{ \bmk_s \bm m i}  Y_{\bmk_s\bm m j},
\end{split}\end{align}
where $a_s\geq 0$, $\bmk_s$ is an index set,  $\sum_{\bmk_s}|X_{\bmk_s \bm m i}^2|,\sum_{\bmk_s} |Y_{\bmk_s\bm m j}^2|=\OO_\prec(1)$. 
    We take $I_s$ to be \eqref{e:x12} or its complex conjugate.
   \item \label{e:3} $I_s$ is bounded by $\OO_\prec\left(\frac{\Im[m_N](|P'|+\Im[m_N])}{(N\eta)^2} \right)$ corresponding to the last term in \eqref{e:dP}.
\end{enumerate}

From the construction, we have $|I_s|\prec \Im[m_N]/N\eta$ for all $1\leq s\leq t$. 
If there is one $I_s$ corresponding to Item \ref{e:1} with $a_s\geq 1$, then
\begin{align*}
    |I_s|\prec |G_{ij}|\frac{\Im[m_N] |P'|}{N\eta}.
\end{align*}
And noticing $|P'|\prec 1$, we have 
\begin{align*}
 \frac{1}{N^{r+2}}\sum^*_{ij\bm m} \bE[G_{ij}R_{ij\bm m} I_1 I_2\cdots I_t P^{r-t_1-1}\bar P^{r-t_2}]
    &\prec \frac{1}{N^{2}}\sum^*_{ij} \bE\left[|G_{ij}|^{2}\left(\frac{\Im[m_N]}{N\eta}\right)^{t}|P^{2r-t-1}|\right]\\
    &=  \bE\left[\left(\frac{\Im[m_N]}{N\eta}\right)^{t+1}|P^{2r-t-1}|\right]\leq \bE[\Phi_r],
\end{align*}
where we used Ward identity \eqref{e:WdI}.

If  there is one $I_s$ corresponding to Item \ref{e:2} with $a_s\geq 1$, then
\begin{align*}
    |I_s|\prec \left(\frac{\Im[m_N]|P'|}{N\eta}+\frac{1}{N}\right)|h_{ij}|^{a_s}.
\end{align*}
Then by the Cauchy-Schwartz inequality, we have 
\begin{align*}
 &\phantom{{}={}}\frac{1}{N^{r+2}}\sum^*_{ij\bm m} \bE[G_{ij}R_{ij\bm m} I_1 I_2\cdots I_t P^{r-t_1-1}\bar P^{r-t_2}]\\
    &\prec \frac{1}{N^{2}}\left(\frac{\Im[m_N]|P'|}{N\eta}+\frac{1}{N}\right)\sum^*_{ij} \bE\left[|G_{ij}||h_{ij}|^{a_s}\left(\frac{\Im[m_N]}{N\eta}\right)^{t-1}|P^{2r-t-1}|\right]\\
    &\leq  \frac{1}{N^{2}}\left(\frac{\Im[m_N]|P'|}{N\eta}+\frac{1}{N}\right) \bE\left[\sqrt{\sum^*_{ij}|G_{ij}|^2\sum^*_{ij}|h_{ij}|^{2a_s}}\left(\frac{\Im[m_N]}{N\eta}\right)^{t-1}|P^{2r-t-1}|\right]\\
    &=  \bE\left[\frac{1}{N^{1/2}}\left(\frac{\Im[m_N]|P'|}{N\eta}+\frac{1}{N}\right)\left(\frac{\Im[m_N]}{N\eta}\right)^{t-1/2}|P^{2r-t-1}|\right]\leq \bE[\Phi_r],
\end{align*}
where in the last to second line we used 
\begin{align*}
    \sum_{ij}|G_{ij}|^2= \frac{ N \Im[m_N]}{\eta},\quad \sum_{ij}|h_{ij}|^{2a_s}\prec N.
\end{align*}

If there is one $I_s=\OO_\prec\left(\frac{\Im[m_N](|P'|+\Im[m_N])}{(N\eta)^2} \right)$ as in Item \ref{e:3}, we have
\begin{align*}
   &\phantom{{}={}}\frac{1}{N^{r+2}}\sum^*_{ij\bm m} \bE[G_{ij}R_{ij\bm m} I_1 I_2\cdots I_t P^{r-t_1-1}\bar P^{r-t_2}]\\
   &\prec \bE\left[\frac{\Im[m_N](|P'|+\Im[m_N])}{(N\eta)^2}\left(\frac{\Im[m_N]}{N\eta}\right)^{t-1}|P^{2r-t-1}|\right]\leq \bE[\Phi_r].
\end{align*}

In the rest, we can assume that each $I_s$ either corresponds to Item \ref{e:1} with $a_s=0$, or corresponds to Item \ref{e:2} with $a_s=0$. Say $I_1, I_2, \cdots, I_{t_3}$ correspond to Item \ref{e:1} with $a_s=0$, and $I_{t_3+1}, I_{t_3+2},\cdots, I_{t}$ correspond to Item \ref{e:2} with $a_s=0$, where $t_3+t_4=t$. Then we have
\begin{align}\begin{split}\label{e:GijFI}
    &\phantom{{}={}}\frac{1}{N^{r+2}}\sum^*_{ij\bm m} \bE[G_{ij}R_{ij\bm m} I_1 I_2\cdots I_t P^{r-t_1-1}\bar P^{r-t_2}]\\
    &=    \frac{1}{N^{r+2}}\sum_{\bmk'}\sum_{k''}^* \sum_{ij\bm m: ij\bm m \bmk''\atop \text{distinct}}\bE\left[G_{ij}X_{\bmk \bm m i}Y_{\bmk \bm m j}\OO\left(\left(\frac{\Im[m_N]P'}{N\eta}+\frac{1}{N}\right)^{t}\right) P^{r-t_1-1}\bar P^{r-t_2}\right].
\end{split}\end{align}
where $\bmk=\bmk'\cup \bmk''$ with $\bmk'=\cup_{1\leq s\leq t_3}\bmk_s$,  $\bmk''=\cup_{t_3+1\leq s\leq t}\bmk_s$, and 
\begin{align*}
    X_{\bm k \bm m i}=\prod_{s} X_{\bm k_s \bm m i},\quad  Y_{\bm k \bm mj}=\prod_{s} Y_{\bm k_s \bm m j}.
\end{align*}
They are bounded 
\begin{align}\label{e:Xbb}
    \sum_{\bm k}|X^2_{\bm k \bm m i}|=\sum_{\bm k}\prod_{s} |X^2_{\bm k_s \bm m i}|
    =\prod_s \sum_{\bm k_s} |X^2_{\bm k_s \bm m i}|\prec  1, \quad \sum_{\bm k}|Y^2_{\bm k\bm m j}|\prec 1.
\end{align}
Then, we use the norm of $G$ is bounded by $1/\eta$, and $|G_{ii}|\prec 1$, 
\begin{align}\begin{split}\label{e:L2est}
    &\left|\frac{1}{N^2}\sum_{ij: ij\bm m \bmk''\atop\text{distinct}}G_{ij}X_{\bmk \bm m i}Y_{\bmk \bm mj}\right|
    \leq \left|\frac{1}{N^2}\sum_{i: i\bm m \bmk''\atop\text{distinct}}\sum_{j: j\bm m \bmk''\atop\text{distinct}}G_{ij}X_{\bmk \bm m i}Y_{\bmk \bm mj}\right|\\
    &+ \left|\frac{1}{N^2}\sum_{i: i\bm m \bmk''\atop\text{distinct}}G_{ii}X_{\bmk \bm m i}Y_{\bmk \bm mi}\right|
\leq  \frac{1}{N^2\eta} \sqrt{\sum_{i}|X_{\bmk\bm m i}|^2 \sum_{j}|Y_{\bmk\bm m j}|^2}.
\end{split}\end{align}
Further, using the Cauchy-Schwartz inequality, and \eqref{e:Xbb} we have
\begin{align}\begin{split}\label{e:L2est}
    &\phantom{{}={}}\frac{1}{N^r}\left|\sum_{\bmk'}\sum_{k''}^* \sum_{\bm m: \bm m \bmk''\atop \text{distinct}}\frac{1}{N^2}\sum_{ij}G_{ij}X_{\bmk \bm m i}Y_{\bmk \bm m j}\right|
    \leq  \frac{1}{N^r} \sum_{\bm m}\frac{1}{N^2\eta}\sum_\bmk \sqrt{\sum_{i}|X_{\bmk \bm mi}|^2 \sum_{j}|Y_{\bmk\bm m j}|^2}\\
    &\leq \frac{1}{N^r} \sum_{\bm m}\frac{1}{N^2\eta}\sqrt{\sum_\bmk \sum_{j}|X_{\bmk \bm m i}|^2 \sum_\bmk \sum_{j}|Y_{\bmk\bm m j}|^2}
    \prec  \frac{1}{N^r}\frac{1}{N^2\eta}\sum_{\bm m}\sqrt{ \sum_{i}1  \sum_{j}1}=\frac{1}{N\eta}, 
\end{split}\end{align}
where we also used that $|\bm m|=r$. By plugging \eqref{e:L2est} into \eqref{e:GijFI}, we get
\begin{align*}
     \frac{1}{N^{r+2}}\sum^*_{ij\bm m} \bE[G_{ij}R_{ij\bm m} I_1 I_2\cdots I_t P^{r-t_1-1}\bar P^{r-t_2}]\prec \bE\left[\frac{1}{N\eta}\left(\frac{\Im[m_N] P'}{N\eta}+\frac{1}{N}\right)^t |P|^{2r-t-1}\right]\leq \bE[\Phi_r].
\end{align*}
This finishes the proof of the first case that $\sum_{e\in E(\cF)}(p_e-s_e)\geq 1$. 

For the second case, when $\sum_{e\in E(\cF)}(p_e-s_e)=0$. 
We use the identities $G_{ij}=\sum_{k\neq i}G_{ii}h_{ik}G_{kj}^{(i)}$, and $G_{kj}^{(i)}=G_{kj}-G_{ki}G_{ji}/G_{ii}$.
We denote
\begin{align*}
    U=R_{ij\bm m} P^{r-1}\bar P^{r}.
\end{align*}
Then by the cumulant expansion, we have
\begin{align}\begin{split}\label{e:reduce}
&\phantom{{}={}}\frac{1}{N^{r+2}}\sum_{ij\bm m}\bE\left[G_{ij}U
\right]
=\frac{1}{N^{r+2}}\sum^*_{ij\bm m}\bE\left[\sum_{k\neq i}h_{ik}G_{kj}^{(i)}G_{ii}U\right]\\
&=\sum_{p=1}^M\frac{\cC_{p+1}}{N^{r+3}q^{p-1}}\sum_{ij\bm m}\bE\left[\sum_{k\neq i}\del_{ik}^p (G_{kj}^{(i)}G_{ii} R_{ij\bm m} P^{r-1}\bar P^r) \right]+\OO_\prec(\bE[\Phi_r])\\
&=\sum_{p=1}^M\sum_s\frac{\cC_{p+1}}{N^{r+3}q^{p-1}}{p\choose s}\sum^*_{ij\bm m}\bE\left[\sum_{k\neq i}\del_{ik}^s (G_{kj}^{(i)}G_{ii} R_{ij\bm m})\del_{ik}^{p-s} (P^{r-1}\bar P^r) \right]+\OO_\prec(\bE[\Phi_r])\\
&=\sum_{p=1}^M\sum_s\frac{\cC_{p+1}}{N^{r+3}q^{p-1}}{p\choose s}\sum^*_{ij\bm m}\bE\left[\sum_{k\neq i}D_{ik}^s (G_{kj}^{(i)}G_{ii} R_{ij\bm m})\del_{ik}^{p-s}( P^{r-1}\bar P^r) \right]+\OO_\prec(\bE[\Phi_r]).
\end{split}\end{align}
where in the last line we replaced $\del_{ik}^s$ by $D_{ik}^s$ using \eqref{e:impp}, and the error term can be bounded by $\OO_\prec(\bE[\Phi_r])$ using \eqref{e:derbb2}. The sum of terms with $k\in j\bm m$ is also bounded 
\begin{align*}
&\phantom{{}={}}\frac{\cC_{p+1}}{N^{r+3}q^{p-1}}{p\choose s}\sum^*_{ij\bm m}\bE\left[\sum_{k\in j\bm m}D_{ik}^s (G_{kj}^{(i)}G_{ii} R_{ij\bm m})\del_{ik}^{p-s}( P^{r-1}\bar P^r) \right]\\
&\prec \frac{1}{N^{r+2}}\sum^*_{ij\bm m}\bE\left[\sum_{k\in j\bm m}\frac{1}{N}|\del_{ik}^{p-s}( P^{r-1}\bar P^r)| \right]\prec\frac{1}{N^{r+2}}\sum^*_{ij\bm m}\bE[\Phi_r]\leq\bE[\Phi_r].
\end{align*}
Therefore, we can further restrict the summation in \eqref{e:reduce} to $\sum_{ijk\bm m}^*$. 
Since $G_{kj}^{(i)}$ is independent of $h_{ik}$, we can further rewrite \eqref{e:reduce} as
\begin{align}\begin{split}\label{e:reduce2}
&=\sum_{p=1}^M\sum_s\frac{\cC_{p+1}}{N^{r+3}q^{p-1}}{p\choose s}\sum_{ij k\bm m}\bE\left[G_{kj}^{(i)} D_{ik}^s (G_{ii} R_{ij\bm m})\del_{ik}^{p-s}( P^{r-1}\bar P^r) \right]+\OO_\prec(\bE[\Phi_r])\\
&=\sum_{p=1}^M\sum_s\frac{\cC_{p+1}}{N^{r+3}q^{p-1}}{p\choose s}\sum_{ij k \bm m}\bE\left[\left(G_{kj}-\frac{G_{ki}G_{ji}}{G_{ii}}\right)D_{ik}^s (G_{ii} R_{ij\bm m})\del_{ik}^{p-s}( P^{r-1}\bar P^r) \right]+\OO_\prec(\bE[\Phi_r])\\
&=\sum_{p=1}^M\sum_s\frac{\cC_{p+1}}{N^{r+3}q^{p-1}}{p\choose s}\sum_{ij k \bm m}\bE\left[G_{kj}D_{ik}^s (G_{ii} R_{ij\bm m})\del_{ik}^{p-s}( P^{r-1}\bar P^r) \right]+\OO_\prec(\bE[\Phi_r]).
\end{split}\end{align} 
Here for the last line we used that if there are two off-diagonal Green's function entries, the sum is bounded by $\OO_\prec(\bE[\Phi_r])$ using Ward identity \eqref{e:WdI} and \eqref{e:derbb2}. The derivative $D_{ik}^s (G_{ii} R_{ij\bm m})$ is again a sum of monomials in Green's function entries. If $p>s$, the last line of \eqref{e:reduce} is in the form of \eqref{e:one-off} with $\sum_{e}(p_e-s_e)\geq 1$. Thus, from the discussion of the first case, they are bounded by $\OO_\prec(\bE[\Phi_r])$. The terms in \eqref{e:reduce2} with $p=s$ are given by 
\begin{align*}
\sum_{p=1}^M\frac{\cC_{p+1}}{N^{r+3}q^{p-1}}\sum^*_{ij k \bm m}\bE\left[G_{kj}D_{ik}^p (G_{ii} R_{ij\bm m})( P^{r-1}\bar P^r) \right].
\end{align*}
The derivative $D_{ik}^s (G_{ii} R_{ij\bm m})$ is a sum of monomials in Green's function entries. For these monomials containing at least one off-diagonal Green's function entry,  the total number of off-diagonal Green's function is at least two.  The sum is bounded by $\OO_\prec(\bE[\Phi_r])$. If there are monomials containing only diagonal Green's function entries, it is necessary that $p\geq 2$ ($D_{ik}G_{xx}=-2G_{xi}G_{kx}$ with $x\in ij\bm m$ contains at least one off-diagonal Green's function entry $G_{kx}$). These terms are again in the form of \eqref{e:one-off} with some extra $1/q$ factors. They are of higher order (recall from \eqref{e:ordero}). Then we can repeat the above procedure, until all the terms have order bigger than $M$. Then we can trivially bound them by $1/q^M\prec \Phi_r$ as in \eqref{e:ordero}. 

\end{proof}

\begin{proof}[Poof of Proposition \ref{p:diag}]
By the definition of the Green's function, we have
\begin{align}
\label{e:GHi'term}&1=-zG_{jj}+(HG)_{jj},\\
\label{e:GHiterm}&1=-zG_{ii}+(HG)_{ii}.
\end{align}
Multiplying \eqref{e:GHi'term} and \eqref{e:GHiterm} by $G_{ii}$ and $G_{jj}$ respectively, averaging over the indices, and then taking the difference, we get 
\begin{align}\label{e:iitokk}
    G_{ii}=m_N+\frac{1}{N}\sum_{j=1}^N (G_{ii} (HG)_{jj} - m_N(HG)_{ii}).
\end{align}
We will use the above relation \eqref{e:iitokk} to replace a copy of $G_{ii}$ on the lefthand side of \eqref{e:diag} to $m_N$.

%We will use the following relation to replace a copy of $G_{ii}$ on the lefthand side of \eqref{e:diag} to $G_{jj}$,
%\begin{align}\label{e:iitokk}
%    G_{ii}=G_{jj}+ (G_{ii} (HG)_{jj} - G_{jj} (HG)_{ii}).
%\end{align}
Denote \begin{align*}
U:=R_{\bm m} V, 
\quad V:=\left(\prod_{e\in E(\cF)}\del_{\al_e \beta_e}^{p_e-s_e}\right)(P^{r-1}\bar P^{r}).
\end{align*}
Then using \eqref{e:iitokk}, we can rewrite each term  on the righthand side of \eqref{e:nooff} (up to some constant) as
\begin{align}\begin{split}\label{e:diffc00}
\frac{1}{N^{r+1}}\sum^*_{i\bm m}\bE\left[G_{ii}^{\al}U\right]
&=\frac{1}{N^{r+1}}\sum^*_{i\bm m}\bE\left[m_NG_{ii}^{\al-1}U\right]\\
&+\frac{1}{N^{r+2}}\sum^*_{i\bm m}\sum_j \bE[(G_{ii} (HG)_{jj} - m_N (HG)_{ii})G_{ii}^{\al-1} U].
\end{split}\end{align}
It turns out the second term on the righthand side of \eqref{e:diffc00} is of order at least $1$ (recall from \eqref{e:ordero}). Using the cumulant expansion, we get
\begin{align}\begin{split}\label{e:diffc0}
    &\phantom{{}={}}\frac{1}{N^{r+2}}\sum^*_{i\bm m}\sum_j\bE[(G_{ii} (HG)_{jj} - m_N (HG)_{ii})G_{ii}^{\al-1}U]\\
    &= \frac{1}{N^{r+2}}\sum^*_{i\bm m}\sum_{jk}\sum_{p=1}^M \frac{\cC_{p+1}}{Nq^{p-1}}
    \bE[\del^{p}_{j k}(G_{ii}^\al G_{j k}U)-\del^{p}_{ik}(m_NG_{ik}G_{ii}^{\al-1}U)]+\OO_\prec \left(\frac{1}{q^{M}}\right)\\
    &=\frac{1}{N^{r+2}}\sum^*_{i \bm m}\sum_{j k }\sum_{p=1}^M\sum_{s} \frac{\cC_{p+1}}{Nq^{p-1}}{p\choose s}
    \bE[\del^{s}_{jk}(G^\al_{ii}G_{j k} R_{\bm m})\del^{p-s}_{j k} V)-\del^{s}_{ik}(m_NG_{ik}G_{ii}^{\al-1}R_{\bm m})\del^{p-s}_{ik} V)]+\OO_\prec \left(\frac{1}{q^{M}}\right).
%        &=\frac{1}{N^{r+2}}\sum^*_{i \bm m}\sum_{jk}\sum_{p=1}^M\sum_{s} \frac{\cC_{p+1}}{Nq^{p-1}}{p\choose s}.
%    \bE[D^{s}_{jk}(G^\al_{ii}G_{j k} U_1)\del^{p-s}_{j k} V)-D^{s}_{ik}(G_{jj}G_{ik}G_{ii}^{\al-1}U_1)\del^{p-s}_{ik} V)]+\OO_\prec \left(\Phi_r\right),
\end{split}\end{align}
%where in the last line, we used Proposition \ref{p:dermN} that replacing $\del^{s}_{jk}, \del^{s}_{ik}$ by $D^{s}_{jk}, D^{s}_{ik}$ gives an error $\Im[m_N]/N\eta$, namely $\del^{s}_{jk}(G^\al_{ii}G_{j k} U_1)=D^{s}_{jk}(G^\al_{ii}G_{j k} U_1)+\OO(\Im[m_N]/N\eta)$. Thanks to \eqref{e:onet}, the error term is bounded by $\OO_\prec(\Phi_r)$.
%\begin{align}
%\frac{1}{N^{r+2}}\sum_{ij k \bm m} 
%    \bE\left[\frac{\Im[m_N]}{N\eta}(|\del^{p-s}_{jk} V|+|\del^{p-s}_{ik} V|)\right]\prec \Phi_r.
%\end{align}

For the first term on the righthand side of \eqref{e:diffc0}, using Proposition \ref{p:one-off},  we can replace $\del_{jk}^s$ by $D_{jk}^s$, $\del^{s}_{jk}(G^\al_{ii}G_{j k} R_{\bm m})=D^{s}_{jk}(G^\al_{ii}G_{j k} R_{\bm m})+\OO(\Im[m_N]/N\eta)$.  Thanks to \eqref{e:derbb2}, the error term is bounded by $\OO_\prec(\bE[\Phi_r])$.
\begin{align*}
\frac{1}{N^{r+2}}\sum^*_{i \bm m} \sum_{jk}
   \frac{\Im[m_N]}{N\eta}|\del^{p-s}_{jk} V|\prec \Phi_r.
\end{align*}
Using \eqref{e:derbb2} again, for the sum when $j=k$, $j\in i\bm m$, or $k\in i\bm m$, the sum is bounded by
\begin{align*}
\frac{1}{N^{r+2}}\sum^*_{i \bm m} \sum_{jk: jki\bm m \atop\text{not distinct}}
    |\del^{p-s}_{jk} V|\prec \Phi_r.
\end{align*}
Thus we can restrict the summation to the case that $ijk\bm m$ are distinct. Finally, using Proposition \ref{p:one-off}, terms in $D^{s}_{jk}(G^\al_{ii}G_{j k} V)$ with at least one off-diagonal Green's function entries can be bounded by $\OO_\prec(\bE[\Phi_r])$, and
\begin{align*}
D^{s}_{jk}(G^\al_{ii}G_{j k} R_{\bm m})
=-\bm1(s \text{ is odd})s! G^\al_{ii}G_{jj}^{\frac{s+1}{2}}G_{kk}^{\frac{s+1}{2}} R_{\bm m}+\{\text{terms with off-diagonal entries}\}.
\end{align*}
Therefore, the leading terms in first term on the righthand side of \eqref{e:diffc0} are those which do not contain any off-diagonal Green's function terms,
\begin{align}\begin{split}\label{e:ftt1}
&\phantom{{}={}}\frac{1}{N^{r+2}}\sum^*_{ij k \bm m}\sum_{p=1}^M\sum_{s} \frac{\cC_{p+1}}{Nq^{p-1}}{p\choose s}
    \bE[D^{s}_{jk}(G^\al_{ii}G_{j k} R_{\bm m})\del^{p-s}_{j k} V)]\\
&=-\frac{1}{N^{r+2}}\sum^*_{ij k \bm m}\sum_{p=1}^M\sum_{s \text{ odd}} \frac{\cC_{p+1}s!}{Nq^{p-1}}{p\choose s}
    \bE[G^\al_{ii}G_{j j}^{\frac{s+1}{2}}G_{kk}^{\frac{s+1}{2}} R_{\bm m}\del^{p-s}_{j k} V)]
    +\OO_\prec(\bE[\Phi_r]).
\end{split}\end{align}
Similarly for the second term on the righthand side of \eqref{e:diffc0}, we have
\begin{align}\begin{split}\label{e:ftt2}
&\phantom{{}={}}-\frac{1}{N^{r+1}}\sum^*_{i \bm m}\sum_k\sum_{p=1}^M\sum_{s} \frac{\cC_{p+1}}{Nq^{p-1}}{p\choose s}
    \bE[\del^{s}_{ik}(m_NG_{ik}G_{ii}^{\al-1}R_{\bm m})\del^{p-s}_{ik} V)],\\
&= \frac{1}{N^{r+1}}\sum^*_{i k \bm m}\sum_{p=1}^M\sum_{s \text{odd}} \frac{\cC_{p+1}}{Nq^{p-1}}{p\choose s}s! {\frac{s-1}{2}+\al-1\choose \al-1}
    \bE[m_NG^{\al+\frac{s-1}{2}}_{ii}G^{\frac{s+1}{2}}_{kk}R_{\bm m}\del^{p-s}_{ik} V)]+\OO_\prec(\bE[\Phi_r]),
\end{split}\end{align}
where we used that 
\begin{align*}
D^{s}_{ik}(G_{ik}G_{ii}^{\al-1})
=-\bm1(s \text{ is odd})s!{\frac{s-1}{2}+\al-1\choose \al-1}G^{\al+\frac{s-1}{2}}_{ii}G^{\frac{s+1}{2}}_{kk}+\{\text{terms with diagonal entries}\}.
\end{align*}

By comparing \eqref{e:ftt1} and \eqref{e:ftt2}, the terms corresponding to $p=1,s=1$ cancel out:
\begin{align*}\begin{split}
-\frac{1}{N^{r+2}}\sum^*_{ij k \bm m} \frac{\cC_{2}}{N}
    \bE[G^\al_{ii}G_{j j}G_{kk} R_{\bm m}V)]
   &= -\frac{1}{N^{r+2}}\sum^*_{ik \bm m}\sum_j \frac{\cC_{2}}{N}
    \bE[G^\al_{ii}G_{j j}G_{kk} R_{\bm m}V)]+\OO_\prec(\bE[\Phi_r])\\
    &=-\frac{1}{N^{r+1}}\sum^*_{ik \bm m} \frac{\cC_{2}}{N}
    \bE[m_N G^\al_{ii}G_{kk} R_{\bm m}V)]+\OO_\prec(\bE[\Phi_r]).
\end{split}\end{align*}
Then the claim \eqref{e:hod} follows from combining \eqref{e:ftt1} and \eqref{e:ftt2}.

\end{proof}

\subsection{Proof of Theorem \ref{t:rigidity} }
\label{s:proofrigidity}

In this section we prove Theorem \ref{t:rigidity} by analyzing the high order moment estimates of $P(z, m_N(z))$ from Proposition \ref{p:DSE}. We recall the shifted spectral domain $\cD$ from \eqref{e:recalldefD}, and the following Proposition from \cite[Proposition 2.11]{huang2020transition}.

\begin{proposition}\label{p:stable}
There exists a constant $\varepsilon>0$ such that the following holds. Suppose that $\delta:\cD\rightarrow \mathbb{R}$ ($\cD$ is as defined in \eqref{e:recalldefD}) is a function so that 
\begin{align*}
|P({\tilde L}+w, m_N({\tilde L}+w))|\leq \delta(w),\quad w\in \cD.
\end{align*}
Suppose that $N^{-2}\leq \delta(w)\ll1$ for $w\in \cD$, that $\delta$ is Lipschitz continuous with Lipschitz constant $N$ and moreover that for each fixed $\kappa$ the function $\eta\mapsto \delta(\kappa+\ri\eta)$ is nonincreasing for $\eta>0$.   Then,
\begin{align*}
|m_N({\tilde L}+ w)-  \tilde m({\tilde L}+ w)|=\OO\left(\frac{\delta(w)}{\sqrt{|\kappa|+\eta+\delta(w)}}\right),
\end{align*}
where $\tilde m(z)$ is from Proposition \ref{p:minfty}, the implicit constant is independent of $N$.
\end{proposition}

Before proving Theorem \ref{t:rigidity}, we first prove a weaker estimate.
\begin{proposition}\label{p:firstbb}
Let $H$ be as in Definition \ref{asup} with $N^{\varepsilon}\leq q\lesssim N^{1/2}$.  Let $m_N(z)$ be the Stieltjes transform of its eigenvalue density, and $\tilde m(z)$ as defined in Proposition \ref{p:minfty}. Uniformly for any $z={\tilde L}+w$, $w=\kappa+\ri\eta\in \cD$ as defined in \eqref{e:recalldefD}, we have
\begin{align*}
|m_N(  z)-\tilde m(  z)|\prec \sqrt{|\kappa|+\eta}.
\end{align*}
\end{proposition} 
\begin{proof}
By Proposition \ref{p:minfty}, we have
\begin{align*}
\Im[  \tilde m(  z)]\asymp \Phi(w)\deq \left\{
\begin{array}{cc}
\sqrt{|\kappa|+\eta}, & \kappa\leq 0,\\
\eta/\sqrt{|\kappa|+\eta}, & \kappa\geq 0.
\end{array}
\right.
\end{align*}
and 
\begin{align*}
|\del_2 P(  z,   \tilde m(  z))|\asymp \sqrt{|\kappa|+\eta}.
\end{align*}
We denote 
\begin{align*}
\Lambda_N(w)\deq |m_N({\tilde L}+ w)-  \tilde m({\tilde L}+ w)|.
\end{align*}
%We roughly have the estimate 
%\begin{align*}\begin{split}
%\bE[|P(  z)|^{2r}]
%&\lesssim \bE\left[\frac{1}{q^3} \frac{\Im[m_N(  z)]}{N\eta}|\del_2 P(  z)||P(  z)|^{2r-2}\right]\\
%&+\bE\left[ \frac{\Im[m_N(  z)]}{N\eta}|P(  z)|^{2r-1}\right]+\bE\left[ \frac{\Im[m_N(  z)]}{(N\eta)^2}|\del_2 P(  z)||P(  z)|^{2r-2}\right].
%\end{split}\end{align*}
Then we have
\begin{align*}
\Im[m_N(  z)]\lesssim \Phi(w)+\Lambda_N(w),
\end{align*}
and by Proposition \ref{p:minfty}
\begin{align*}
\del_2 P(  z, m_N(  z))
=\del_2 P(  z, \tilde m(  z))+\OO(|m_N(  z)-\tilde m(  z)|)
=\OO(\sqrt{|\kappa|+\eta}+\Lambda_N(w)).
\end{align*}
By H{\" o}lder's inequality we obtain from Proposition \ref{p:DSE}, 
\begin{align}\label{e:momentbound}
\bE[|P(  z, m_N(  z))|^{2r}]
\prec\frac{1}{(N\eta)^{2r}}\bE\left[\Lambda_N(w)^{2r}+(|\kappa|+\eta)^{r/2}\left(\Phi(w)^r+\Lambda_N(w)^r\right)\right].
\end{align}
With overwhelming probability we have the following Taylor expansion,
\begin{align}
\nonumber P(  z, m_N(  z))
&=P(  z, \tilde m(  z))
+\del_2 P(  z, \tilde m(  z))(m_N(  z)-\tilde m(  z))+\frac{\del_2^2 P(  z, \tilde m(  z))+\OO(1)}{2}(m_N(  z)-\tilde m(  z))^2\\
&=\del_2 P(  z, \tilde m(  z))(m_N(  z)-\tilde m(  z))+(1+\oo(1))(m_N(  z)-\tilde m(  z))^2,
\label{e:taylorexpand}
\end{align}
where we used that $\del_2^2 P(  z, \tilde m(  z))=2+\OO(1/q)$ and $\Lambda_N(w)\ll1$ with overwhelming probability.   Rearranging the last equation and using the definition of $\Lambda_N(w)$, we have arrived at 
\begin{align}\label{e:Labda2bound}
\Lambda_N(w)^2\lesssim\Lambda_N(w)\sqrt{|\kappa|+\eta}+ |P(  z,m_N(  z))|, 
\end{align}
and thus 
\begin{align}\label{e:Lambda4r}
\bE[\Lambda_N(w)^{4r}]\lesssim(|\kappa|+\eta)^r\bE[\Lambda_N(w)^{2r}]+ \bE[|P(  z,m_N(  z))|^{2r}].
\end{align}
On the domain $\cD$, we have $1/N\eta\leq \sqrt{\kappa+\eta}$. We replace $\bE[|P(  z,m_N(  z))|^{2r}]$ in \eqref{e:Lambda4r} by \eqref{e:momentbound}, and get
\begin{align*}
\bE[\Lambda_N(w)^{4r}]\prec (|\kappa|+\eta)^{2r},
\end{align*}
It follows from Markov's inequality that $\Lambda_N(w)\prec \sqrt{|\kappa|+\eta}$.

\end{proof}

\begin{proof}[Proof of Theorem \ref{t:rigidity}]

We assume that there exists some deterministic control parameter $\Lambda(w)$ such that the prior estimate holds
\begin{align*}
|m_N({\tilde L}+ w)-  \tilde m({\tilde L}+ w)|\prec \Lambda(w)\lesssim \sqrt{|\kappa|+\eta}.
\end{align*}
%
%and Markov's inequality, we get
%\begin{align*}\begin{split}\label{e:estimateP}
%|P(  z, m_N(  z))|
%&\prec \sqrt{\left(\frac{1}{q^3}+\frac{1}{N\eta}\right) \frac{\Phi(w)+\Lambda(w)}{N\eta}(\sqrt{|\kappa|+\eta}+\Lambda(w))}+\frac{\Phi(w)+\Lambda(w)}{N\eta}.
%\end{split}
%\end{align*}
Since $\Phi(w)\gtrsim\sqrt{|\kappa|+\eta}$ and $\Lambda(w)\prec \sqrt{|\kappa|+\eta}$ from Proposition \ref{p:firstbb}, \eqref{e:momentbound}  combining with Markov's inequality leads to
\begin{align}\label{e:estimateP}
|P(  z, m_N(  z))|
&\prec \frac{1}{N\eta}\left((\Lambda(w)+\Phi(w))\sqrt{|\kappa|+\eta}\right)^{1/2}.
\end{align}

If $\kappa\geq 0$, then $\Phi(w)=\eta/ \sqrt{|\kappa|+\eta}$, and \eqref{e:estimateP} simplifies to
\begin{align}\label{e:boundoutside}
|P(  z, m_N(  z))|
\prec \frac{1}{N\eta^{1/2}}+\frac{(|\kappa|+\eta)^{1/4}\Lambda(w)^{1/2}}{N\eta}.
\end{align}
Thanks to Proposition \ref{p:stable}, by taking $\delta(w)$ the righthand side of \eqref{e:boundoutside} times $N^{\varepsilon}$ with arbitrarily small $\varepsilon$, we have
\begin{align}\begin{split}\label{e:outS0}
|m_N(  z)-\tilde m(  z)|\prec \frac{1}{\sqrt{|\kappa|+\eta}}\left( \frac{1}{N\eta^{1/2}}+\frac{(|\kappa|+\eta)^{1/4}\Lambda(w)^{1/2}}{N\eta}\right).
\end{split}\end{align}
%On the region $|\kappa|+\eta\gg 1/q^3N^{1/2}$, $N\eta\sqrt{|\kappa|+\eta}\gg1$ and $N\eta(|\kappa|+\eta)q^3\gg1$, the coefficients of $\Lambda^{1/2}$ is less than one. 
By iterating \eqref{e:outS0}, we get
\begin{align}\label{e:outS}
|m_N(  z)-\tilde m(  z)|\prec \frac{1}{\sqrt{|\kappa|+\eta}}\left(\frac{1}{N\eta^{1/2}}+\frac{1}{(N\eta)^2}\right).
\end{align}

If $\kappa\leq 0$, then $\Phi(w)=\sqrt{|\kappa|+\eta}$ and $\Lambda(w)\prec \sqrt{|\kappa|+\eta}$, \eqref{e:estimateP} simplifies to
\begin{align}\label{e:boundinside}
|P(  z,m_N(  z)|
\prec \frac{(|\kappa|+\eta)^{1/2}}{N\eta}.
\end{align}
It follows from Proposition \ref{p:stable}, by taking $\delta(z)$ the righthand side of \eqref{e:boundinside} times $N^{\varepsilon}$ with arbitrarily small $\varepsilon$, we have
\begin{align}\label{e:inS}
|m_N(  z)-\tilde m(  z)|\prec \frac{1}{N\eta}.
\end{align}
The claim \eqref{e:largeeig} follows from the estimates of the Stieltjes transform \eqref{e:outS} and \eqref{e:inS}, see \cite[Section 11]{erdHos2017dynamical}.

\end{proof}

\section{Edge statistics of $H(t)$} \label{sec:ht}
Let $H$ be as in Definition \ref{asup}. In this section we consider the Gaussian divisible ensemble
\beq\label{e:Ht}
H(t) := \e^{ - t/2}H + \left( 1 - \e^{ -t} \right)^{1/2} W,
\eeq
where $H(0)=H$ and $W$ is an independent GOE matrix. We denote the eigenvalues of $H(t)$ as $\la_1(t), \la_2(t),\cdots, \la_N(t)$, and the Stieltjes transform of its empirical eigenvalue distribution as
\begin{align*}
    m_t(z)=\frac{1}{N}\sum_{i=1}^N \frac{1}{\la_i(t)-z},\quad m_0(z)=m_N(z).
\end{align*}

Conditioning on $H$, the matrix ensemble $H(t)$ has the same law as the matrix Brownian motion starting from $H$ with each entry given by an Ornstein--Uhlenbeck process. The dynamic of the eigenvalues of the matrix Brownian motion is given by Dyson's Brownian motion
\begin{align}\label{e:DBM}
\rd \tilde{\la}_i(t)=\frac{\rd B_i(t)}{\sqrt N}+\frac{1}{N}\sum_{j:j\neq i}\frac{\rd t}{\tilde \la_i(t)-\tilde \la_j(t)}-\frac{\tilde \la_i(t)}{2}\rd t,
\end{align}
where for one time slice, $(\tilde \la_1(t), \tilde \la_2(t),\cdots, \tilde \la_N(t))$ has the same law as $(\la_1(t), \la_2(t),\cdots,  \la_N(t))$.

For sufficiently regular initial data, it has been proven in \cite{landon2017edge}, after short time the eigenvalue statistics at the spectral edge of \eqref{e:DBM} agree with GOE. A modified version of this theorem was proven in \cite{adhikari2020dyson}, which assumes that the initial data is sufficiently close to a nice profile. To use these results, we need to restrict  $H$ to a subset, on which the optimal rigidity holds.
We denote $\cA$ to be the set of sparse random matrices $H$, such that \eqref{e:outS} and \eqref{e:inS} hold at edges $\pm   {\tilde L}$:
\begin{align}\label{e:defA}
\cA\deq \{H: \text{\eqref{e:outS} and \eqref{e:inS} hold.} \}.
\end{align}
By Theorem \ref{t:rigidity}, we know that the event $\cA$ holds with probability $\bP(\cA)\geq 1-N^{-D}$ for any $D\geq 0$. 

We denote $\rhosc(x)$ the semicircle law which is the limit eigenvalue density of a Gaussian orthogonal ensemble $W$. The limit eigenvalue density of $(1-e^{-t})^{1/2}W$, is given by $(1-e^{-t})^{-1/2}\rhosc((1-e^{-t})^{-1/2}x)$, and the empirical eigenvalue distribution of $e^{-t/2}H$ concentrates around $e^{-t/2}  \tilde \rho(e^{t/2}x)$ (as defined in Proposition \ref{p:minfty}).
We denote the free convolution of $(1-e^{-t})^{-1/2}\rhosc((1-e^{-t})^{-1/2}x)$ and $e^{-t/2}  \tilde \rho(e^{-t/2}x)$ by $\tilde \rho_t$, and its Stieltjes transform by  $\td m_t$. Then $\tilde m_t$ satisfies the functional equations
\beq\label{e:defmt}
e^{-t/2}\td m_t(z) = \int \frac{ \tilde \rho (x)\rd x}{  x -  \xi_t (z)}=\tilde m( \xi_t(z)),\quad  \xi_t (z) := e^{t/2}z + e^{t/2}(1-\e^{-t} ) \td m_t (z).
\eeq
By the definition we have $\tilde m=\tilde m_0$. Recall from Proposition \ref{p:minfty}, $\tilde m_0(z)$ satisfies the functional equation 
\begin{align}\label{e:m0eq} 
    1+z\tilde m_0(z)+\tilde m_0(z)^2+Q(\tilde m_0(z))=0,
\end{align}
The next proposition states that $\tilde m_t$ satisfies a similar equation
\begin{proposition}\label{p:mtEQ}
Adapt the assumptions in Theorem \ref{t:rigidity}, and recall $\tilde m_t(z)$ from \eqref{e:defmt}. It is the Stieltjes transform of a measure $\tilde \rho_t$, which is the free convolution of $(1-e^{-t})^{-1/2}\rhosc((1-e^{-t})^{-1/2}x)$ and $e^{-t/2}  \tilde \rho(e^{-t/2}x)$ by $\tilde \rho_t$. The measure $\tilde \rho_t$ is symmetric and supported on $[-{\tilde L}_t, {\tilde L}_t]$. Moreover, $\tilde m_t(z)$ satisfies the following equation
\begin{align*}
    1+z\tilde m_t(z)+\tilde m_t^2(z)+Q(e^{-t/2}\tilde m_t(z))=0.
\end{align*}
\end{proposition}
\begin{proof}
By taking $z$ to be $ \xi_t(z)$ in \eqref{e:m0eq}, and using the relation \eqref{e:defmt}
\begin{align}\label{e:ppmt}
1+ \xi_t(z)e^{-t/2} \tilde m_t(z) +e^{-t}\tilde m^2_t(z)+Q(e^{-t/2}\tilde m_t(z))=0.
\end{align}
From the definition of $ \xi_t(z)$, we have
\begin{align*}
    \xi_t(z)e^{-t/2} +(e^{-t}-1)\tilde m_t(z)=z,
\end{align*}
and \eqref{e:ppmt} simplifies to 
\begin{align*}
    1+z \tilde m_t(z) +\tilde m^2_t(z)+Q(e^{-t/2}\tilde m_t(z))=0.
\end{align*}
For any $t\geq 0$,  The same argument as for Proposition \ref{p:minfty}, $\tilde m_t(z)$ is the Stieltjes transform of a measure $\tilde \rho_t$, which is symmetric and supported on $[-{\tilde L}_t, {\tilde L}_t]$. Moreover, it has square root behavior.
\end{proof}

We remark that $Q$ is a random polynomial which depends on certain averaged quantities of $h_{ij}$, so $\tilde L_t$ is also a random. But once we condition on $H$, both of them are deterministic. Next we prove the following theorem. It states that for time $t\gg N^{-1/3}$ the fluctuations of extreme eigenvalues of $H(t)$ conditioning on $H(0)\in \cA$ as in \eqref{e:defA} are given by the Tracy-Widom distribution.

\begin{theorem}\label{thm:htedge}
Let $H$ be as in Definition \ref{asup} with $N^{\varepsilon}\leq q\lesssim N^{1/2}$.  Conditioning on $H\in \cA$ as  in \eqref{e:defA}, let $H(t)$ be as in \eqref{e:Ht}, with eigenvalues denoted by $\la_1(t),\la_2(t),\cdots,\lambda_N(t)$, and $t=N^{-1/3+\fd}$.  Let $k\geq 1$ and $F : \rr^k \to \rr$ be a bounded test function with bounded derivatives.  There is a universal constant $\fc>0$ depending on $\fd$, it holds
\begin{align}\begin{split}\label{eqn:htedge2}
&\phantom{{}={}} \ee_{W}[ F (N^{2/3} ( \lambda_1(t) - {\tilde L}_t ), \cdots , N^{2/3} ( \lambda_k(t) - {\tilde L}_t) | H] \\
&= \ee_{GOE}[ F (N^{2/3} ( \mu_1 - 2), \cdots, N^{2/3} ( \mu_k - 2  ) ) +\OO\left(N^{-\fc}\right),
\end{split}\end{align}
where the expectation on the righthand side is with respect to a GOE matrix with eigenvalues denoted by  $\mu_1\geq \mu_2\geq \cdots\geq\mu_N$.
\end{theorem}

\begin{proof}
Take $\eta_*=N^{-2/3+\fd/2}$, and $  z={\tilde L}+ w$, where ${\tilde L}$ is from Proposition \ref{p:minfty} and $w=\kappa+\ri\eta\in \cD$ from \eqref{e:recalldefD}. For any $H\in \cA$, from the defining relations of $\cA$, i.e. \eqref{e:outS} and \eqref{e:inS}, we have 
\begin{align*}
|m_N(  z)-  \tilde m(  z)|\prec \frac{1}{N\eta},
\end{align*}
for $0\leq \kappa\leq 1$ and $\eta_*\leq \eta\leq 1$, and 
\begin{align*}
|m_N(  z)-  \tilde m(  z)|\prec \frac{1}{N\eta^{1/2}\sqrt{|\kappa|+\eta}}+\frac{1}{(N\eta)^2},
\end{align*}
for $-1\leq \kappa\leq 0$ and $\eta_*\leq \eta\leq 1$.
Moreover, \eqref{e:outS} also implies \eqref{e:largeeig} such that $\lambda_1 (0) -  {{\tilde L}}  \leq N^{-2/3+\fd/2}$.   Hence, $H$ is $\eta_*$-regular in the sense of \cite[Assumption 4.1]{adhikari2020dyson}, and the result of \cite[Theorem 6.1]{adhikari2020dyson} applies for $t=N^{-1/3+\fd}$ as above. This result gives the limiting distribution of the extreme eigenvalues of $H(t)$, and Theorem \ref{thm:htedge} follows.   
\end{proof}

It was also proven in \cite[Proposition 4.6]{adhikari2020dyson} that the Stieltjes transform $m_t(z)$ concentrates around $\tilde m_t(z)$. We collect the result in the following Proposition, which will be used in the next Section.
\begin{proposition}\label{p:mt}
Adapt the assumptions in Theorem \ref{thm:htedge}. Conditioning on $H\in \cA$ as  in \eqref{e:defA}, let $H(t)$ be as in \eqref{e:Ht}, with Stieltjes transform $m_t(z)$. Then for any $0\leq t\leq (\log N)^{-3}$, the following holds uniformly for $z=\tilde L_t+w$, with $w\in \cD$ from \eqref{e:recalldefD}
\begin{align}\label{e:mtrigid}
    |m_t(z)-\tilde m_t(z)|\prec \frac{1}{N\Im[w]}.
\end{align}
 \end{proposition}

\section{Comparison}\label{s:comparison}

We recall $H(t)$ from \eqref{e:Ht}. We denote the Stieltjes transform of its empirical eigenvalue density as $m_t(z)$, and the Stieltjes transform of the eigenvalue density of $H$ as $m_0(z)=m_N(z)$. 
%\subsection{Rigidity of Extreme Eigenvalues}
%\subsection{Comparison}
%The same as in \cite{Lee2016}, we use the interpolation, 
%\begin{align}
%H(t)=e^{-t/2}H+(1-e^{-t})^{1/2}W,
%\end{align}
%where $W$ is a Gaussian orthogonal ensemble.
In this section we prove the following theorem, which states that for $t\ll N^{-1/3}q$, the rescaled extreme eigenvalues of $H$ and $H(t)$  have the same distribution. 
Then Theorem \ref{thm:Tracy-Widom} follows from combining Theorem \ref{thm:htedge} and Theorem \ref{thm:comp}.

\begin{theorem}\label{thm:comp}
Let $H$ be as in Definition \ref{asup} with $N^{\varepsilon}\leq q\lesssim N^{1/2}$,  and $H(t)$ be as in \eqref{e:Ht}, with eigenvalues denoted by $\la_1(t),\la_2(t),\cdots,\lambda_N(t)$, and $t=N^{-1/3+\fd}$ with $\fd\leq \varepsilon/20$.  Fix $k\geq 1$ and numbers $s_1,s_2,\cdots, s_k$, there is a universal constant $\fc>0$ so that,
\begin{align}\begin{split}\label{e:comp1}
&\phantom{{}={}} \bP_{H}\left( N^{2/3} ( \lambda_i(0) - {\tilde L} )\geq s_i,1\leq i\leq k \right)\\
&= \bP_{H(t)}\left( N^{2/3} ( \lambda_i(t) - {\tilde L}_t )\geq s_i,1\leq i\leq k \right) +\OO\left(N^{-\fc}\right),
\end{split}\end{align}
where ${\tilde L}_t$ is as defined in Proposition \ref{p:mtEQ}. 
The analogous statement holds for the smallest eigenvalues.
\end{theorem}

Theorem \ref{thm:comp} is a consequence of the following Green's function comparison result. 
\begin{proposition}\label{p:comp}
Adapt the assumptions in Theorem \ref{thm:comp}. We 
fix $\fc>0$, $E_1,E_2,\cdots, E_k=\OO( N^{-2/3})$, $\eta_0=N^{-2/3-\fc}$ and $F:\bR^k\mapsto \bR$ a bounded test function with bounded derivatives. For $t\ll1$ we have
\begin{align}\begin{split}\label{e:comparison} 
&\phantom{{}={}}\bE_{H}\left[F\left(\left\{\Im\left[N\int_{E_i}^{N^{-2/3+\fc}} m_N( {\tilde L}+y+\ri \eta_0)\right]\right\}_{i=1}^k\right)\right]\\
&=\bE_{H(t)}\left[\left\{F\left(\Im\left[N\int_{E_i}^{N^{-2/3+\fc}} m_t( {\tilde L}_t+y+\ri \eta_0)\right]\right\}_{i=1}^k\right)\right]+\OO\left( N^{10\fc}\left(\frac{N^{1/3}t}{q}\right)\right).
\end{split}\end{align}
\end{proposition}

Next we prove Theorem \ref{thm:comp} using Proposition \ref{p:comp} as an input. The proof of Proposition \ref{p:comp} will occupy the remaining of this section. 
\begin{proof}[Proof of Theorem \ref{thm:comp}]
We need to first introduce some notations. For any $E\in \bR$, we define 
\begin{align*}
\cN_t(E)\deq |\{i:\la_i(t)\geq {\tilde L}_t+E\}|,
\end{align*}
and we write $\cN_0(E)$ as $\cN(E)$. We fix $\fc>0$, and take $\ell=N^{-2/3-\fc/3}$ and $\eta_0=N^{-2/3-\fc}$. Both are smaller than $N^{-2/3}$. Then with overwhelming probability, from \eqref{e:largeeig}, we have that $\la_1(t)\leq {\tilde L}_t+N^{-2/3+\fc}$. We define:
\begin{align*}
\chi_E(x)={\bf 1}_{[E, N^{-2/3+\fc}]}(x-{\tilde L}_t),\quad \theta_\eta(x)\deq \frac{\eta}{\pi(x^2+\eta^2)}=\frac{1}{\pi}\Im\frac{1}{x+\ri \eta}.
\end{align*}
From the same argument as in \cite[Lemma 2.7]{MR3034787}, we get that 
\begin{align}\label{e:resolventexp}
\Tr(\chi_{E+\ell}*\theta_\eta)(H(t))-N^{-\fc/9}\leq \cN_t(E)\leq \Tr(\chi_{E-\ell}*\theta_\eta)(H(t))+N^{-\fc/9},
\end{align}
hold with overwhelming probability. Let $K_i: \bR\mapsto[0,1]$ be a monotonic smooth function satisfying,
\begin{align*}
K_i(x)=\left\{\begin{array}{cc}
0 & x\leq i-2/3,\\
1& x\geq i-1/3.
\end{array}\right.
\end{align*}
We have that ${\bf 1}_{\cN_t(E)\geq i}=K_i(\cN_t(E))$, and since $K_i$ is monotonically increasing, and so 
\begin{align*}
K_i\left(\Tr(\chi_{E+\ell}*\theta_\eta)(H(t))\right)+\OO(N^{-\fc/9})\leq {\bf 1}_{\cN_t(E)\geq i}\leq K_i\left(\Tr(\chi_{E-\ell}*\theta_\eta)(H(t))\right)+\OO(N^{-\fc/9}),
\end{align*}

In this way we can express the locations of eigenvalues in terms of the integrals of the Stieltjes transform of the empirical eigenvalue densities.  We have,
\begin{align}\begin{split}\label{e:sandwich}
&\phantom{{}={}}\bE_{H(t)}\left[\prod_{i=1}^k K_i\left(\Im\left[\frac{N}{\pi}\int_{s_iN^{-2/3}+\ell}^{N^{-2/3+\fc}} m_t({\tilde L}_t+y+\ri \eta)\right]\rd y\right)\right]+\OO\left(N^{-\fc/9}\right)\\
&\leq \bP_{H(t)}\left( N^{2/3} ( \lambda_i(t) - {\tilde L}_t  )\geq s_i,1\leq i\leq k \right) =\bE\left[\prod_{i=1}^k{\bf 1}_{\cN_t(s_iN^{-2/3})\geq i}\right] \\
&\leq \bE_{H(t)}\left[\prod_{i=1}^k K_i\left(\Im\left[\frac{N}{\pi}\int_{s_iN^{-2/3}-\ell}^{N^{-2/3+\fc}} m_t( {\tilde L}_t+y+\ri \eta)\right]\rd y\right)\right]+\OO\left(N^{-\fc/9}\right).\end{split}
\end{align}

Since $q\geq N^\varepsilon$ and $t=N^{-1/3+\fd}$, we can take $\fc$  and $\fd$ smaller than $\varepsilon/20$, and then the error terms in \eqref{e:comparison} are of order $\OO(N^{-\fc})$.  By combining \eqref{e:sandwich} and \eqref{e:comparison}, we get
\begin{align*}\begin{split}
&\leq \bP_{H(t)}\left( N^{2/3} ( \lambda_i(t) - {\tilde L}_t  )\geq s_i+2N^{2/3}\ell,1\leq i\leq k \right)+\OO(N^{-\fc/9})\\
&\leq \bE_{H(t)}\left[\prod_{i=1}^k K_i\left(\Im\left[\frac{N}{\pi}\int_{s_iN^{-2/3}+\ell}^{N^{-2/3+\fc}} m_t( {\tilde L}_t+y+\ri \eta)\right]\rd y\right)\right]+\OO\left(N^{-\fc/9}\right)\\
&\leq \bP_{H}\left( N^{2/3} ( \lambda_i(0) - {\tilde L}  )\geq s_i,1\leq i\leq k \right) \\
&\leq \bE_{H(t)}\left[\prod_{i=1}^k K_i\left(\Im\left[\frac{N}{\pi}\int_{s_iN^{-2/3}-\ell}^{N^{-2/3+\fc}} m_t( {\tilde L}_t+y+\ri \eta)\right]\rd y\right)\right]+\OO\left(N^{-\fc/9}\right)\\
&\leq \bP_{H(t)}\left( N^{2/3} ( \lambda_i(t) - {\tilde L}_t  )\geq s_i-2N^{2/3}\ell,1\leq i\leq k \right)+\OO(N^{-\fc/9}).
\end{split}
\end{align*}
Since $N^{2/3}\ell=N^{-\fc/3}\ll1$, \eqref{e:comp1} follows. 
\end{proof}

For simplicity of notation we only prove Proposition \ref{p:comp} in the case $k=1$.  The general case can be proved in the same way. Let,
\begin{align*}
X_t:=X_t(H(t), \tilde L_t)=\Im\left[N\int_{E}^{N^{-2/3+\fc}}m_t( {\tilde L}_t +y+\ri \eta_0) \rd y\right].
\end{align*}
We prove the $k=1$ case of \eqref{e:comparison} 
\begin{align}\label{e:onecomp}
|\bE[F(X_t)]-\bE[F(X_0)]|\prec N^{10\fc}\left(\frac{N^{1/3}t}{q}\right).
\end{align}

In the rest of this section, we recall $H(t)$ from \eqref{e:Ht}
\beq\label{e:Htcopy}
H(t) := \e^{ - t/2}H + \left( 1 - \e^{ -t} \right)^{1/2} W.
\eeq
We denote the Green's function of $H(t)$ by $G(z;t)=(H(t)-z)^{-1}$. If the context is clear, we will simply write $G(z;t)$ as $G$ or $G(z)$. We write the derivatives $\del_{ij}=\del_{h_{ij}}$. For the remaining of this section, we will take $z=\tilde L_t+w$, with $w=y+\ri\eta$, $|y|\leq N^{-2/3+\fc}$ and $\eta\geq N^{-2/3-\fc}$. 
Then $z$ depends on $h_{ij}$ through $\tilde L_t$. From Theorem \ref{locallaw}, we have that $|G_{ij}(z,t)|\prec 1$. 
Both $z$ and $\tilde L_t$ are independent of $W$. The derivative $\del_{ij}$ in $\del_{ij}G(z;t)$ may hit $G$ or $z$. We introduce the notation 
$D_{ij}G(z;t):=\del_{h_{ij}(t)}G(z;t)=-G(z;t)(E_{ij}+E_{ji})G(z;t)$, where the derivative does not hit $z$, and $E_{ij}$ is the $N\times N$ matrix whose $(i,j)$-th entry is one and other entries are zero.  With this notation, we have 
\begin{align*}
\del_{ij}G(z;t)=\del_{h_{ij}}G(z;t)=e^{-t/2}D_{ij}G(z;t)+(\del_{h_{ij}}\tilde L_t)\del_z G(z;t).
\end{align*}
In the rest of this section, we prove the following proposition on the time derivative of $\bE[F(X_t)]$. The claim \eqref{e:onecomp} follows from plugging \eqref{e:dFX0} into \eqref{e:dFXt} and integrating from time $0$ to $t$.
\begin{proposition}\label{p:DSEt}
Adapt the assumptions in Proposition \ref{p:comp}.
Let $w=y+\ri\eta_0$ with $\eta_0=N^{-2/3-\fc}$ and $z=\tilde L_t+w$, we have the following estimates
\begin{align}\begin{split}\label{e:dFXt}
     \frac{\rd}{\rd t}\bE[F(X_t)]&= \left.\sum_{p=2}^M \frac{e^{-t/2}\cC_{p}}{2Nq^{p-1}} \sum_{ij}\bE[\del_{ij}^p( G_{ij}(z;t) F'(X_t))]+\frac{N}{2}\bE[ Q(e^{-t/2}m_t(z))F'(X_t)]\right|_{y=E}^{y=N^{-2/3+\fc}}\\
     &+\OO_\prec\left( \frac{N^{10\fc+1/3}}{q}\right),
\end{split}\end{align}
and uniformly for any $ |y|\leq N^{-2/3+\fc}$, ,
\begin{align}\label{e:dFX0}
 \frac{1}{N}\sum_{p=2}^M \frac{e^{-t/2}\cC_{p}}{N q^{p-1}} \sum_{ij}\bE[\del_{ij}^p( G_{ij}(z;t) F'(X_t))]+\bE[ Q(e^{-t/2}m_t(z))F'(X_t)]\prec \frac{N^{8\fc}}{N^{2/3}q}.
\end{align}
\end{proposition}

\subsection{Proof of \eqref{e:dFXt}}\label{s:p1}
We compute the time derivative of $\bE[F(X_t)]$
\begin{align}\begin{split}\label{e:dFt}
&\phantom{{}={}}\frac{\rd}{\rd t}\bE[F(X_t)]
=\bE\left[F'(X_t)\frac{\rd X_t}{\rd t}\right]\\
&=\bE\left[F'(X_t)\Im\int_{E}^{N^{-2/3+\fc}}\left.\left(\sum_{ijk}\dot{h}_{ij}(t) D_{ij}G_{kk}(z,t)+N\del_t{{\tilde L}}_t\del_z m_t(z)\right)\right|_{z=\tilde L_t+y+\ri \eta_0}\rd y\right],
\end{split}\end{align}
where from the definition \eqref{e:Htcopy} of $H(t)$,
\begin{align*}
\dot{h}_{ij}(t)=-\frac{1}{2}e^{-t/2}h_{ij}+\frac{e^{-t}}{2\sqrt{1-e^{-t}}}w_{ij}.
\end{align*}

%By the cumulant expansion formula, (we get $p=1$ term cancels)  \hty{remove the last  sentence in (... )}:
%\begin{align}
%\bE[h_{ij}(t)F(H(t))]
%&=\bE[(e^{-t/2}h_{ij}+\sqrt{1-e^{-t}}w_{ij})F(e^{-t/2}H+\sqrt{1-e^{-t}}W)]\\
%&=\frac{1}{N}\bE[\del_{ij}F(H(t))]+\sum_{p\geq 2} \frac{e^{-(p+1)t/2}\cC_{p}}{2Nq^{p-1}}\sum_{ijk}\bE[\del_{ij}^p F(H(t))].
%\end{align}
%where $\del_{ij}$ is with respect to $h_{ij}(t)$. 

The key to understand the righthand side of \eqref{e:dFt} is to compute the time derivative of ${\tilde L}_t$, which is given by  the following Proposition.
\begin{proposition}\label{p:Ltder}
Adapt the assumptions in Proposition \ref{p:comp}. 
We have the following estimate with high probability, uniformly for any $w=y+\ri\eta_0$ with $ |y|\leq N^{-2/3+\fc}, \eta_0=N^{-2/3-\fc}$
\begin{align}\label{e:Ltder}
\left|\del_t {\tilde L}_t -\frac{1}{2}\del_m Q(e^{-t/2}m_t({\tilde L}_t+w))\right|\prec \frac{N^{4\fc}}{N^{1/3}q},
\end{align}
where the polynomial $Q(m)$ is from \eqref{e:defQ}, and ${\tilde L}_t$ is constructed in Proposition \ref{p:mtEQ}.
\end{proposition}
Before proving Proposition \ref{p:Ltder}, we first state some useful estimates, which will be used repeatedly in the rest of this section. Their proofs are postponed to the next section. 

\begin{proposition}\label{p:Ltder2}
Adapt the assumptions in Proposition \ref{p:comp}.  Uniformly for any $w=y+\ri\eta$ with $ |y|\leq N^{-2/3+\fc}, \eta\geq N^{-2/3-\fc}$, we have the following estimates \begin{align}\label{e:dermtha}
   |\Im[m_t({\tilde L}_t+w)]|\prec N^{-1/3+3\fc}, \quad |\del_w m_t({\tilde L}_t+w)|\prec N^{1/3+4\fc}, \quad |\del_wG_{ab}(\tilde L_t+w)|\prec N^{1/3+4\fc}.
\end{align}
Fix distinct indices $i,j, \bm m=\{m_1, m_2,\cdots, m_r\}$, we consider the differential operators 
\begin{align}\label{e:Db}
    \del^{\bm \beta}=\prod_{u,v\in ij\bm m}\del^{\beta_{uv}}_{uv},\quad  D^{\bm \beta}=\prod_{u,v\in ij\bm m}D^{\beta_{uv}}_{uv}, \quad |\bm\beta|=\sum_{u,v\in ij\bm m} \beta_{uv}\geq 1.
\end{align}
The following holds
\begin{enumerate}
\item The  derivatives of $\tilde L_t$ satisfy: $\del^{\bm \beta} {\tilde L_t}\prec 1/N$. 
In the special case $\beta_{ij}+\beta_{ji}=1$, we have slightly stronger estimate $\del^{\bm \beta} {\tilde L_t}\prec |h_{ij}|/N$. 

\item The derivatives of the Green's function $G$ and Stieltjes transform $m_N$ satisfy
\begin{align}\label{e:impp2}
    \del^{\bm \beta} G_{ab}(\tilde L_t+w)=e^{-|\bm\beta|t/2}D^{\bm\beta}G_{ab}(\tilde L_t+w)+\OO_\prec \left(N^{-2/3+4\fc}\right)\prec 1,
\end{align}
and 
\begin{align}\label{e:impp3}
    |\del^{\bm \beta}m_N(\tilde L_t+w)|\prec N^{-2/3+4\fc},\quad |\del^{\bm \beta}X_t|\prec N^{-1/3+5\fc},\quad |\del^{\bm \beta} F'(X_t)|\prec N^{-1/3+5\fc}.
\end{align}
\item For any monomial $R_{ij\bm m}$ of the Green's function entries $G(\tilde L_t+w, t)$  as in Definition \ref{d:evaluation}, if $\beta_{ij}=\beta_{ji}=0$, we have the estimates
\begin{align}\begin{split}\label{e:dtLt}
&\phantom{{}={}}\frac{1}{N}\sum_{ij:ij\bm m\atop \text{distinct}}\bE\left[ \del_{w_{ij}}(R_{ij\bm m}\del^{\bm\beta}(F'(X_t)))\right]\\
&= \frac{e^{t/2}\sqrt{1-e^{-t}}}{N}\sum_{ij:ij\bm m\atop \text{distinct}}\bE\left[ \del_{ij}(R_{ij\bm m}\del^{\bm\beta}(F'(X_t)))\right]+\OO_\prec\left(\frac{N^{1/3+5\fc}\sqrt{1-e^{-t}}}{q}\right).
\end{split}\end{align}

\end{enumerate}
\end{proposition}

As an easy consequence of Proposition \ref{p:Ltder2}, for any $w=y+\ri\eta$ with $|y|\leq N^{-2/3+\fc}, \eta\geq N^{-2/3-\fc}$ by Ward identity \eqref{e:WdI}
\begin{align}\begin{split}\label{e:derofGG}
&\frac{1}{N^2}\sum_{ij}|G_{ij}(\tilde L_t+w)|^2= \frac{\Im[  m_t(\tilde L_t+w)]}{N\eta}\prec N^{-2/3+4\fc},\\
& \frac{1}{N^2}\sum_{ij}|G_{ij}(\tilde L_t+w)|\leq  \sqrt{\frac{\Im[  m_t(  \tilde L_t+w)]}{N\eta}}\prec N^{-1/3+2\fc}.
\end{split}\end{align}

\begin{proof}[Proof of Proposition \ref{p:Ltder}]
The spectral edge ${\tilde L}_t$ is characterized by
\begin{align}\label{e:dettLt}
    {\tilde L}_t=-\frac{1}{\zeta_t}-\zeta_t-\frac{ Q_t(\zeta_t)}{\zeta_t},\quad \left.\del_m \left(-\frac{1}{m}-m-\frac{ Q_t(m)}{m}\right)\right|_{m=\zeta_t}=0,
\end{align}
where $\zeta_t=\tilde m_t({\tilde L}_t)$, and $\tilde m_t$ is the solution of $1+z\tilde m_t+\tilde m_t^2+Q_t(\tilde m_t)=0$. 
By taking time derivative on both sides of \eqref{e:dettLt}, we get
\begin{align}\begin{split}\label{e:mtbb00}
    \del_t {\tilde L}_t&=\left.\del_m \left(-\frac{1}{m}-m-\frac{ Q_t(m)}{m}\right)\right|_{m=\zeta_t}\del_t \zeta_t -\frac{(\del_t Q_t)(\zeta_t)}{\zeta_t}\\
    &=-\frac{\del_t Q_t(\zeta_t)}{\zeta_t}
    =\sum_{\ell=1}^L \ell a_{2\ell}e^{-t \ell}\zeta_t^{2\ell-1}=\frac{1}{2}\del_m Q_t(m)|_{m=\zeta_t}.
\end{split}\end{align}
Let $w'=y+N^{-2/3+\fc}\ri$, then it is easy to see that $w'\in \cD$ (recall from \ref{e:recalldefD}).
By \eqref{e:dermtha} and the optimal rigidity estimates \eqref{e:mtrigid}
\begin{align*}
    |m_t({\tilde L}_t+w)-\zeta_t|
    &\leq |m_t({\tilde L}_t+w')-m_t({\tilde L}_t+w)|+
    |m_t({\tilde L}_t+w')-\tilde m_t(\tilde L_t)|\\
    &\leq|w'-w| N^{1/3+4\fc}+|\tilde m_t({\tilde L}_t+w')-\tilde m_t({\tilde L}_t)|+|m_t({\tilde L}_t+w')-\tilde m_t({\tilde L}_t+w')|\\
    &\leq N^{-1/3+5\fc}+|\tilde m_t({\tilde L}_t+w)-\tilde m_t({\tilde L}_t)|+\OO_\prec\left( \frac{1}{N\Im[w']}\right),
\end{align*}
where $\Im[w]=\eta_0=N^{-1/3-\fc}$.
Thanks to the square root behavior of $\tilde m_t$, close to the spectral edge we have $|\tilde m_t({\tilde L}_t+w)-\tilde m_t({\tilde L}_t)|\lesssim |w|^{1/2}\lesssim N^{-1/3+\fc/2}$.
Therefore it follows that 
\begin{align}\label{e:mtbb}
    m_t({\tilde L}_t+w)-\zeta_t\prec  \frac{N^{5\fc}}{N^{1/3}}.
\end{align}
By plugging \eqref{e:mtbb} into \eqref{e:mtbb00}, and using that $\del_m Q_t(m)$ is a finite polynomial in $m$ with coefficients bounded by $\OO_\prec(1/q)$, we conclude
\begin{align*}
    \del_t {\tilde L}_t=\frac{1}{2}\del_m Q_t(m)|_{m=\zeta_t}
    &=\frac{1}{2}\del_m Q_t(m_t({\tilde L}_t+w))+\OO_\prec\left(\frac{N^{5\fc}}{N^{1/3}q}\right).
\end{align*}
This finishes the proof of Proposition \eqref{p:Ltder}.
\end{proof}

\begin{proof}[Proof of \eqref{e:dFXt}]
Let $w=y+\ri\eta_0$, for the first term on the righthand side of \eqref{e:dFt}, $\sum_{k}D_{ij}G_{kk}=-\sum_k G_{ik}G_{jk}=\del_w G_{ij}$, we can rewrite it as
\begin{align}\begin{split}\label{e:fft}
\sum_{ij}\bE\left[\dot h_{ij}(t)F'(X_t) \del_w G_{ij}(\tilde L_t +w)\right].
\end{split}\end{align}
By using Proposition \ref{p:Ltder}, we can rewrite the second term on the righthand side of \eqref{e:dFt} as
\begin{align}\begin{split}\label{e:fft2}
N\bE\left[F'(X_t)\del_t{{\tilde L}}_t\del_z m_t(z)\right]
&=\frac{N}{2}\bE\left[F'(X_t)(\del_mQ_t)(m_t)\del_z m_t(z)\right]
+\OO_\prec\left(\frac{N^{2/3+5\fc}}{q}\right)\bE[|\del_z m_t(z)|]\\
&=\frac{N}{2}\bE\left[F'(X_t)\del_z(Q_t(m_t(z)))\right]
+\OO_\prec\left(\frac{N^{1+9\fc}}{q}\right),\\
&=\frac{N}{2}\bE\left[F'(X_t)\del_z(Q(e^{-t/2}m_t(z)))\right]
+\OO_\prec\left(\frac{N^{1+9\fc}}{q}\right),
\end{split}\end{align}
where we used Proposition \ref{p:Ltder} in the first equality, and \eqref{e:dermtha} in the second line. 

By plugging \eqref{e:fft} and \eqref{e:fft2} into \eqref{e:dFt}, we get
\begin{align}\begin{split}\label{e:dFt2}
    &\frac{\rd}{\rd t}\bE[F(X_t)]=\OO_\prec\left( \frac{N^{10\fc+1/3}}{q}\right)\\
    &+\Im\int_E^{N^{-2/3+\fc}}\bE\left[F'(X_t)\left(\sum_{ij}\dot h_{ij}(t)\del_w G_{ij}(\tilde L_t+w)+\frac{N}{2}\del_w (Q(e^{-t/2}m_t(\tilde L_t+w)))\right)\right]\rd y\\
    &=\Im\left.\bE\left[F'(X_t)\left(\sum_{ij}\dot h_{ij}(t) G_{ij}(\tilde L_t+w)+\frac{N}{2} Q(e^{-t/2}m_t(\tilde L_t+w))\right)\right]\right|_{y=E}^{y=N^{-2/3+\fc}}+\OO_\prec\left( \frac{N^{10\fc+1/3}}{q}\right).
\end{split}\end{align}
For the first term on the righthand side of \eqref{e:dFt2}, by the cumulant expansion formula, we have
\begin{align}\begin{split}\label{e:ffterm}
&-\sum_{ij}\bE\left[\dot h_{ij}(t)F'(X_t) G_{ij}\right]=\frac{1}{2}\sum_{ij}\bE\left[e^{-t/2}h_{ij}F'(X_t) G_{ij}\right]
-\sum_{ij}\bE\left[\frac{e^{-t}}{2\sqrt{1-e^{-t}}}w_{ij}F'(X_t) G_{ij}\right]\\
&=\sum_{p= 1}^M \frac{e^{-t/2}\cC_{p}}{2Nq^{p-1}}\sum_{ij}\bE[\del_{ij}^p(F'(X_t)G_{ij})]
-\frac{e^{-t}}{2N\sqrt{1-e^{-t}}}\sum_{ij}\bE[\del_{w_{ij}}(F'(X_t)G_{ij})]
+\OO_{\prec}\left(\frac{1}{q^M}\right)\\
&=\sum_{p= 1}^M \frac{e^{-t/2}\cC_{p}}{2Nq^{p-1}}\sum^*_{ij}\bE[\del_{ij}^p(F'(X_t)G_{ij})]
-\frac{e^{-t}}{2N\sqrt{1-e^{-t}}}\sum^*_{ij}\bE[\del_{w_{ij}}(F'(X_t)G_{ij})]
+\OO_{\prec}\left(1\right)\\
&=\sum_{p= 2}^M \frac{e^{-t/2}\cC_{p}}{2Nq^{p-1}}\sum^*_{ij}\bE[\del_{ij}^p(F'(X_t)G_{ij})]
-\frac{e^{-t/2}}{2N}\sum^*_{ij}\bE[\del_{ij}(F'(X_t)G_{ij})]
+\OO_{\prec}\left(\frac{ N^{1/3+5\fc}}{q }\right)\\
&=\sum_{p= 2}^M \frac{e^{-t/2}\cC_{p}}{2Nq^{p-1}}\sum^*_{ij}\bE[\del_{ij}^p(F'(X_t)G_{ij})]
+\OO_\prec \left(\frac{ N^{1/3+5\fc}}{q }\right),
%&=\sum_{p= 1}^M \frac{e^{-(p+1)t/2}\cC_{p}}{2Nq^{p-1}}\sum_{iij}\bE[\del_{ij}^p(F'(X_t)G_{ij}G_{ki})]
%-\frac{e^{-t}}{2N}\sum_{iij}\bE[\del_{ij}(F'(X_t)G_{ij}G_{ki})]
%+\OO\left(\frac{1}{q^M}\right)\\
%&=-\sum_{p=2}^M \frac{e^{-(p+1)t/2}\cC_{p}}{2Nq^{p-1}}\sum_{jk}\bE[\del_{jk}^p(F'(X_t)\del_z G_{jk})]+\OO\left(\frac{1}{q^M}\right)
\end{split}\end{align}
where in the third line, we used that the contribution from terms corresponding to $i=j$ is of oder $\OO_\prec(1)$; in the fourth line, we used 
the relation \eqref{e:dtLt} between $\del_{ij}$ and $\del_{w_{ij}}$.
%that $h_{jk}(t)=e^{-t/2}h_{jk}+\sqrt{1-e^{-t}}w_{jk}$, so $\del_{h_{jk}}=e^{-t/2}\del_{jk}$, and $\del_{w_{jk}}=\sqrt{1-e^{-t}}\del_{jk}$.
The claim \eqref{e:dFXt} follows from plugging \eqref{e:ffterm} into \eqref{e:dFt2}.
\end{proof}

\subsection{Proof of \eqref{e:dFX0}}\label{s:p2}
If we replace $F'(X_t)$ by $P^{r-1}\bar P^r$, the expression on the righthand side of \eqref{e:dFX0} is essentially the same as \eqref{e:Pmoment}, up to some $e^{-t/2}$ factors (In \eqref{e:Pmoment}, the term corresponds to $p=1$ cancels with $\bE[m_N P^{r-1}\bar P^r]$).
We have these $e^{-t/2}$ factors in \eqref{e:dFX0}, because the cumulant expansion formula with respect to $h_{ij}(t)$ is slightly different from the cumulant expansion formula with respect to $h_{ij}$. We record the
 cumulant expansion formula with respect to $h_{ij}(t)$. Take $U=R_{ij\bm m} \del^{\bm\beta}(F'(X_t))$, for any monomial $R_{ij\bm m}$ of the Green's function entries $G(\tilde L_t+w, t)$  as in Definition \ref{d:evaluation} and $\bm\beta=\{\beta\}_{u,v\in ij\bm m}$ with $\beta_{ij}=\beta_{ji}=0$, then
\begin{align}\begin{split}\label{e:hijt}
    &\phantom{{}={}}\frac{1}{N^{r+2}}\sum^*_{ij\bm m}\bE[h_{ij}(t)U]
=\frac{1}{N^{r+2}}\sum^*_{ij\bm m}\bE\left[e^{-t/2}h_{ij}U\right]
+\frac{1}{N^{r+2}}\sum^*_{ij\bm m}\bE\left[\sqrt{1-e^{-t}}w_{ij}U\right]\\
&=\frac{1}{N^{r+2}}\sum^*_{ij\bm m}\sum_{p=1}^M \frac{e^{-t/2}\cC_{p}}{Nq^{p-1}}\bE[\del_{ij}^pU]
+\frac{\sqrt{1-e^{-t}}}{N^{r+3}}\sum^*_{ij\bm m}\bE[\del_{w_{ij}}U]+\OO_\prec\left(\frac{1}{q^M}\right)\\
&=\frac{1}{N^{r+2}}\sum^*_{ij\bm m}\sum_{p\geq 1} \frac{e^{-t/2}\cC_{p}}{Nq^{p-1}}\bE[\del_{ij}^pU]
+\frac{\sqrt{1-e^{-t}}}{N^{r+3}}\sum^*_{ij\bm m}\bE[e^{t/2}\sqrt{1-e^{-t}}\del_{ij}U]+\OO_\prec\left(\frac{N^{5\fc}}{q N^{2/3}}\right)\\
&=\frac{e^{t/2}}{N^{r+3}}\sum^*_{ij\bm m}\bE[\del_{ij}U]+\frac{1}{N^{r+2}}\sum^*_{ij\bm m}\sum_{p\geq 2} \frac{e^{-t/2}\cC_{p}}{Nq^{p-1}}\bE[\del_{ij}^pU]+\OO_\prec\left(\frac{N^{5\fc}}{q N^{2/3}}\right),
%&=\sum_{p\geq 1} \frac{e^{-(p+1)t/2}\cC_{p}}{Nq^{p-1}}\sum_{ij}\bE[\del_{ij}^p U]
%+\frac{1-e^{-t}}{N}\sum_{ij}\bE[\del_{ij}U]\\
%&=\sum_{p\geq 1} \frac{e^{-(p+1)t/2}\cC_{p}}{Nq^{p-1}}\sum_{ij}\bE[\del_{ij}^p U]
%+\frac{1-e^{-t}}{N}\sum_{ij}\bE[\del_{ij}U]\\
%&=\frac{1}{N}\bE[\del_{ij}U]+\sum_{p\geq 2} \frac{e^{-(p+1)t/2}\cC_{p}}{Nq^{p-1}}\sum_{ij}\bE[\del_{ij}^pU]
\end{split}\end{align}
where we used \eqref{e:dtLt} to replace $\del_{w_{ij}}$ in the second line.

Similarly to \eqref{e:ordero}, all the terms we will get  in the expansion are in the form
%\begin{align}\label{e:ordero2}
%    \frac{1}{q^\fo}\frac{1}{N^r}\sum_{\bm m}\bE\left[\prod_{x\in \bm m}G_{xx}^{\al_{x}} \prod_{x,y\in \bm m}\del_{xy}^{\beta_{xy}}F'(X_t)\right]\prec \frac{1}{q^\fo},
%\end{align}
%where $\bm m$ is an index set with size $|\bm m|=r$, $\al_{x},\beta_{xy}\geq 1$ for $x,y\in \bm m$, and  $\fo\geq 0$ is the order. For terms with order at least $M$, we will trivially bound them by $\OO_\prec(1/q^M)$.
\begin{align}\label{e:ordero2}
    \frac{1}{q^\fo}\times \frac{1}{N^r}\sum_{\bm m}^*\bE\left[R_{\bm m} \left(\prod_{e\in E(\cF)}\del_{\al_e\beta_e}^{p_{e}-s_e}
    \right)(F'(X_t))\right],
\end{align}
where $\cF$ is a weighted forest  with vertex set $V(\cF)=\bm m=\{m_1, m_2,\cdots, m_r\}$, $R$ is a monomial as in Definition \ref{d:evaluation},   $\bm p=\{p_{e}\}_{e\in E(\cF)}$ are nonnegative integers, and  $\fo\geq 0$ is the order parameter. Since the second factor in \eqref{e:ordero} can be trivially bounded by $\OO_\prec(1)$, the whole expression can be bounded by $\OO_\prec(1/q^\fo)$.  For terms with order at least $M$, we will trivially bound them by $\OO_\prec(1/q^M)$.

\begin{proof}[Proof of \eqref{e:dFX0}]
We follow the three step strategy as in the proof of Proposition \ref{p:DSE}. 

\noindent\textit{Step $1$ (eliminate off-diagonal Green's function terms) }
The first term on the righthand side of \eqref{e:dFX0} is in the following form
\begin{align}\label{e:dFX}
\frac{1}{N}\sum_{ij}\sum_{p=2}^{ M}\sum_{s= 0}^p\frac{e^{-t/2}\cC_{p+1}}{Nq^{p-1}}{p\choose s}\bE[(\del_{ij}^{s}G_{ij}(\tilde L_t+w;t))\del_{ij}^{p-s}(F'(X_t))].
\end{align}
Thanks to Proposition \ref{p:Ltder2}, we can replace $\del_{ij}^{s}G_{ij}(\tilde L_t+w;t)=e^{-st/2}D_{ij}^{s}G_{ij}(\tilde L_t+w;t)+\OO_\prec(N^{-2/3+4\fc})$. The error term is bounded by  
\begin{align*}
\frac{1}{N}\sum_{ij}\sum_{p=2}^{ M}\sum_{s= 0}^p\frac{e^{-t/2}\cC_{p+1}}{Nq^{p-1}}{p\choose s}\bE[N^{-2/3+4\fc}|\del_{ij}^{p-s}(F'(X_t))|]\prec 
\frac{1}{N^2 }\sum_{ij}\sum_{p=2}^{ M}\sum_{s= 0}^p\frac{N^{-2/3+4\fc}}{q^{p-1}}
\lesssim \frac{N^{4\fc}}{N^{2/3}q},
\end{align*}
where we used that $|\del_{ij}^{p-s}(F'(X_t))|\prec 1$, from \eqref{e:impp3}. For the term in \eqref{e:dFX} with $i=j$, we can similarly bound them as
\begin{align*}
\frac{1}{N}\sum_{i}\sum_{p=2}^{ M}\sum_{s= 0}^p\frac{e^{-(s+1)t/2}\cC_{p+1}}{Nq^{p-1}}{p\choose s}\bE[(D_{ii}^{s}G_{ii})\del_{ii}^{p-s}(F'(X_t))]\prec \frac{1}{N^2}\sum_{i}\sum_{p=2}^{ M}\sum_{s= 0}^p\frac{1}{q^{p-1}}\lesssim \frac{1}{Nq}.
\end{align*}
Therefore we can further restrict the summation in \eqref{e:dFX} to $i\neq j$. 
\begin{align}\label{e:dFXa}
\frac{1}{N}\sum^*_{ij}\sum_{p=2}^{ M}\sum_{s= 0}^p\frac{e^{-(s+1)t/2}\cC_{p+1}}{Nq^{p-1}}{p\choose s}\bE[(D_{ij}^{s}G_{ij})\del_{ij}^{p-s}(F'(X_t))].
\end{align}

The derivative $D_{ij}^{s}G_{ij}$ is a sum of monomials in the form $G_{ii}^a G_{jj}^b G_{ij}^c$, 
\begin{align*}
D_{ij}^s G_{ij}=-\bm1(s \text{ is odd})s! G_{ii}^{\frac{s+1}{2}}G_{jj}^{\frac{s+1}{2}}+\{\text{terms with off-diagonal entries}\}.
\end{align*}
If the monomial contains at least two off-diagonal terms, i.e. $c\geq 2$, then it is bounded by $|G_{ij}|^2$ and \eqref{e:derofGG} gives
\begin{align*}
    \frac{1}{N^2q^{p-1}}\sum_{ij}\bE[|G_{ij}|^2 |\del_{ij}^{p-s}(F'(X_t))|]
    \prec
    \frac{1}{N^{2}q^{p-1}}\sum_{ij}\bE[|G_{ij}|^2]
    \prec 
    \frac{N^{4\fc}}{N^{2/3}q}.
\end{align*}
Analogous to Propostion \ref{p:one-off}, terms with exactly one off-diagonal term are negaligible. We have the following 
\begin{proposition}\label{p:one-off2}
Adopt the assumptions of Proposition \ref{p:comp}. Given a weighted forest $\cF$ with vertex set $V(\cF)=ij\bm m$, where $\bm m=\{m_1, m_2,\cdots, m_r\}$. Then for any monomial $R_{ij\bm m}$ of Green's function entries $G(\tilde L_t+w,t)$ as in Definition \ref{d:evaluation} and nonnegative integers  $\bm p=\{p_{e}\}_{e\in E(\cF)}$ such that $p_e\geq s_e$, we have
\begin{align}\label{e:one-off2}
    \frac{1}{N^{r+2}}\sum_{ij\bm m}^*\bE\left[ G_{ij}R_{ij\bm m}\left(\prod_{e\in E(\cF)}\del_{\al_e\beta_e}^{p_{e}-s_e}\right)(F'(X_t)) \right]\prec \frac{N^{7\fc}}{N^{2/3}}.
\end{align}
\end{proposition}

As a consequence of Proposition \ref{p:one-off2}, we have that the  terms in $D_{ij}^s G_{ij}$ with exactly one off-diagonal term, i.e. $c=1$, are bounded by
\begin{align*}
    \frac{1}{N^2q^{p-1}}\sum_{ij}\bE[G_{ij}G^a_{ii}G^b_{jj}\del_{ij}^{p-s}(F'(X_t))]\prec \frac{N^{7\fc}}{N^{2/3}q},
\end{align*}
where $p\geq 2$.

Combining the discussion above, the leading term in \eqref{e:dFXa} comes from the monomials with only diagonal Green's functions entries, 
\begin{align}\label{e:dFX2}
-\frac{1}{N}\sum_{p=2}^{ M}\sum_{s\text{ odd}}\sum^*_{ij}\frac{e^{-(s+1)t/2}\cC_{p+1}s!}{Nq^{p-1}}{p\choose s}\bE[G_{ii}^{\frac{s+1}{2}}G_{jj}^{\frac{s+1}{2}}\del_{ij}^{p-s}(F'(X_t))]+\OO_\prec\left(\frac{N^{7\fc}}{N^{2/3}q}\right).
\end{align}

\noindent\textit{Step $2$ (replace diagonal Green's function entries by $m_t$) }
Analogous to Proposition \ref{p:diag}, we can use the following proposition to replace diagonal Green's function entries
by $m_t$.
\begin{proposition}\label{p:diag2}
Adopt the assumptions of Proposition  \ref{p:comp}. Given a weighted forest $\cF$ with vertex set $V(\cF)=i\bm m$, where $\bm m=\{m_1, m_2,\cdots, m_r\}$. Then for any monomial $R_{\bm m}$ as in Definition \ref{d:evaluation} with no off-diagonal Green's function entries, i.e. $\chi(R_{\bm m})=0$, and  integers  $\bm p=\{p_{e}\}_{e\in E(\cF)}$ with $p_e-s_e\geq 0$, we have
\begin{align}\begin{split}\label{e:diag2}
    &\frac{1}{N^{r+1}}\sum^*_{i\bm m}\bE\left[ G^\al_{ii}R_{\bm m}V \right]=\frac{1}{N^{r+1}}\sum^*_{i\bm m}\bE\left[ G^{\al-1}_{ii} m_N R_{\bm m}V \right] + \Omega_1+\Omega_2+\OO_{\prec}(N^{6\fc-2/3}),\\ 
    &V=\left(\prod_{e\in E(\cF)}\del_{\al_e \beta_e}^{p_e-s_e}\right)(F'(X_t)),
\end{split}\end{align}
where the $\Omega_1, \Omega_2$ are given by
\begin{align}\begin{split}\label{e:hod2}
&\Omega_1=-\frac{1}{N^{r+2}}\sum_{p=2}^M\sum_{s \text{ odd}} \frac{e^{-(s+1)t/2}\cC_{p+1}s!}{Nq^{p-1}}{p\choose s}
   \sum^*_{ij k \bm m} \bE[G^\al_{ii}G_{j j}^{\frac{s+1}{2}}G_{kk}^{\frac{s+1}{2}} R_{\bm m}\del^{p-s}_{j k} V)],\\
& \Omega_2=  \frac{1}{N^{r+1}}\sum_{p=2}^M\sum_{s \text{odd}} \frac{e^{-(s+1)t/2} \cC_{p+1}s!}{Nq^{p-1}}{p\choose s} {\frac{s-1}{2}+\al-1\choose \al-1}
    \sum^*_{i k \bm m}\bE[m_NG^{\al+\frac{s-1}{2}}_{ii}G^{\frac{s+1}{2}}_{kk}R_{\bm m}\del^{p-s}_{ik} V)].\end{split}\end{align}
Comparing these terms $\Omega_1, \Omega_2$ with \eqref{e:diag2}, since $p\geq 2$ they are of order at least $1$ (recall from \eqref{e:ordero2}).
\end{proposition}

\begin{remark}
These terms $\Omega_1, \Omega_2$ in \eqref{e:hod2} are in the same form as in \eqref{e:diag2}. With given $p$, the term in $\Omega_1$  is associated with a weighted forest $\cF_1$, which is from $\cF$ by adding vertices $j,k$ and an edge $\{j,k\}$ with weight $s$. In total $\Omega_1$ has $r+2$ vertices.  $\Omega_2$ is associated with a weighted forest $\cF_2$, which is from $\cF$ by adding one vertex $k$ and an edge $\{i,k\}$ with weight $s$. In total $\Omega_2$ has $r+1$ vertices.
For given $s$, both $\Omega_1, \Omega_2$ have an extra factor $e^{-(s+1)t/2}$ and derivative $\del^{p-s}$. Moreover, the total number of diagonal Green's function entries increases by $s+1$. 
\end{remark}

By repeatedly using Proposition \ref{p:diag2}, we can  replace the product $G_{ii}^{(s+1)/2}G_{jj}^{(s+1)/2}$ in \eqref{e:dFX2} by $m_t^{s+1}$, to get the leading terms
\begin{align}\label{e:nextt2}\begin{split}
    &\phantom{{}={}}-\frac{1}{N}\sum_{p=2}^{ M}\sum_{s\text{ odd}}\sum^*_{ij} \frac{e^{-(s+1)t/2}\cC_{p+1}s!}{Nq^{p-1}}{p\choose s}\bE\left[m_t^{s+1}\del_{ij}^{p-s}F'(X_t)\right],
\end{split}\end{align}
with higher order terms (which have at least one more copies of $1/q$) as linear combination of terms (with bounded coefficients) in the form 
\begin{align}\label{e:final2}
  \frac{1}{N^{r}}\sum_{ 2\leq p_{e_1},\cdots, p_{e_{|E(\cF)|}}\leq M}
  \sum^*_{\bm m}\bE\left[R_{\bm m}
  \left(\prod_{e\in E(\cF)}\frac{e^{-(s_e+1)t/2}\cC_{p_{e+1}}s_e!}{q^{p_e-1}}{p_e\choose s_e}  \del_{\al_e \beta_e}^{p_e-s_e}\right)(F'(X_t))\right],
\end{align}
where  $\cF$ is a weighted forest as in Definition \ref{d:forest} with vertex set $\bm m=\{m_1,m_2,\cdots, m_r\}$; the monomial $R_{\bm m}$ has no off-diagonal entries, i.e. $\chi(R_{\bm m})=0$ and $\deg(R_{\bm m})=\sum_{e\in E(\cF)}(s_e+1)$.  

We can repeat Step 2 for these higher order terms \eqref{e:final2}. If a term has order bigger than $M$, we can trivially bound it by $\OO_\prec(1/q^M)\prec N^{8\fc}/qN^{2/3}$ as in \eqref{e:ordero2}.  
The final expression is a linear combination (with bounded coefficients) of terms in the form: 
\begin{align}\begin{split}\label{e:final3}
   &\bE[m_t^{2\ell} L^t_\cF(P^{r-1}\bar P^{r})],\\
   &L^t_\cF= \sum_{2\leq p_{e_1},\cdots, p_{e_{|E(\cF)|}}\leq M}\frac{1}{N^{|V(\cF)|}}\sum_{ x_1,\cdots, x_{|V(\cF)|}}^*\prod_{e\in E(\cF)}\frac{e^{-(s_e+1)t/2}\cC_{p_{e+1}} s_e!}{Nq^{p_e-1}}{p_e\choose s_e}  \del_{\al_e\beta_e}^{p_e-s_e}
\end{split}\end{align}
where  $\cF$ is a forest as in Defintion \ref{d:forest}, $x_1, x_2,\cdots, x_{|V(\cF)|}$ enumerate the vertices of $\cF$. Moreover, all the weights $s_e$ are odd positive integers, and the total weights satisfies $\sum_e (s_e+1)=2\ell$. The above discussion leads to the following claim.
\begin{claim}\label{c:step22}
Under the assumptions of Proposition \ref{p:diag2}, \eqref{e:dFXa} is a finite sum of terms in the form, with an error $\OO_\prec(N^{8\fc}/qN^{2/3})$:
\begin{align}\label{e:finalt1}
    \bE[m_t^{2\ell} L^t_\cF F'(X_t)],
\end{align}
where $L^t_\cF$ is as defined in \eqref{e:final3}.
\end{claim}

\noindent\textit{Step $3$ (rewrite differential operators as an expectation) }
Finally use the cumulant expansion, the same as in \eqref{e:ediff}, we have
\begin{align}\begin{split}\label{e:ediffc}
    &\phantom{{}={}}\bE\left[e^{-(s+1)t/2}\left(h_{ij}^{s+1}-\frac{\bm1(s=1)}{N}\right)m_N^{2\ell}F'(X_t)\right]\\
    &=\sum_{p=2}^M\frac{e^{-(s+1)t/2}\cC_{p+1}s!}{Nq^{p-1}}{p\choose s} \bE[m_N^{2\ell}\del_{ij}^{p-s}(F'(X_t))]+\OO_{\prec}\left(\frac{1}{q^M}+\frac{N^{4\fc}}{N^{5/3} q}\right).
\end{split}\end{align}
By repeatedly using the relation \eqref{e:ediffc}, and in the same way as in Claim \ref{c:step3}, we can rewrite \eqref{e:finalt1} in the following form
\begin{align}\label{e:final2t}\begin{split}
    &\phantom{{}={}}\bE[m_t^{2\ell}{ L}^t_\cF(F'(X_t))]\\
    &=\prod_{e\in E(\cF)} \bE\left[\sum^*_{  x_1,\cdots, x_{|V(\cF)|}}\left(\frac{1}{N^{\theta(\cF)}}  \prod_{e\in E(\cF)} \left(h_{\al_e\beta_e}^{s_e+1}-\frac{\bm1_{s_e=1}}{N}\right)\right) e^{-\ell t}m_t^{2\ell}F'(X_t)\right]+\OO_{\prec}\left(\frac{N^{4\fc}}{qN^{2/3}}\right)\\
    &=\prod_{e\in E(\cF)} \bE\left[w(\cF) e^{-\ell t/2}m_t^{2\ell }F'(X_t)\right]+\OO_{\prec}\left(\frac{N^{4\fc}}{qN^{2/3}}\right).
\end{split}\end{align}

So far every estimate is parallel to those from Proposition \ref{p:one-off}, except for the extra factor $e^{-\ell t/2}$ in \eqref{e:final2t}. 
Thanks to Claim \ref{c:step22} and \eqref{e:final2t}, the first term on the righthand side of \eqref{e:dFX0} is in the form 
\begin{align}\label{e:amat}
    -\bE[(e^{-t}a_2m_t^2+e^{-2t} a_4 m_t^4+\cdots+e^{-\ell t}a_{2L}m_t^{2L})F'(X_t)]+\OO_\prec\left(\frac{N^{8\fc}}{qN^{2/3}}\right),
\end{align}
where
$a_{2\ell}$ is a sum of terms in the form $w(\cF)$ as in \eqref{e:tt}, where  $\cF$ is a forest as in Defintion \ref{d:forest}. Moreover, all the weights $s_e$ are odd positive integers with $\sum_e (s_e+1)=2\ell$.   
The expression \eqref{e:amat} is precisely the definition of the polynomial $Q(e^{-t}m_t) F'(X_t)$. Thus 
the term \eqref{e:amat} cancels with $\bE[Q(e^{-t}m_t) F'(X_t)]$ in \eqref{e:dFX0}, and we conclude Proposition \ref{p:DSEt}.

\end{proof}

\subsection{Proof of Propositions from Sections \ref{s:p1} and \ref{s:p2} }

\begin{proof}[Proof of Proposition \ref{p:Ltder2}]
    Let $w'=y+N^{-2/3+\fc}\ri$, then it is easy to see that $w'\in \cD$ (recall from \ref{e:recalldefD}), and Proposition \ref{p:mt} gives
    \begin{align*}
        \Im[m_t({\tilde L}_t+w')]\leq \Im[\tilde m_t({\tilde L}_t+w')]+\OO_\prec\left(\frac{1}
        {N\Im[w']}\right)
        \lesssim \sqrt{|w'|}+\OO_\prec\left(\frac{1}
        {N\Im[w']}\right)\lesssim \frac{N^{\fc/2}}{N^{1/3}},
    \end{align*}
    where we used that $\tilde m_t$ has square root behavior.
    The derivative of $m_t$ satisfies
    \begin{align}\label{e:dermt}
        |\del_z \Im[m_t(z)]|\leq |\del_z m_t(z)|\leq \frac{\Im[m_t(z)]}{\Im[z]},
    \end{align}
    which gives that $\Im[m_t(E+\ri \eta/M)\leq M\Im[m_t(E+\ri \eta)]$ for any $M\geq 1$. In particular, we have
    $\Im[m_t({\tilde L}_t+w)]\leq (\Im[w']/\Im[w])\Im[m_t({\tilde L}_t+w')]\leq N^{2\fc}N^{-1/3+\fc/2}\leq N^{-1/3+3\fc}$. Using \eqref{e:dermt} again, we have 
    $|\del_z m_t({\tilde L}_t+w)|\leq N^{1/3+4\fc}$. 
For the derivative of the Green's function, Ward identity \eqref{e:WdI} implies 
\begin{align}\begin{split}\label{e:GG}
|\del_w G_{ab}(\tilde L_t+w)|
&\leq \sum_{i=1}^N |G_{ai}(\tilde L_t+w)G_{bi}(\tilde L_t+w)|
\leq \frac{1}{2}\sum_{i=1}^N (|G_{ai}(\tilde L_t+w)|^2+|G_{bi}(\tilde L_t+w)|^2)\\
&\prec \frac{\Im[m_N(\tilde L_t+w)]}{\eta}\leq N^{1/3+4\fc}.
\end{split}\end{align}

Since $\tilde L_t$ is the spectral edge of $\tilde \rho_t$, which is characterized by $1+z\tilde m_t(z)+\tilde m_t^2(z)+Q(e^{-t/2}\tilde m_t(z))=0$. The same as in Proposition \ref{p:minfty}, $\tilde L_t$  depends smoothly on the coefficients of $Q$. In particular, its derivatives with respect to $a_2, a_4,\cdots, a_{2L}$ are bounded. Thus the bounds on $|\del^{\bm\beta}{\tilde L_t}|$ follow from Proposition \ref{p:a2kstr}.  The estimates in \eqref{e:impp2} can be proven the same way as \eqref{e:impp}, and using that $\Im[m_t(\tilde L_t+w)]/N\eta\prec N^{-1/3+3\fc}$. 

For \eqref{e:impp3}, we can rewrite the derivative $\del^{\bm \beta}m_t$ as
\begin{align*}
\del^{\bm\beta}m_t=\frac{1}{N}\sum_{i=1}^N e^{|\bm\beta|t/2}D^{\bm\beta}G_{ii} +\OO_\prec \left(N^{-2/3+4\fc}\right).
\end{align*}
$D^{\bm\beta}G_{ii}$ is a monomial of Green's function entries, and each contains at least two off-diagonal entries. We can bound the sum using the Ward identity \eqref{e:WdI} as in \eqref{e:GG},
\begin{align*}
\frac{1}{N}\sum_{i=1}^N e^{|\bm\beta|t/2}D^{\bm\beta}G_{ii} \prec \frac{N^{1/3+4\fc}}{N} =\frac{N^{4\fc}}{N^{2/3}}. 
\end{align*}
The second and third relation in \eqref{e:impp3} follows from 
\begin{align*}
|\del^{\bm\beta}X_t|
=\left| \Im\left[N\int_{E}^{N^{-2/3+\fc}}\del^{\bm\beta} m_t( {\tilde L}_t +w) \rd y\right]\right|\lesssim N^{-1/3+5\fc},
\end{align*}
and the fact that $F$ has bounded derivatives. 

For either $U=R_{ij\bm m}$ or $U=F'(X_t)$, $U=U(H(t), \tilde L_t)$ is a function of both $H(t)$ and $\tilde L_t$. Since $\tilde L_t$ depends only on $H$ but not on $W$, the derivatives of $\del_{ij}$ and $\del_{w_{ij}}$ are related by the following relation
\begin{align}\label{e:relatewh}
    \del_{w_{ij}}U= e^{t/2}\sqrt{1-e^{-t}}(\del_{ij}U-(\del_{ij}\tilde L_t) \del_{\tilde L_t}U).  
\end{align}
Using the relation \eqref{e:relatewh}, we can rewrite the lefthand side of \eqref{e:dtLt}
\begin{align}\begin{split}\label{e:dtLt2}
&\phantom{{}={}}\frac{1}{N}\sum_{ij:ij\bm m\atop \text{distinct}}\bE\left[ \del_{w_{ij}}(R_{ij\bm m}\del^{\bm\beta}(F'(X_t)))\right]= \frac{e^{t/2}\sqrt{1-e^{-t}}}{N}\sum_{ij:ij\bm m\atop \text{distinct}}\bE\left[ \del_{ij}(R_{ij\bm m}\del^{\bm\beta}(F'(X_t)))\right]\\
&+\frac{e^{t/2}\sqrt{1-e^{-t}}}{N}\sum_{ij:ij\bm m\atop \text{distinct}}\bE\left[ (\del_{ij}\tilde L_t)(\del_{\tilde L_t}R_{ij\bm m})\del^{\bm\beta}(F'(X_t))+R_{ij\bm m}\del^{\bm\beta}((\del_{ij}\tilde L_t) \del_{\tilde L_t}F'(X_t))\right].
\end{split}\end{align}
Using \eqref{e:dermtha},  \eqref{e:impp2}, \eqref{e:impp3} for any $\bm\beta'$ with $\beta'_{ij}=\beta'_{ji}=0$, we have
\begin{align*}
|\del_{ij}\tilde L_t|\prec \frac{|h_{ij}|}{N},
\quad 
|\del^{\bm\beta'}\del_{ij}\tilde L_t|\prec \frac{|h_{ij}|}{N},
\quad
|R_{ij\bm m}|\prec 1,
\quad
|\del_{\tilde L_t}R_{ij\bm m}|\prec N^{1/3+4\fc},
\quad
 \del^{\bm \beta'}X_t\prec N^{1/3+5\fc}.
\end{align*}
Moreover we also have 
\begin{align*}
 \del^{\bm \beta'}\del_{\tilde L_t}X_t=\left| \Im\left[N\int_{E}^{N^{-2/3+\fc}}\del^{\bm\beta}\del_{\tilde L_t} m_t( {\tilde L}_t +w) \rd y\right]\right|\lesssim N^{1/3+5\fc}.
 \end{align*}
Then we can bound the second term on the righthand side of \eqref{e:dtLt2} as
\begin{align*}
\frac{1}{N}\sum_{ij: ij \bm m\atop \text{distinct}}N^{1/3+5\fc}\frac{\bE[|h_{ij}|]}{N}
\prec  \frac{N^{1/3+5\fc}}{q}.
\end{align*}
The claim \eqref{e:dtLt} follows.
\end{proof}

\begin{proof}[Proof of Proposition \ref{p:one-off2}]
If $\sum_e (p_e-s_e)\geq 1$, then \eqref{e:impp3} implies that
\begin{align*}
    \left|\left(\prod_{e\in E(\cF)}\del_{\al_e\beta_e}^{p_{e}-s_e}\right)(F'(X_t))\right|\prec N^{-1/3+5\fc},
\end{align*}
and
\begin{align}\begin{split}\label{e:atleastoned}
  \frac{1}{N^{r+2}}\sum_{ij\bm m}^*\bE\left[ G_{ij}R_{ij\bm m}\left(\prod_{e\in E(\cF)}\del_{\al_e\beta_e}^{p_{e}-s_e}\right)(F'(X_t)) \right]\prec \frac{1}{N^{r+2}}\sum_{ij\bm m}\bE[|G_{ij}|N^{-1/3+5\fc}]
   \prec \frac{N^{7\fc}}{N^{2/3}},
\end{split}\end{align}
where in the last inequality we used \eqref{e:derofGG}.
If $\sum_e (p_e-s_e)=0$, then we use the identities $G_{ij}=\sum_{k\neq i}G_{ii}h_{ik}G_{kj}^{(i)}$, and $G_{kj}^{(i)}=G_{kj}-G_{ki}G_{ji}/G_{ii}$.  
We denote
\begin{align*}
    U=R_{ij\bm m}F'(X_t).
\end{align*}
Then by the cumulant expansion \eqref{e:hijt},  we have
\begin{align}\begin{split}\label{e:redu}
&\phantom{{}={}}\frac{1}{N^{r+2}}\sum^*_{ij\bm m}\bE\left[G_{ij}U
\right]
=\frac{1}{N^{r+2}}\sum^*_{ij\bm m}\bE\left[\sum_{k\neq i}h_{ik}(t)G_{kj}^{(i)}G_{ii}U\right]\\
&=\frac{1}{N^{r+2}}\sum^*_{ij\bm m}\sum_{p=1}^M\sum_{k:k\neq i} \frac{e^{-t/2}\cC_{p}}{Nq^{p-1}}\bE[\del_{ik}^p(G_{kj}^{(i)}G_{ii}U)]
+\frac{\sqrt{1-e^{-t}}}{N^{r+3}}\sum^*_{ij\bm m}\sum_{k:k\neq i}\bE[\del_{w_{ik}}(G_{kj}^{(i)}G_{ii}U)]+\OO_\prec\left(\frac{1}{q^M}\right)\\
&=\frac{1}{N^{r+2}}\sum^*_{ijk\bm m}\sum_{p=1}^M\frac{e^{-t/2}\cC_{p}}{Nq^{p-1}}\bE[\del_{ik}^p(G_{kj}^{(i)}G_{ii}U)]
+\frac{\sqrt{1-e^{-t}}}{N^{r+3}}\sum^*_{ijk\bm m}\bE[\del_{w_{ik}}(G_{kj}^{(i)}G_{ii}U)]+\OO_\prec\left(\frac{1}{N}+\frac{1}{q^M}\right)\\
&=\frac{e^{t/2}}{N^{r+3}}\sum^*_{ijk\bm m}\bE\left[\del_{ik}(G_{kj}^{(i)} G_{ii}U) \right]
+
\sum_{p=2}^M\frac{e^{-t/2}\cC_{p+1}}{N^{r+3}q^{p-1}}\sum^*_{ijk\bm m}\bE\left[\del_{ik}^p(G_{kj}^{(i)} G_{ii} U) \right]+\OO_{\prec}\left(\frac{N^{5\fc}}{N^{2/3}}+\frac{1}{q^M}\right).
\end{split}\end{align} 
where to get the third line, we used that the summation for terms with $k\in j\bm m$ is bounded by $1/N$. Then \eqref{e:redu} can be analyzed in the same way as for \eqref{e:reduce}, by using $|F'(X_t)|\lesssim 1$ and  $|\del^{\bm\beta}F'(X_t)|\prec N^{-1/3+4\fc}$ from \eqref{e:impp3}. This leads to the claim  \eqref{e:one-off2}.

\end{proof}

\begin{proof}[Poof of Proposition \ref{p:diag2}] 
The proof is similar to that of Proposition \ref{p:diag}. We will use \eqref{e:diag} to replace a copy of $G_{ii}$ to $m_t$. 
Denote \begin{align*}
U:=R_{i\bm m} V,
\quad V:=\left(\prod_{e\in E(\cF)}\del_{\al_e\beta_e}^{p_{e}-s_e}\right)(F'(X_t)).
\end{align*}
Then we have exactly the same expression as in \eqref{e:diffc00}. For the second term on the righthand side of \eqref{e:diffc00}, we can rewrite it as
\begin{align}\begin{split}\label{e:diffc1}
    &\phantom{{}={}}\frac{1}{N^{r+2}}\sum^*_{ij\bm m}\bE[(G_{ii} (HG)_{jj} - m_t (HG)_{ii})G_{ii}^{\al-1}U]\\
   &= \frac{1}{N^{r+2}}\sum^*_{ij\bm m}\sum_{p=1}^M\sum_k \frac{e^{-t/2}\cC_p}{Nq^{p-1}}
    \bE[\del^p_{jk}(G^\al_{ii}G_{j k} R_{i\bm m}V)-\del^p_{ik}(G_{jj}G_{ik}G_{ii}^{\al-1}R_{i\bm m}V)]\\
          &+\frac{\sqrt{1-e^{-t}}}{N^{r+3}}\sum^*_{ij\bm m}\sum_k
    \bE[\del_{w_{jk}}(G^\al_{ii}G_{j k} R_{i\bm m}V)-\del_{w_{ik}}(G_{jj}G_{ik}G_{ii}^{\al-1}R_{i\bm m}V)]+\OO\left(\frac{1}{q^M}\right)\\
    &= \frac{1}{N^{r+2}}\sum^*_{ijk\bm m}\sum_{p=1}^M \frac{e^{-t/2}\cC_p}{Nq^{p-1}}
    \bE[\del^p_{jk}(G^\al_{ii}G_{j k} R_{i\bm m}V)-\del^p_{ik}(G_{jj}G_{ik}G_{ii}^{\al-1}R_{i\bm m}V)]\\
          &+\frac{\sqrt{1-e^{-t}}}{N^{r+3}}\sum^*_{ijk\bm m}
    \bE[\del_{w_{jk}}(G^\al_{ii}G_{j k} R_{i\bm m}V)-\del_{w_{ik}}(G_{jj}G_{ik}G_{ii}^{\al-1}R_{i\bm m}V)]+\OO\left(\frac{1}{N}+\frac{1}{q^M}\right),
\end{split}\end{align}
where we used that the summation for terms with $k\in ij\bm m$ is bounded by $1/N$.
Then we can replace the derivatives $\del_{w_{jk}}$ and $\del_{w_{ik}}$ by $\del_{jk}$ and $\del_{ik}$ using \eqref{e:dtLt}, and get
\begin{align}\begin{split}\label{e:diffc1}
    &\phantom{{}={}}\frac{1}{N^{r+2}}\sum^*_{ij\bm m}\bE[(G_{ii} (HG)_{jj} - m_t (HG)_{ii})G_{ii}^{\al-1}U]\\
   &= \frac{1}{N^{r+2}}\sum^*_{ij k \bm m}\frac{e^{t/2}}{N}
    \bE[\del_{jk}(G^\al_{ii}G_{j k} R_{i\bm m}V)-\del_{ik}(G_{jj}G_{ik}G_{ii}^{\al-1}R_{i\bm m}V)]+\OO_{\prec}\left(\frac{1}{q^M}\right)\\
          &+\frac{1}{N^{r+2}}\sum^*_{ij k \bm m}\sum_{p=2}^M \sum_{s} \frac{e^{-(s+1)t/2}\cC_{p+1}}{Nq^{p-1}}{p\choose s}
    \bE[\del^{s}_{jk}(G^\al_{ii}G_{j k} R_{i\bm m})\del^{p-s}_{j k} V)-\del^{s}_{ik}(G_{jj}G_{ik}G_{ii}^{\al-1}R_{i\bm m})\del^{p-s}_{ik} V)].
\end{split}\end{align}
Then \eqref{e:diffc1} can be analyzed in the same way as for \eqref{e:diffc0}, by using $|F'(X_t)|\lesssim 1$ and  $|\del^{\bm\beta}F'(X_t)|\prec N^{-1/3+5\fc}$ from \eqref{e:impp3}. This leads to the claim  \eqref{e:hod2}.
\end{proof}

\bibliography{References.bib}{}
\bibliographystyle{abbrv}

\end{document}